\documentclass[fleqn,trsc,nonblindrev]{informs3noheader}

\OneAndAHalfSpacedXI 
\usepackage{endnotes}
\usepackage{booktabs}
\usepackage{etoolbox}  
\usepackage{array}      

\newcommand{\tablefontsize}{10}      
\newcommand{\tablebaselineskip}{8}

\AtBeginEnvironment{tabular}{%
  \fontsize{\tablefontsize}{\tablebaselineskip}\selectfont
  \renewcommand{\arraystretch}{\tablearraystretch}%
}
\AtBeginEnvironment{tabularx}{%
  \fontsize{\tablefontsize}{\tablebaselineskip}\selectfont
  \renewcommand{\arraystretch}{\tablearraystretch}%
}
\AtBeginEnvironment{longtable}{%
  \fontsize{\tablefontsize}{\tablebaselineskip}\selectfont
  \renewcommand{\arraystretch}{\tablearraystretch}%
}

\usepackage[linesnumbered,ruled,vlined]{algorithm2e}
\let\footnote=\endnote
\let\enotesize=\normalsize
\def\notesname{Endnotes}%
\def\makeenmark{$^{\theenmark}$}
\def\enoteformat{\rightskip0pt\leftskip0pt\parindent=1.75em
	\leavevmode\llap{\theenmark.\enskip}}

\usepackage{graphicx}
\usepackage{multirow}
\usepackage{hhline}
\usepackage[short,nocomma]{optidef}
\usepackage{pifont}
\usepackage[small, margin=1cm]{caption}
\usepackage{bm}
\usepackage{appendix}
\usepackage{color}
\usepackage{caption}
\usepackage{url}
\usepackage{subcaption}
\usepackage{graphicx}
\usepackage{booktabs}
\usepackage{filecontents}
\usepackage{tikz}
\usepackage{lineno}
\usepackage{pgfplots}
\usepackage{csquotes}
\usepackage{hyperref}
\hypersetup{
    colorlinks=true, linkcolor=blue,
    citecolor=blue,urlcolor=blue, 
}
\definecolor{strcolor}{rgb}{0.6, 0.2, 0.6}
\definecolor{commentcolor}{rgb}{0.3125, 0.5, 0.3125}
\definecolor{keycol}{rgb}{0, 0, 1}
\usepackage[utf8]{inputenc}
\usepackage{booktabs}    
\usepackage{multirow}    
\usepackage{caption}
\usepackage{adjustbox}   
\usepackage{siunitx}     
\usepackage{array}
\usepackage{threeparttable}
\usepackage{adjustbox}
\usepackage{multirow}

\usepackage{bbm}

\usepackage{listings}
\usepackage{url}

\newcommand {\bea}{\begin{eqnarray}}
	\newcommand {\eea}{\end{eqnarray}}

\newcounter{mycounter}
\def\blot{\quad \mbox{$\vcenter{ \vbox{ \hrule height.4pt
				\hbox{\vrule width.4pt height.9ex \kern.9ex \vrule width.4pt}
				\hrule height.4pt}}$}}

\usepackage{natbib}
\bibpunct[, ]{(}{)}{,}{a}{}{,}%
\def\bibfont{\fontsize{9}{10.5}\selectfont}%

\TheoremsNumberedThrough     
\ECRepeatTheorems

\EquationsNumberedThrough   
\gdef\AQ#1{}
\gdef\CQ#1{}
\usepackage{xcolor}
\usepackage{xpatch}
\makeatletter
\def\changeBibColor#1{%
  \in@{#1}{}
  \ifin@\color{blue}\else\normalcolor\fi
}
 
\xpatchcmd\@bibitem
  {\item}
  {\changeBibColor{#1}\item}
  {}{\relax}
 
\xpatchcmd\@lbibitem
  {\item}
  {\changeBibColor{#2}\item}
  {}{\relax}
\makeatother

\begin{document}

\def\COPYRIGHTHOLDER{INFORMS}%
\def\COPYRIGHTYEAR{2017}%
\def\DOI{\fontsize{7.5}{9.5}\selectfont\sf\bfseries\noindent https://doi.org/10.1287/opre.2017.1714\CQ{Word count = 9740}}

\RUNAUTHOR{Xia, Ma, Sharif Azadeh} %

\RUNTITLE{Integrated timetabling and scheduling of modular autonomous vehicles under uncertainty}

\TITLE{Integrated timetabling and scheduling of modular autonomous vehicles under uncertainty}

\ARTICLEAUTHORS{
\AUTHOR{Dongyang Xia\textsuperscript{a,b}, Jihui Ma\textsuperscript{b}, Shadi Sharif Azadeh\textsuperscript{a,*}} 
\AFF{$^{a}$ Department of Transport \& Planning, Delft University of Technology, Netherlands\\
$^{b}$ Key Laboratory of Transport Industry of Big Data Application Technologies for Comprehensive Transport, Ministry of Transport, Beijing Jiaotong University, 100044, Beijing, China }%
}

\ABSTRACT{Addressing the Integrated Timetabling and Vehicle Scheduling (TTVS) problem is important for improving transit operations. Recently, the emerging modular autonomous vehicles composed of modular autonomous units have made it possible to dynamically adjust on-board capacity to better match space-time imbalanced passenger flows. This paper introduces an integrated framework for the TTVS problem in a dynamically capacitated and modularized bus network, considering time-varying and uncertain passenger demand. In this network, units can be (de)coupled and rerouted across different lines within the network at various times and locations, providing passengers with the opportunity to make in-vehicle transfers—that is, to transfer between lines while remaining onboard. We formulate a stochastic programming model to jointly determine the optimal robust timetable, dynamic formations of vehicles, and cross-line circulations of units, aiming to minimize the weighted sum of operators' and passengers' costs. To solve realistic instances, we propose a tailored integer L-shaped method to solve the formulated model dynamically through a rolling-horizon (RH) optimization algorithm. Furthermore, we extend our approach into a novel learning-based real-time decision-making framework that fine-tunes timetables and re-optimizes vehicle schedules in response to evolving and new demand realizations during practical operations. At its core is a scenario-retention method that selects a representative subset of scenarios using a machine learning model trained on scenario-level features. This subset is then incorporated into the optimization, ensuring both computational scalability and solution quality. To validate the effectiveness of our methods on realistic instances, we conduct experiments based on the Beijing bus network involving two bidirectional lines, 89 stops, up to 50 trips, and a 4-hour operational horizon. Our integrated optimization method outperforms the sequential approach. Compared to fixed-formation vehicles, our approach generates timetables and vehicle schedules that require fewer units. Additionally, the learning-based real-time decision-making framework outperforms benchmark algorithms in solution quality within a one-minute computation time limit.}

\KEYWORDS{Modularized bus network; Cross-line circulation; Rolling horizon framework; Integer L-shaped method; Machine learning}

\maketitle
	
\section{Introduction}
\label{sec:introduction}
In operations of the bus system, timetabling and vehicle scheduling link the quality of services with operational costs \citep{Ibarra-Rojas2014}. Passenger demand for bus networks is time-varying and uncertain due to factors such as the time of day, day of the week, season, special events, and unexpected disruptions. Within this context, the emerging technology of modular autonomous vehicles (MAVs) provides an innovative solution in bus operations, enhancing flexibility in catering to spatiotemporally varied passenger demand.

Specifically, MAVs represent a new generation of vehicles composed of interchangeable modular autonomous units (MAUs). The number of MAUs forming an MAV is referred to as its \textit{formation} in this study. For example, MAVs developed by Next Future Transportation Inc. with formations of 1 and 2 are shown in Figure A.1 in Online Appendix A \citep{nextMAV}. A key feature of MAVs is their \textit{flexible formations}, allowing dynamic composition adjustments en route and at depots through (de)coupling MAUs. Bus operations can benefit from the flexibility of MAVs by dynamically adjusting their capacity through (de)coupling units both en route and at depots and rerouting decoupled units to high-demand lines. We refer to this capacity adjustment as \textit {dynamic capacity allocation}. Additionally, MAVs feature an open, bus-like space between coupled MAUs, allowing passengers to move freely between units \citep{nextMAV}. This design enables \textit{in-vehicle transfers} of passengers, which refers to the feature allowing passengers to transfer between lines by walking from one on-board MAU to another without having to disembark, walk to other platforms, and board again. Several pilot studies have demonstrated promising results for modularized transit networks \citep{nextCompanyDubai}.

Timetabling at the network level has been proven to be NP-hard \citep{Caprara2002}. Integrating the vehicle scheduling decisions increases this computational challenge significantly. The literature has mostly explored the timetabling problem and the vehicle scheduling problem independently and sequentially, typically determining vehicle schedules only after timetables have been optimized. However, this sequential approach often results in suboptimal solutions, requiring more vehicles to realize the optimized timetables, and in extreme cases, generating timetables that exceed available fleet capacities. This limitation has already been highlighted in \cite{laporte2017multi}.

Therefore, the integration of timetabling and vehicle scheduling has gained attention in the railway domain. \cite{Rolf2021} demonstrated that timetables generated by the sequential optimization approach are strongly inferior to those produced through the integrated optimization method in realistic instances. This is because sequential approaches design timetables under the assumption of unlimited resources, often resulting in misalignment with practical operations and leading to inefficiencies. This issue is even more critical for MAVs, as their decision-making time windows are significantly shorter than those of railways. Moreover, overestimating or underestimating the required number of vehicles to meet an optimized timetable has a greater impact at the network level, further underlining the necessity of integrated optimization in modularized bus networks with MAVs.

However, even research addressing the timetabling and vehicle scheduling problem with fixed-formation vehicles at the network level remains limited, despite focusing on static passenger demand. Existing studies typically assume fixed-formation vehicles and have not fully utilized the advantages of MAVs that enable dynamic and flexible unit allocation at various locaions and times, and across multiple lines. This flexibility allows better alignment of supply with dynamic passenger demand. Recently, \citet{XIA2023} explored this framework on a single bus line, and \citet{Xia2024} extended it to an intermodal urban transit network with fixed-line and demand-responsive transport services. However, neither study accounts for in-vehicle passenger transfers or vehicle scheduling in a multi-line network. Additionally, their fleet sizing assumptions are overly optimistic, as they allow MAUs to appear or disappear at any station.  

In this paper, we explore the integrated timetabling, vehicle scheduling, and dynamic capacity allocation (TT-VS-DCA) problem in a network with MAVs under time-varying and uncertain demand, with explicit consideration of both the tactical level and the operational level. Given that the aforementioned timetable and vehicle schedule are obtained from optimization under historical demand scenarios, we further introduce a learning-based real-time decision-making framework to fine-tune them in response to newly emerging and time-varying demand realizations during practical operations, aiming to maintain high service quality. In the integrated optimization approaches, we examine the trade-off between passengers’ costs and operational costs (fleet size and dynamic-capacity usage of MAUs). Passengers’ costs are measured as a sum of waiting costs at origin stops and transfer waiting costs at transfer corridors, incurred during initial boarding and conventional transfers that require alighting, walking to change platform, and reboarding. The transfer waiting costs serve as the transfer penalties, reflecting the inconvenience associated with conventional transfers. 

The key idea of our integrated optimization approach is to jointly determine robust timetables, fleet sizes, and dynamic capacity allocations at different times and stops, while enabling \textit{cross-line vehicle circulations} to flexibly accommodate uneven spatiotemporal demand across the network. A \textit{vehicle circulation} in our study is defined as the complete working sequence of an MAU within the studied time horizon. A \textit{cross-line vehicle circulation} refers to a circulation in which an MAU serves trips on multiple lines, made possible by coupling, decoupling, and rerouting both at depots and at transfer stops. In comparison, traditional \textit{one-line circulations} can be defined as restricting each MAU to operate on a single bidirectional line and to return to the same depot in this study. This vehicle scheduling strategy can be seen as a special case of our cross-line method. A second special case is \textit{depot-based cross-line circulation}, where an MAU may be reassigned to a different line, but only at depots. Our cross-line circulation enables \textit{in-vehicle transfers}, where passengers remain onboard while transferring between lines. This avoids the inconvenience of exiting, walking, and re-boarding, thereby improving the passengers' traveling experience. By leveraging the flexible (de)coupling and re-routing capabilities, our approach supports more flexible fleet operations and enhances service responsiveness across the network. 

To model this, we first formalize and analyze the integrated optimization problem, explicitly incorporating the reconfigurability of MAUs and cross-line vehicle circulation, two features that have received limited attention in existing studies on timetabling and vehicle scheduling with fixed-formation vehicles.
We then propose a stochastic mixed-integer linear programming (MILP) formulation that accommodates scenarios with non-uniform probabilities. From the operator’s perspective, our formulation enables MAVs to change their formations both at depots and at transfer stops en route, and allows MAUs to execute trips across multiple lines at both depots and transfer stops through coupling and decoupling operations. From the passenger’s perspective, the model facilitates in-vehicle transfers by taking advantage of the modularity of MAUs, thereby avoiding the need to exit, walk to change platforms, and reboard when making transfers.

The investigated TT-VS-DCA problem in this study is as hard as the traditional timetabling and vehicle scheduling problem with fixed-formation vehicles at its minimum because it involves the joint design of timetables and vehicle schedules with flexible formations and cross-line circulations at various locations and times. Different from existing studies, the decoupling, coupling, and re-routing of MAUs at depots and en route are all allowed in our problem, requiring additional decision variables and coupling constraints. Thus, the solution approaches from previous studies cannot be directly applied to the TT-VS-DCA problem. To solve this problem, we develop an integer L-shaped method with tailored valid inequalities to solve the model. However, results show that the proposed exact method is ineffective for solving large-scale real-world instances. To address this computational challenge, we further design a rolling-horizon algorithm that solves the problem dynamically. In addition, we incorporate the proposed integer L-shaped method into the rolling-horizon optimization framework after finding that GUROBI, despite using the rolling-horizon approach, is unable to solve the TT-VS-DCA problem at the network level. As a result, instead of solving a large MILP, the developed algorithm solves a series of smaller MILPs. 

Furthermore, as operational conditions evolve (e.g., due to evolving passenger demand patterns that may fall outside the historical demand scenarios), the optimized timetables and vehicle schedules derived from the historical demand scenario set may become suboptimal, resulting in inefficient resource utilization or service shortages. To address this, we extend our approach to the real-time control level by proposing a novel learning-based real-time decision-making framework that fine-tunes tactical timetables within a small adjustment range and re-optimizes vehicle schedules accordingly under strict computational time limits. The effectiveness of our learning-based real-time decision-making algorithm stems from introducing a new and fundamental idea that carefully selects and incorporates a subset of scenarios into optimization. The first key concept of this algorithm involves a learning-based scenario-retention method, which selects a representative subset of scenarios using machine learning models trained on scenario-level features. This allows high-quality solutions to be generated in real time. The second key concept is the identification of features, which determines the most relevant scenario characteristics for predicting their criticality. Our learning-based real-time decision-making algorithm is novel because it provides a thorough method to construct a scenario subset to be integrated into the optimization of the stochastic programs to improve the scalability of the algorithm, and it provides a framework that can be generalized or extended in broader real-time optimization contexts. We conduct extensive computational experiments on real-world instances to validate the effectiveness, robustness, and scalability of our proposed methods. The results highlight the practical benefits of our integrated optimization approach, including the (de)coupling and cross-line vehicle circulation strategies at both depots and transfer stops, and demonstrate superior solution quality, scalability, and real-time applicability of our algorithms compared to benchmark approaches.

To summarize, the main contributions of this paper are (\roman{mycounter}\stepcounter{mycounter}) conceptualizing the TT-VS-DCA problem within a dynamically capacitated and modularized bus network that enables in-vehicle passenger transfers and cross-line unit circulations, (\roman{mycounter}\stepcounter{mycounter}) developing a stochastic MILP model, (\roman{mycounter}\stepcounter{mycounter}) designing effective solution methods, (\roman{mycounter}\stepcounter{mycounter}) developing a learning-based real-time decision-making algorithm to fine-tune schedules in response to new demand realizations during practical operations,  (\roman{mycounter}\stepcounter{mycounter}) illustrating that considerable savings in operational costs can be achieved by allowing flexible formations of MAVs and cross-line circulations of MAUs, and showing the benefits of our approaches.

The remainder of this paper is structured as follows. In Section \ref{sec:Literature}, we summarize related literature. In Section \ref{sec:problemDescription}, we formalize the investigated problem. In Section \ref{sec:mathematicalFormulation}, we describe the mathematical formulation. In Section \ref{sec:properties}, we analyze the properties of the proposed model and present the tailored valid inequalities. In Section \ref{sec:algorithm}, we describe the integer L-shape algorithm, the rolling-horizon optimization approach, the integrated solution method, and the learning-based real-time decision-making framework. In Section \ref{sec:experiments}, we discuss the results of the computational experiments. Finally, in Section \ref{sec:conclusion}, we summarize our findings.

\section{Literature review}
\label{sec:Literature}

In this section, we discuss research on the integrated timetabling and vehicle scheduling with fixed-formation vehicles, timetabling of MAVs, and real-time optimization under uncertainty.

\subsection{Integrated timetabling and vehicle scheduling with fixed-formation vehicles}

The traditional planning procedure for public transportation networks, as outlined by \cite{Huisman2004}, can be divided into the following sequential steps: first, the line plans and the corresponding departure frequency are designed; second, timetables are developed to minimize passengers' travel time (or waiting time); subsequently, vehicle schedules are determined to minimize operational costs, and finally, crew schedules are planned. The disadvantages of this sequential planning are mainly two-fold. Firstly, the main components of operational costs, i.e., the number of used vehicles and crews, are only determined at the final stage of the planning process \citep[]{Schobel2009}, which may result in generating inefficient operational plans or suboptimal solutions. Secondly, and perhaps more critically, a solution for the timetabling problem can turn out to be infeasible for the vehicle scheduling problem.

Models for the integrated timetabling and vehicle scheduling (TTVS) problem with fixed-formation vehicles focus on minimizing the weighted sum of  operational and passengers' costs, measured by the number of vehicles required to satisfy the optimal timetable, passengers' waiting time, and transfer penalties. The first category of the integrated TTVS problem is dedicated to the one-line circulation of vehicles with static demand \citep[e.g.,][]{Schobel2009, Ibarra-Rojas2014}. Considering one-ine vehicle circulation, \cite{laporte2017multi} presented an integrated optimization approach to jointly optimize timetables and vehicle schedules while incorporating information about users' preferences. The objectives include the minimization of users’ inconvenience and the operator’s costs. For the TTVS problem with the fixed-formation vehicles, depot-based cross-line circulations of vehicles, and static passenger demand, \cite{Rolf2021} proposed a MILP model to jointly optimize the periodic timetable and the vehicle schedule, where vehicles are allowed to be scheduled across lines at depots. The depot-based cross-line circulations of fixed-formation buses are further addressed by \cite{RolfPartB}. Recently, \cite{Amberg2023} presented a flow-based formulation for the integration of vehicle scheduling and crew scheduling in a bus network, which incorporates the controlled trip-shifting strategy. 

To sum up, only a few studies address the TTVS problem, and most of them consider fixed-formation vehicles and static and deterministic passenger demand. In particular, only \citet{Rolf2021} and \citet{RolfPartB} consider depot-based cross-line circulation, assuming that each trip can be served by a single vehicle and that passenger demand is both static and deterministic. We advance this research stream by integrating cross-line vehicle circulation into the TTVS problem under time-varying and uncertain demand. Our model allows each MAU to change its serving lines both during and after trips, and enables each MAV to adjust formations at both depots and transfer stops. These capabilities also support in-vehicle transfers, allowing passengers to transfer between lines without alighting and reboarding.

\subsection{Timetabling of MAVs with flexible formations}

With technological advancements, scheduling MAVs with flexible formations has attracted increasing attention in the timetabling literature for public transportation systems. In such systems, one benefit of MAVs is the reduction in the number of MAUs required for operations, enabled by flexible formation adjustments.

The majority of existing contributions consider the flexible formation of an MAV at the temporal level on a one-way line, allowing MAUs to be (de)coupled at depots while taking into account deterministic passenger demand \citep[e.g.,][]{CHEN20191, Chen2020, Xiaowei2021, XIA2023}. To further close the substantial gap between the passenger demand and vehicle capacity, a few recent studies have begun to focus on the timetabling of MAVs with station-wise flexible formation, i.e., an MAV can change its formation at any station on a single line \citep[e.g.,][]{Chen2022, TIAN2023103986, Xia2024}. Considering the uncertain passenger demand, related studies have employed stochastic programming, robust optimization, and distributionally robust optimization (DRO) approaches, aiming to improve the robustness of timetables and vehicle schedules \citep[e.g.,][]{AN2020572}. In the context of the timetabling of MAVs, \cite{XIA2023} proposed a DRO model on a modularized bus line and designed an exact solution method to solve it. They do not consider the vehicle scheduling aspect and assume that the fleet sizing is fixed. \cite{Xia2024} addressed the timetabling and vehicle scheduling problem within an integrated fixed-line and demand-responsive transport system. Despite considering the scheduling of units, they focus on fixed-line services on a bus line and direct services between two specific stops on this line, considering the vehicle scheduling within a bus network and passengers' in-vehicle transfers out of scope.  

In summary, existing literature considering MAVs have conceptually demonstrated that flexible coupling and decoupling strategies can reduce both the number of used units and associated operational costs. However, most of these studies focus on timetabling and neglect vehicle scheduling, resulting in an overly optimistic estimation of the potential reduction in fleet size. There is a noticeable gap concerning the joint optimization of timetabling and vehicle scheduling with flexible formations of MAVs and cross-line circulations of units to further reduce operational costs at the network level, which is significantly more challenging than at the line level due to the complexity of additional variables and the coupling constraints among timetables, vehicle schedules, and MAV formations within the entire network.

\subsection{Real-time optimization under uncertainty}
\label{sec:rt-review}

To better position our study, we briefly review two major streams of real-time optimization under uncertainty in broader domains such as public transport and ride-hailing. Robust and stochastic model predictive control (RMPC and SMPC) are the widely used approaches, repeatedly solving short-horizon problems within the rolling-horizon solution framework while accounting for uncertainty and operational constraints \citep[e.g.,][]{Ma2021,Bart2022,Wang2022}. For example, \citet{Ma2021} proposed an RMPC framework for bus regulation under demand uncertainty, where timetable adjustments are made in a rolling-horizon fashion using CPLEX.

Beyond RMPC and SMPC, another line of studies explores online algorithms that make per-epoch decisions with performance guarantees such as competitive ratios or regret bounds \citep[e.g.,][]{Cho2023,Lyu2024}. For instance, \citet{Lyu2024} developed an adaptive-weighting method that solves a single weighted-sum subproblem each epoch and uses a randomized policy across scenarios and time.

We contribute to this research line by proposing a learning-based decision-making framework for solving real-time optimization problems under uncertainty. In this framework, we develop a learning-based scenario-retention method that ranks all stochastic scenarios based on their criticality using machine learning and carefully selects a representative subset for optimization. This framework greatly enhances computational efficiency while maintaining high solution quality under strict real-time limits. This framework is broadly applicable to real-time stochastic optimization problems beyond the timetabling and vehicle scheduling problem studied in this study.

\subsection{Contributions against the state-of-the-art studies}

Table~\ref{tab:literature} in Appendix~\ref{sec:research} overviews closely related studies and highlights the current research gaps that motivate our study and contributions. In the existing literature on MAVs, most studies focus on timetabling for a single line, whereas vehicle scheduling, particularly at the network level, has received limited attention. The broader TTVS literature for fixed-formation fleets is likewise limited and typically assumes single-line circulations under static, deterministic demand. When depot-based cross-line circulation is considered, line changes of vehicles are usually restricted to depots; mid-trip line changing at transfer stops of units and passengers’ in-vehicle transfers remain largely unexplored. Besides, solution approaches from previous studies are tailored for TTVS with fixed-formation vehicles and cannot address the unique operational features of MAVs, such as in-vehicle transfers and cross-line circulation of MAUs. Our main contributions are as follows:

$\bullet$ We tackle the TT-VS-DCA problem in a dynamically capacitated and modularized bus network with time-varying and uncertain demand. Our approach allows MAVs to dynamically adjust formations and enables MAUs to be (de)coupled and rerouted both mid-trip and after trips, supporting each MAU in changing its serving line at depots and transfer stops. This operational flexibility enables in-vehicle transfers of passengers, allowing passengers to transfer without disembarking and reboarding. In contrast to existing methods that typically assume fixed formations of vehicles, restrict line changes of vehicles to depots, or do not allow every vehicle to change its serving line, our approach enables more efficient utilization of MAUs across the entire network.

$\bullet$ We extend the formulations proposed by \citet{XIA2023} and \citet{Xia2024} to a stochastic MILP model at the network level, explicitly incorporating in-vehicle transfers of passengers and cross-line vehicle circulation of MAUs. Key decisions include fleet sizing, departure times, buffer times at stops, the optimal number of MAUs comprising each MAV as it leaves the depot and each transfer stop, (de)coupling and rerouting of units to other lines at transfer stops. We also formulate strict capacity constraints, coupling constraints between MAVs and passengers, and vehicle scheduling constraints. 

$\bullet$ We develop an integer L-shaped method with tailored valid inequalities grounded in the theoretical properties of the proposed model, outperforming GUROBI on moderate-size instances. To tackle real-world instances with time-varying and uncertain passenger demand, where the flexibility of MAUs introduces high operational complexity, we propose a rolling-horizon (RH) optimization approach within a double-decomposition framework, integrating the integer L-shaped method. This framework first decomposes the planning horizon into smaller decision-making stages and then incorporates the integer L-shaped method within each stage. In this framework, RH is used to bring dynamism into the solution process. By decomposing the problem into a series of smaller MILPs, our approach improves tractability compared to solving a large MILP with GUROBI or the integer L-shaped method alone.

$\bullet$ To handle the evolving operational conditions in practical operations, we extend our approach to a novel learning-based real-time decision-making algorithm embedded in an RH framework with the integer L-shaped method. At each decision-making stage, this framework fine-tunes tactical timetables within a small adjustment range and re-optimizes vehicle schedules based on real-time conditions (e.g., time-varying passenger demand). To ensure computational tractability, we solve a stochastic MILP using the integer L-shaped method, incorporating a representative subset of scenarios. This subset is carefully selected through our proposed learning-based scenario-retention procedure, which ranks all scenarios based on predefined features using machine learning techniques and selects the top-$k$ most representative ones. By limiting the number of scenarios and ensuring that only the most informative ones are retained in the optimization at each decision-making stage, our framework achieves high solution quality while satisfying real-time computational requirements.

$\bullet$ 
Lastly, we validate our model and solution methods using real-world instances from the Beijing bus network. Extensive experiments demonstrate the effectiveness of our algorithms and the practical advantages of the integrated optimization approach over a sequential optimization method. Moreover, our proposed approach with flexible (de)coupling and cross-line circulation achieves considerable savings in the number of required MAUs, compared to approaches using fixed-formation vehicles or those allowing (de)coupling only at depots and depot-based cross-line circulation. Post hoc evaluations further show the robustness of our integrated optimization method in handling out-of-sample scenarios, as well as the adaptability of the fine-tuning strategy under demand surges and redistributions. In addition, we show that the learning-based real-time decision-making framework consistently produces high-quality solutions within one minute per decision-making stage. It outperforms both full scenario enumeration and the greedy scenario-retention method, demonstrating its scalability for real-time and large-scale applications.

\section{Problem description}
\label{sec:problemDescription}

To provide a formal definition of the TT-VS-DCA problem from a mathematical perspective, we describe the network structure and operational framework within this network in Section \ref{sec:networkStructure}. This is followed by the introduction of the time-dependent and uncertain passenger demand in Section \ref{sec:passengerDemand}. The operators' decisions and costs are summarized in Section \ref{sec:operatorDecision}. Lastly, the assumptions adopted in our model are introduced in Section \ref{sec_assimptions}.

\subsection{Network structure and operational framework}
\label{sec:networkStructure}

We first introduce the bus network. The structure of lines and trips within the network is then defined, followed by an introduction of MAVs. Lastly, we describe the operations of MAVs, including timetabling, vehicle scheduling, and dynamic capacity allocation. 

\subsubsection{Network structure.}  
We consider a bus network, represented as a directed graph $\mathcal{G}=\{\mathcal{S}, \mathcal{A}\}$, where $\mathcal{S}$ represents the set of line-dependent stops, and $\mathcal{A}$ represents the set of directed links between stops. The set $\mathcal{A}$ consists of two types of directed connections: \textit{sections} and \textit{transfer corridors}. Sections, denoted as $\tilde{\mathcal{A}}$, represent directed links between two consecutive stops along the same line. Transfer corridors, denoted as $\dot{\mathcal{A}}$, connect stops that correspond to the same physical location but belong to different lines, facilitating passenger transfers. Formally, we have  $\mathcal{A}=\tilde{\mathcal{A}} \bigcup \dot{\mathcal{A}}$. 

The network includes multiple lines, denoted as $\mathcal{L}$. Each \textit{line} $l \in \mathcal{L}$ is defined as an ordered sequence of line-dependent stops and directed sections that follows a specific operational direction. Each line is represented as a directed subgraph $\mathcal{G}_{l}=\{\mathcal{S}_{l}, \tilde{\mathcal{A}}_{l}\}$, where $\mathcal{S}_{l}$ denotes the set of stops on line $l$, and $\tilde{\mathcal{A}}_{l}$ represents the directed sections along the line. The set of stops on line $l$ is indexed as $\mathcal{S}_{l} =\{s_{1}, s_{2},\dots,s_{\mid \mathcal{S}_{l} \mid}\}$, ensuring that each stop is uniquely associated with its respective line. The set of sections along line $l$ is given by $\tilde{\mathcal{A}_{l}}=\{(l, i, l, i+1): i, i+1 \in \mathcal{S}_{l}, l \in \mathcal{L} \}$. 

Stops in the network are categorized into \textit{transfer stops} and \textit{non-transfer stops}. A transfer stop is a stop that appears on at least two different lines, allowing passengers to transfer between lines. The set of transfer stops on line $l$ is denoted as $\mathcal{R}_{l}$ while the set of non-transfer stops is denoted as $\hat{\mathcal{R}}_{l}$. Thus, the set of all stops on line $l$ satisfies $\mathcal{S}_{l}=\mathcal{R}_{l} \bigcup \hat{\mathcal{R}}_{l}, \forall l \in \mathcal{L}$. At transfer stops, there are transfer corridors to facilitate passenger transfers. The set of transfer corridors is denoted as $\dot{\mathcal{A}}=\{g: g=(l, i, l^{'}, i^{'}), i \in \mathcal{R}_{l}, i^{'} \in \mathcal{R}_{l'}, l \neq l^{'}, l, l^{'} \in \mathcal{L} \}$. In addition, the network includes depots, where MAVs are stored or (de)coupled. The set of depots is denoted as $\mathcal{M}$. Each bidirectional line is associated with two depots, located adjacent to its terminal stops. These depots are shared by the corresponding pair of unidirectional lines used to represent the bidirectional line in our model.

\subsubsection{Trips and MAVs.}
A \textit{trip} is defined as a scheduled movement of an MAV along a specific line $l \in \mathcal{L}$, starting at a depot and ending at the other depot. The set of trips on line $l$ is denoted as $\mathcal{K}_l=\{1,2,\dots,\left|\mathcal{K}_l\right|\}$. Each trip is executed by a single MAV, whose formation can be changed at both depots and transfer stops.

An MAV consists of a specific number of MAUs, each with a fixed capacity $\text{CAP}$. The MAV formation for a trip depends on the time-varying and uncertain passenger demand and the trip’s departure time. As a result, the relationship between MAVs and their formations plays a crucial role in operational planning.

\subsubsection{Timetabling, vehicle scheduling, and dynamic capacity allocation.}
The TT-VS-DCA problem studied here integrates tactical-level timetabling with operational-level vehicle scheduling and dynamic capacity allocation to generate optimal operating plans for a bus network using MAVs. The timetable specifies the departure and arrival times of MAVs at all stops on their assigned trips. Vehicle scheduling and dynamic capacity allocation determine, at each stop, the formation of the operated MAV, given passenger demand and operational constraints. At the start of any trip, an MAU operates as part of an MAV. From this perspective, two operational patterns arise: (i) an MAU may complete an entire trip from the first to the last stop as part of an MAV; after the trip, the MAU may be stored in a depot or reassigned to another trip on the same or a different line; (ii) an MAU may operate from the first stop to a transfer stop on the same line, after which it is decoupled from its original MAV and rerouted to another line, enabling efficient cross-line circulation of units throughout the network.

We allow each MAU to change its serving line at depots and transfer stops through (de)coupling and rerouting. We do not allow deadheading, so MAUs never move empty between depots or between a depot and a transfer stop. Rerouting at transfer stops is not limited to lines that share a depot; an MAU decoupled at a transfer stop may be rerouted to any connected line, even if this line has different depots. Example~\ref{ex:2} in Appendix~\ref{sec:examples} illustrates vehicle scheduling and dynamic capacity allocation with our proposed cross-line circulations of MAUs.

\subsection{Time-dependent and uncertain path-based passenger demand}
\label{sec:passengerDemand}

It is well known that passenger demand in bus networks is typically time-dependent and uncertain. Hence, we use the historical data to construct a set of stochastic passenger flow scenarios to capture the passenger demand in real-world operating environments. To facilitate modeling, we first discretize the studied time horizon into a finite set of time intervals with a duration $\Delta$, denoted by $\mathcal{T}=\{1, 2, \dots, \left|\mathcal{T}\right|\}$ and indexed by $t$. The set of stochastic passenger flow scenarios is represented as $\mathcal{W} = \{1, 2, \dots, w, \dots, \left|\mathcal{W} \right|\} $, where $\rho_{w}$ denotes the probability of each scenario $w \in \mathcal{W}$. 

In each scenario, passenger demand is time-dependent. To describe the dynamic evolution of passengers in the bus network, we adopt path-based passenger groups following \citet{LEE2022103505}. Each \textit{passenger group} $p$ belongs to the set $\mathcal{P}_w$ under scenario $w \in \mathcal{W}$ and consists of $n_p$ passengers who arrive at the origin stop $i^o_p$ at time $u_p$ and travel to the destination stop $i^d_p$ following the same sequence of lines and transfer facilities. We denote the ordered sequence of lines used by group~$p$ as $\mathcal{L}_p = (l_0, l_1, \ldots, l_{|\mathcal{L}_p|}) \subseteq \mathcal{L}$, and the ordered sequence of transfer corridors as $\mathcal{R}_p = (g_1, g_2, \ldots, g_{|\mathcal{R}_p|})$, where each $g$ belongs to the set of transfer corridors $\dot{\mathcal{A}}$.

We define a \textit{journey} as the complete movements of passenger group $p$ from its origin $i^o_p$ to its destination $i^d_p$. The journey consists of one or more \textit{legs}, where each leg is the travel along a single line between two consecutive transfers (or between the origin/destination and the nearest transfer). If no transfer is required, the journey has a single leg and the set of the ordered sequence of transfer corridors (i.e., $\mathcal{R}_p$) is empty. If transfers occur, the first leg starts at the origin and ends at the first transfer stop, subsequent legs connect successive transfer stops, and the final leg ends at the destination. The elements of $\mathcal{R}_p$ provide the links between consecutive legs of the journey.

Passenger movements within a modularized bus network consist of four phases. Passengers first arrive at their origin stop and wait for an MAV to begin their first leg. Once an MAV arrives, they board. Upon reaching the destination stop of their final leg, they alight. If a transfer is required, passengers alight at the designated transfer stop and begin the next leg of their journey by boarding another MAV. Unlike conventional bus systems, where passengers must disembark and wait for another scheduled vehicle, the modularized system allows for \textit{in-vehicle transfers}. If an MAU within the passenger group’s on-board MAV is scheduled to reroute to their transfer-in line, the transfer can occur without requiring passengers to disembark and reboard, thereby improving transfer convenience.

\subsection{Decisions and costs of operators}
\label{sec:operatorDecision}

Bus operators make three key decisions: timetabling, vehicle scheduling, and dynamic capacity allocation. These decisions are highly interdependent and must adapt swiftly to the time-dependent and uncertain nature of passenger demand, while also accounting for the dynamics of passengers within the network. 

\subsubsection{Timetabling decisions.} To provide passengers with a reliable and easy-to-remember schedule, operators must determine a robust timetable that remains consistent across all passenger flow scenarios occurring on different operational days. The timetable specifics arrival, departure, and dwell times at each stop for every trip on all lines. 

\subsubsection{Vehicle scheduling and dynamic-capacity allocations decisions.} Operators have to make decisions on the formation of the MAV assigned to each trip leaving the depot and each transfer stop, whether the MAV assigned to each trip is decoupled or coupled at transfer stops, as well as how the decoupled MAUs should be rerouted to serve new trips and lines. 

\subsubsection{Operational costs.} The costs of purchasing MAUs are formulated as the budget constraints in this study. Additionally, operational costs are formulated as a \textit{cost function of operators}, which correlate positively with the increase in the number of used MAUs assigned to execute all trips. We also explore the impact of the fleet size on the trade-off between operational and passengers' costs.

We can now formally define the investigated TT-VS-DCA problem as follows. 
\begin{definition}
    (The TT-VS-DCA problem). Given a bus network $\mathcal{G} = (\mathcal{S}, \mathcal{A})$, the stochastic set $\mathcal{W}$, the set of time-dependent and uncertain passenger group $\mathcal{P}_w$ under scenario $w$, and the group of trips $\mathcal{K}_l$ on line $l \in \mathcal{L}$, the goal of the TT-VS-DCA problem is to determine a robust timetable that specifies the arrival time and departure time of each trip on every stop along each line (determined by $z_{k, i, t}^{l}$,  which indicates whether an MAV assigned to trip $k$ on line $l$ departs stop $i$ at time interval $t$). Additionally, the problem involves determining the formation of the MAV assigned to each trip $k$ at each stop $i$ on line $l$ in scenario $w$ (denoted as $x^l_{k,i}(w)$), the cross-line circulation of decoupled MAUs at transfer corridor $g$ (denoted as $h^{g}_{k,k'}(w)$), and the number of MAUs stored at each depot $m$ at the beginning of the studied time horizon $\kappa_m$. The goal is to minimize the weighted sum of passengers' and operational costs while satisfying strict capacity constraints.

\end{definition}

\subsection{Assumptions}
\label{sec_assimptions}

To formulate the investigated problem, we adopt the following assumptions without loss of generality:

\textbf{Assumption 1.} We assume that the formation of an MAV can be reconfigured at depots and transfer stops through the coupling and decoupling of MAUs. These operations are assumed to require minimal time and are fully accommodated within the dwell time at transfer stops and the preparation time at depots, without introducing any additional time. Both dwell time and preparation time are explicitly incorporated into the model.

\textbf{Assumption 2.} We assume that when coupled, MAUs form an open, bus-like interior space that allows passengers to stand and move freely between units.

\textbf{Assumption 3.} We assume that in-vehicle transfers are completed before the MAV arrives at transfer stops, guided by onboard announcements. According to the operational plan, eligible passengers are informed in advance and can relocate to the appropriate MAU within the same coupled MAV. Upon arrival at the transfer stop, the MAV decouples, and the appropriate MAU with passengers will be rerouted to a different line.

Assumption 1 has been widely adopted in existing studies \citep[e.g.,][]{CHEN2021102388, Chen2022, Xia2024}. Assumption 2 reflects the practical design features of emerging MAU technologies \citep[e.g.,][]{nextInvehicle}. As for Assumption 3, unlike conventional transfers where passengers must alight, walk to another platform, wait for a different MAV, and board again, in-vehicle transfers involve no such actions. From the passenger’s perspective, the rerouting of the MAUs to another line feels no different from a typical stop at a non-transfer bus stop. They remain seated throughout the decoupling and coupling process. It is important to note that during a regular stop at a non-transfer bus stop, passengers also remain on board while the MAV dwells at this bus stop. This period is also not considered part of the passenger’s waiting time in our objective function. In the same spirit, the period associated with in-vehicle transfers is also excluded from the passenger’s waiting time in the objective function, as it does not impose additional effort on passengers, unlike conventional transfers. In contrast, conventional transfers incur waiting time due to the inconvenience of alighting, walking, and re-boarding, which we account for as a transfer penalty. This assumption highlights the superior comfort of in-vehicle transfers and supports the efficient design of timetables and (de)coupling strategies of MAUs.

\section{Mathematical formulation}
\label{sec:mathematicalFormulation}

In this section, we propose a stochastic programming formulation for the TT-VS-DCA problem. We begin with a concise formulation to outline the model framework in Section \ref{sec:conciseFormulation}. Section \ref{sec:objfunction} presents the detailed formulations of the objective function. In Section \ref{sec:timetable}, we delve into the timetabling constraints. In Section \ref{sec:passenger}, we formulate constraints on interactions between passengers' dynamics and MAVs' movements. Section \ref{sec:vehicleScheduling} presents the formulations related to vehicle scheduling. Finally, formulations for dynamic capacity allocations are provided in Section \ref{sec:dynamicAllocation}. All  notations are summarized in Tables \ref{tab:notations} and \ref{tab:dependentVariables} in Appendix ~\ref{sec:Notations}. 

\subsection{Concise formulation of the stochastic programming model for the TT-VS-DCA problem}
\label{sec:conciseFormulation}

Due to the complexity of the TT-VS-DCA problem, a large number of constraints are incorporated into the developed stochastic programming model. These constraints cover various aspects, including passengers' and operational costs (PO constraints), timetabling constraints (TT constraints), interactions between passengers’ dynamics and MAVs' movements (PM constraints), vehicle scheduling constraints (VS constraints), and interactions between vehicle scheduling and dynamic capacity allocation (VD constraints). To provide readers with a clear overview of the model structure, we first present the following concise formulation, which outlines the types of constraints and decision variables to illustrate the modeling framework:
\begin{mini!}|l|[2]
{\mathbf{z}, \mathbf{h}, \mathbf{x}, \boldsymbol{\kappa}}
{\varphi_1 \Psi^{pass}+ \varphi_2 \Psi^{oper}}{}{\label{eq:objective_brief}}
\addConstraint
{\text{(PO constraints \eqref{passenger_cost_sum}-\eqref{operation_cost_z})},}
{\label{eq:othercosts}}
{}
\addConstraint
{\text{(TT constraints \eqref{binary_z}-\eqref{cons_headway})},}
{\label{eq:otherTimetabling}}
{}
\addConstraint
{\text{(PM constraints \eqref{cons_binary_xi}-\eqref{cons_capacity}}),}
{\label{eq:otherPM}}
{}
\addConstraint
{\text{(VS constraints \eqref{cons_totalnumber_mv}-\eqref{con:downMv}}),}
{\label{eq:othervehicleScheduling}}
{}
\addConstraint
{\text{(VD constraints \eqref{cons_e}-\eqref{eq:dominY}}),}
{\label{eq:otherVD}}
{}
\addConstraint
{z_{k,i,t}^l}
{\in\{0, 1\} \label{zbinary}}
{\ \forall l\in \mathcal{L}, k \in \mathcal{K}_l, i \in \mathcal{S}_l, t \in \mathcal{T},} 
\addConstraint
{h^{g}_{k,k'}(w)}
{\in\{0, 1\} \label{eq:dominH} }
{\ \forall w\in\mathcal{W}, g=(l,i,l',i')\in \dot{\mathcal{A}}, k\in\mathcal{K}_l, k'\in\mathcal{K}_{l'}}
\addConstraint
{x_{k,i}^l(w)}
{\in \mathbb{N} \label{eq:dominX}}
{\ \forall w\in \mathcal{W}, l\in{\mathcal{L}}, k \in {\mathcal{K}}_l, i \in {\mathcal{S}}_l,  }
\addConstraint
{\kappa_m}
{\in \mathbb{N} \label{eq:dominV}}
{\ \forall m \in \mathcal{M}.}
\end{mini!}
The first term in the objective function (\ref{eq:objective_brief}), $\Psi^{pass}$, represents passengers’ costs in monetary units. Passengers’ costs include the waiting costs at origin stops and the transfer waiting costs. The transfer waiting costs account for the inconvenience of conventional transfers, such as disembarking, walking, waiting, and reboarding. The second term, $\Psi^{oper}$, denotes operational costs, also expressed in monetary units. The coefficients $\varphi_1$ and $\varphi_2$ serve as weighting factors. Constraints \eqref{zbinary} - \eqref{eq:dominV} define the domains of the decision variables. In the following sections, we provide detailed formulations of constraints \eqref{eq:othercosts}, \eqref{eq:otherTimetabling}, \eqref{eq:otherPM}, \eqref{eq:othervehicleScheduling}, and \eqref{eq:otherVD}, respectively.

\subsection{The objective function}\label{sec:objfunction}

The goal of this study is to minimize the weighted sum of passengers' and operational costs from a system-wide perspective. Below, we introduce the specific formulations for each cost component. First, from the perspective of passengers, the objective is to minimize the expected value of waiting costs (denoted as $\Psi^{pass}$) both at the origin stop of a passenger group $p$ in scenario $w$ (denoted as $\Psi_p^1(w)$) and at the transfer corridors (denoted as $\Psi_p^2(w)$), which can be expressed as follows:
\begin{align}
\Psi^{pass}=\sum_{w\in\mathcal{W}}\sum_{p\in\mathcal{P}_w}n_p\cdot\rho_w \cdot \vartheta^T \left(\varepsilon_1\cdot \Psi_p^1(w)+\varepsilon_2\cdot \Psi_p^2(w) \right),\label{passenger_cost_sum}
\end{align}
where $\varepsilon_1$ and $\varepsilon_2$ are weighting coefficients. $\vartheta^T$ represents the equivalent monetary value per unit of passengers’ waiting costs (unit: \$). $\rho_w$ and $n_p$ indicates the probability of each scenario $w$ and the number of passengers within the passenger group $p$.

We can now formulate the specific mathematical expressions of $\Psi_p^1(w)$ and $\Psi_p^2(w)$. Let a binary variable $\chi_{p,k}^l(w)$ denote whether passenger group $p$ boards the MAV assigned to trip $k$ on line $l$ at their origin stop in scenario $w$. If yes, $\chi_{p,k}^l(w)=1$; otherwise, $\chi_{p,k}^l(w)=0$. The waiting time $\Psi_p^1(w)$ at the origin stop $i^o_p$ of passenger group $p \in \mathcal{P}_{w}$ for all $w \in \mathcal{W}$ is formulated as follows:
\begin{align}
&\Psi_p^1(w) = \sum\limits_{k\in {\mathcal{K}}_l}\chi^l_{p,k}(w)(d_{k,i}^{l}-u_p) &\forall w\in{\mathcal{W}}, p\in\mathcal{P}, l=l_0\in\mathcal{L}_p,i=i^o_p,\label{z_wait}
\end{align}
where $u_p$ represents the arrival time of passenger group $p$ and $d_{k,i}^{l}$ denotes the departure time of trip $k$ at stop $i$ on line $l$.

To compute the waiting costs due to conventional transfers, we define two binary variables $\zeta_{p,k,k'}^{g}(w)$ and $h^{g}_{k,k'}(w)$. The variable $\zeta_{p,k,k'}^{g}(w)$ indicates whether passenger group $p$ transfers from the MAV assigned to trip $k$ to trip $k'$ at transfer corridor $g$ in scenario $w$. The variable $h^{g}_{k,k'}(w)$ represents whether MAUs on the MAV assigned to trip $k$ on line $l$ will be decoupled at transfer corridor $g=(l,i,l',i')\in \dot{\mathcal{A}}$ and rerouted to execute the $k'$th trip on line $l'$ in scenario $w$. The waiting costs at transfer corridors depend on the type of transfer, which can be categorized as either an in-vehicle transfer or a conventional transfer that requires alighting and reboarding.

1. \textit{In-vehicle transfers.} An in-vehicle transfer occurs when the MAUs on the MAV assigned to trip $k$ on line $l$ are scheduled to be decoupled at transfer corridor $g$ and rerouted to serve the $k'$th trip on line $l'$, i.e., when $h^{g}_{k,k'}(w) = 1$, and there exits an on-board passenger group $p$ requiring a transfer to that line, i.e., $\zeta_{p,k,k'}^{g}(w)=1$. In this case, passenger group $p$ can make an in-vehicle transfer, benefiting from its convenience, with zero transfer waiting cost.

2. \textit{Conventional transfers requiring alighting and reboarding.} If either of the aforementioned conditions is not met, passengers must alight from the MAV assigned to trip $k$ on line $l$, walk to and wait at the transfer corridor $g$, and board another MAV assigned to trip $k'$ on line $l'$. In this case, the waiting time in the transfer corridor is given by the difference between the departure time of the transfer-in MAV assigned to trip $k'$ (i.e., $d_{k',i'}^{l'}$) and the arrival time of the transfer-out MAV assigned to trip $k$ (i.e., $a_{k,i}^{l}$), formulated as  $d_{k',i'}^{l'}-a_{k,i}^{l}$.

To sum up, transfer waiting costs can be formulated in a single equation as follows:
\begin{align}
&\Psi_p^2(w) = \sum_{g=(l,i,l',i')\in \mathcal{R}_p}\sum_{k\in\mathcal{K}_l}\sum_{k'\in\mathcal{K}_{l'}}\zeta_{p,k,k'}^{g}(w)\cdot(1-h^g_{k,k'}(w)) \cdot (d_{k',i'}^{l'}-a_{k,i}^{l}) &\forall w\in{\mathcal{W}}, p\in\mathcal{P}_w,\label{z_transfer}
\end{align}
where the terms $\zeta_{p,k,k'}^{g}(w)(1-h^g_{k,k'}(w))\cdot d_{k',i'}^{l'}$ and $\zeta_{p,k,k'}^{g}(w)(1-h^g_{k,k'}(w))\cdot a_{k,i}^{l}$ are both nonlinear. We derive the linearized form of nonlinear constraints \eqref{z_wait} and \eqref{z_transfer} in Appendix \ref{sec:Linearization}. 

On the other hand, from the perspective of operators, the objective is to minimize operational costs. These costs are associated with the number of MAUs serving each section between any two adjacent stops on every line in each scenario. We define a binary decision variable $y^l_{k,i,q}(w)$ to indicate whether the number of MAUs composed in the MAV assigned to trip $k$ serving the section between stops $i$ and $i+1$ on line $l$ in scenario $w$ is $q \in \mathcal{Q}$ or not, where $\mathcal{Q}$ denotes the set of the number of MAUs that can be contained in a MAV. A cost-related parameter associated with using $q$ MAUs when an MAV departs from stop~$i$ is denoted by $\vartheta_{q,i}$. The formulation for the expected value of operational costs among all scenarios (denoted as $\Psi^{oper}$) can be presented as follows:
\begin{align}
&\Psi^{oper} = \sum\limits_{w\in\mathcal{W}}\rho_w \Psi_{op}(w)=\sum_{w\in\mathcal{W}}\sum_{l\in\mathcal{L}}\sum_{i\in\mathcal{S}_l}\sum_{k\in\mathcal{K}_l}\rho_w \cdot \sum_{q\in\mathcal{Q}}\vartheta_{q,i}\cdot y_{k,i,q}^l(w),
\label{operation_cost_z}\end{align}
where $\Psi_{op}(w)$ refers to the operational costs in each scenario $w \in \mathcal{W}$.

To facilitate the trade-offs between operators and passengers, the aforementioned objective functions are reformulated as a single
one through the weighting coefficients $\varphi_1$ and $\varphi_2$, shown as below:
\begin{align}
    &\Psi=\varphi_1 \Psi^{pass}+ \varphi_2 \Psi^{oper}.
\end{align}

\subsection{Timetabling constraints}
\label{sec:timetable}

The timetable specifies the arrival and departure times of each trip at every stop along each line. We introduce the parameters $\alpha_{l,k}^{min}$ and $\alpha_{l,k}^{max}$ as the lower and upper bounds of the shifting time at the first stop of trip $k$ on line $l$. Similarly, $\beta_{l,k,i}^{min}$ and $\beta_{l,k,i}^{max}$ define the dwell time limits of trip $k$ on line $l$ at stop $i$, while $h_{min}^l$ and $h_{max}^l$ constrain the headway between two consecutive trips on line $l$. The running time between stops $i-1$ and $i$ on line $l$ for trip $k$ is denoted as $r_{k,i-1}^l$.  

To model timetabling decisions, we define the binary decision variable $z_{k, i, t}^{l}$ to indicate whether an MAV assigned to trip $k$ on line $l$ departs stop $i$ at time interval $t$. $z_{k,i,t}^l=1$ indicates trip $k$ on line $l$ has not departed from stop $i$; otherwise, $z_{k,i,t}^l=0$. The real-valued variables $a_{k,i}^l$ and $d_{k,i}^l$ represent the arrival and departure times at stop $i$, respectively. Additionally, $\alpha_{k,1}^l$ captures the shifting time at the first stop, while $\beta_{k,i}^l$ denotes the dwell time at stop $i$. The timetable is governed by the following constraints:
\begin{align}
    &z_{k,i,t+1}^l\leq z_{k,i,t}^l &&\forall l\in \mathcal{L}, k \in \mathcal{K}, i \in \mathcal{S}, t \in \mathcal{T} \backslash\{\left|\mathcal{T}\right|\}. \label{binary_z}\\
    &d_{k,i}^l=\Delta (1+ \sum_{t\in{\mathcal{T}}}z_{k,i,t}^l) &&\forall l\in{\mathcal{L}}, k \in {\mathcal{K}}_l, i \in {\mathcal{S}_l}. \label{cons_couple}\\
    &a_{k,i}^l = \begin{cases}
o_{k,1}^l+ \alpha_{k,1}^l  &\text{if} \ i = 1\\\
d_{k,i-1}^l + r_{k,i-1}^l &\text{otherwise}\ \ 
\end{cases} &&\forall l\in{\mathcal{L}}, k \in {\mathcal{K}}_l.  \label{cons_arrive} \\
&\beta_{k,i}^l  = d_{k,i}^l - a_{k,i}^l && \forall l\in{\mathcal{L}}, k \in {\mathcal{K}}_l, i \in {\mathcal{S}_l}. \label{cons_buffer}\\
&\alpha_{l,k}^{min}\leq\alpha_{k, 1}^l\leq\alpha_{l,k}^{max} &&\forall l\in{\mathcal{L}}, k \in {\mathcal{K}}_l. \label{cons_shift}\\
&\beta_{l,k,i}^{min}\leq \beta_{k,i}^l \leq\beta_{l,k,i}^{max} &&\forall l\in{\mathcal{L}}, k \in {\mathcal{K}}_l, i \in {\mathcal{S}_l}. \label{cons_dwell}\\
&h_{min}^l\leq d_{k+1,i}^l-d_{k,i}^l\leq h_{max}^l &&\forall l\in{\mathcal{L}}, k \in {\mathcal{K}}_l\backslash\{\left|\mathcal{K}_l\right|\}\}, i \in {\mathcal{S}_l}. \label{cons_headway}
\end{align}

Constraints (\ref{binary_z}) ensure the monotonically non-increasing property of the binary variable related to timetabling. Constraints (\ref{cons_couple}) define the relationship between the binary variable related to timetabling and the real-valued departure times. Constraints (\ref{cons_arrive}) calculate the arrival time at stop $i$ of trip $k$ on line $l$. Constraints~(\ref{cons_buffer}) are formulated to track the dwell time of trip~$k$ at stop~$i$ on line~$l$, denoted by $\beta_{k,i}^{l}$. This dwell time is defined as the difference between the departure time $d_{k,i}^{l}$ and the arrival time $a_{k,i}^{l}$ at stop~$i$ on line~$l$ for trip~$k$. Additionally, constraints (\ref{cons_shift}), (\ref{cons_dwell}), and (\ref{cons_headway}) set the lower and upper limitations on shifting times, dwell times, and headway, respectively. 

\subsection{Constraints on the interactions between passengers' dynamics and MAVs' movements}
\label{sec:passenger}

Passenger dynamics in the bus network are closely linked to MAV movements and capacity constraints. These interactions are formulated in two subsections: the first subsection presents interactions between dynamics on different legs of the passengers' journey and MAVs, while the second one formulates interactions between dynamics of passengers and capacity of MAVs.

\subsubsection{Interactions between dynamics on different legs of the passengers' journey and MAVs} 

Theoretically, a passenger can board all trips arriving after his arrival at the origin station. In real-world operations, a passenger waiting at a stop usually chooses the first available MAV capable of completing their journey. To model relations between passengers' dynamics and movements of MAVs, we define a binary variable $\xi^l_{p,k}(w)$ that indicates whether passenger group $p$ can board the MAV assigned to trip $k$ on line $l$ in scenario $w$ or not for the first leg of their journey. Specifically, this variable takes value 1 if the MAV assigned to trip~$k$ on line~$l$ departs from the origin stop after the arrival time of passenger group~$p$ in scenario~$w$, and 0 otherwise.

In the first leg of the journey, before any transfers, passenger group~$p$ has the opportunity to board any trip~$k$ on line~$l$ whose departure time $d_{k,i^o_p}^l$ is later than the group’s arrival time $u_p$ at the origin stop $i^{o}_{p}$. That is, 
\begin{align}
& M_1(\xi_{p,k}^l(w)-1) \leq d_{k,i^o_p}^l - u_p\leq M_1\xi_{p,k}^l(w) - \epsilon_1 &&\forall w\in{\mathcal{W}}, p\in{\mathcal{P}_w}, l=l_0\in{\mathcal{L}}_p, k\in {\mathcal{K}}_l.\label{cons_binary_xi}     
\end{align}

Thereafter, if passenger group $p$ needs to transfer from the MAV assigned to trip $k$ on line $l$ to trip $k'$ on line $l'$ at transfer corridor $g$, the departure time $d_{k',i'}^{l'}$ of the transfer-in trip has to be later than the arrival time $a_{k,i}^{l}$ of the transfer-out trip. To model this, we introduce a binary variable $\pi^{g}_{k,k'}$, which indicates whether passengers have the opportunity to transfer from the MAV assigned to trip $k$ to the other MAV allocated to trip $k'$ at transfer corridor $g$. This condition is formulated by the following constraints:
\begin{align}
& M_2(\pi^{g}_{k,k'}-1) \leq d_{k',i'}^{l'}-a_{k,i}^{l} - \underline{\theta}_g \leq M_2\pi^{g}_{k,k'}- \epsilon_2 &&\forall g=(l,i,l',i')\in \mathcal{\dot{\mathcal{A}}}, k\in {\mathcal{K}}_l, k'\in {\mathcal{K}}_{l'},\label{cons_pi}
\end{align}
where the parameter $\underline{\theta}_g$ represents the minimum time required for a transfer at transfer corridor $g$.

Note that constraints \eqref{cons_pi} capture all possible transfer options for a passenger group $p$ arriving at a transfer corridor $g$. However, since we assume passengers always board the first available trip, additional constraints should be introduced to enforce this selection. To do so, we formulate nonlinear constraints (\ref{cons_indicator_transfer}) to link passenger group $p$ and its actually boarded trip. Specifically, passenger group $p$ actually transfers from the MAV assigned to trip $k$ to trip $k'$ at transfer corridor $g$ in scenario $w$ ($\zeta_{p,k,k'}^{g}(w) =1$) if and only if this group was on trip $k$ on line $l$ (i.e., $\chi_{p,k}^{l}(w) =1$) and trip $k'$ on line $l'$ has a later departure time than the arrival time of trip $k$ and is the closest among all possible transfer options on line $l'$ (i.e., $\pi^{g}_{k,k'}-\pi^{g}_{k,k'-1}=1$). As a result, passenger group $p$ actually transfer from trips $k$ to $k'$ at the transfer corridor $g$. The linearized form of this set of constraints is given in Appendix \ref{sec:Linearization}. 
\begin{align}
& \zeta_{p,k,k'}^{g}(w)=\chi_{p,k}^l(w)\cdot(\pi^{g}_{k,k'}-\pi^{g}_{k,k'-1}) &&\forall w\in{\mathcal{W}} ,p\in\mathcal{P}_w,g=(l,i,l',i')\in\mathcal{R}_p, k\in {\mathcal{K}}_l, k'\in {\mathcal{K}}_{l'}.\label{cons_indicator_transfer}
\end{align}

Additionally, constraints (\ref{eq:chi}) are proposed to track the MAV assigned to trip \( k \) that passenger group \( p \) boards when transferring from trip \( k' \) on line \( l' \) to trip \( k \) on line \( l \) in scenario \( w \). If line \( l \) is the line where the origin stop of passenger group \( p \)'s journey is located (i.e., $l = l_0$), \( \chi^l_{p,k}(w) \) indicates whether the group boards trip \( k \) in the first leg. If line \( l \) is a transfer-in line within the journey, the variable determines whether the group successfully transfers to trip \( k \) after departing from trip \( k' \). The corresponding formulation is as follows:
\begin{align}
&  \chi^l_{p,k}(w) =
\begin{cases}
\xi_{p,k}^{l}(w) - \xi_{p,k-1}^{l}(w), & \text{if } l = l_0\\
\sum\limits_{k'\in\mathcal{K}_{l'}} \zeta^{g}_{p, k',k}(w), &\begin{subarray}{l} \text{if } l \neq l_0,\\ g=(l',i',l,i)\in R_p \end{subarray}
\end{cases}
&& \forall w\in\mathcal{W}, p\in\mathcal{P}_w,l\in\mathcal{L}_p, k\in\mathcal{K}_l. \label{eq:chi}
\end{align}

\subsubsection{Interactions between dynamics of passengers and capacity of MAVs} 

Let variables $b_{k,i}^{l}(w)$ and $c_{k,i}^{l}(w)$ represent the numbers of passengers boarding and alighting trip $k$ at stop $i$ on line $l$ in scenario $w$. The variable $v_{k,i}^{l}(w)$ denotes the number of in-vehicle passengers when the MAV assigned to trip $k$ departs from stop $i$ on line $l$ in scenario $w$. We define an integer variable $x^l_{k,i}(w)$ to indicate the real-valued number of MAUs composed in the MAV allocated to trip $k$ serving the section between stops $i$ and $i+1$ on line $l$. We then formulate constraints related to passenger movements as follows:
\begin{align}
    &b_{k,i}^l(w) = \sum_{p\in\mathcal{PO}_{i}}n_{p}\cdot\chi^l_{p,k}(w) && \forall w\in{\mathcal{W}}, l\in{\mathcal{L}}, k \in {\mathcal{K}}_l, i \in {\mathcal{S}}_l. \label{cons_board}\\
&c_{k,i}^l(w) = \sum_{p\in\mathcal{PD}_{i}}n_{p}\cdot\chi^l_{p,k}(w) && \forall w\in{\mathcal{W}}, l\in{\mathcal{L}}, k \in {\mathcal{K}}_l, i \in {\mathcal{S}}_l.\label{cons_alight}\\
&v_{k,i}^l(w) = \sum_{j\leq i}(b_{k,j}^l(w)-c_{k,j}^l(w)) &&\forall w\in{\mathcal{W}}, l\in{\mathcal{L}}, k \in {\mathcal{K}}_l, i,j \in {\mathcal{S}}_l. \label{cons_invehicle}\\
&v_{k,i}^l(w)\leq \text{CAP} \cdot  x_{k,i}^l(w) &&\forall w\in{\mathcal{W}}, l\in{\mathcal{L}}, k \in {\mathcal{K}}_l, i \in {\mathcal{S}}_l.\label{cons_capacity}
\end{align}

Constraints (\ref{cons_board}) compute the number of passengers boarding the MAV assigned trip $k$ when it arrives at stop $i$, where set $\mathcal{PO}_{i}$ consists of all the passengers who board at stop $i$. Similarly, constraints (\ref{cons_alight}) are formulated to compute the number of passengers alighting at stop $i$, where the set $\mathcal{PD}_{i}$ includes all the passengers who alight at stop $i$. Constraints (\ref{cons_invehicle}) compute the number of in-vehicle passengers when the MAV assigned to trip $k$ leaves stop $i$ on line $l$, which equals to the difference between the cumulative number of boarding and alighting passengers.  Constraints (\ref{cons_capacity}) ensure that the number of in-vehicle passengers cannot exceed the capacity.

\subsection{Constraints associated with vehicle scheduling}\label{sec:vehicleScheduling}

In this study, we assume that MAVs allocated to trips can be decoupled and coupled at transfer stops and depots. Thus, the number of usable MAUs stored at the depots at each time interval during operations should be sufficient to perform trips. We introduce two variables $AV_{m,t}(w)$ and $DV_{m,t}(w)$ to track the cumulative inflow and outflow of MAUs at depot $m \in \mathcal{M}$ at time $t$ in scenario $w$. The decision variable $\kappa_m$ is defined to indicate the number of MAUs stored at each depot $m \in \mathcal{M}$ at the beginning of the study time horizon. We have the following constraints:
\begin{align}
 &\sum\limits_{m \in\mathcal{M}}\kappa_m\leq \text{B}. \label{cons_totalnumber_mv}\\
    & DV_{m,t}(w)\leq \kappa_m + AV_{m,t}(w) &&\forall m\in\mathcal{M}, t\in\mathcal{T}, w\in\mathcal{W}. \label{cons_number_mv}\\
    & AV_{m,t}(w) =  \sum_{l\in LL_m}\sum_{k\in\mathcal{K}_l}(1-z_{k,\left|\mathcal{S}_l\right|,t-t^{pre}_m}^l)\cdot x_{k,\left|\mathcal{S}_l\right|}^l(w) && \forall m\in\mathcal{M}, t\in\mathcal{T}, w\in\mathcal{W}.\label{cons_up_mv} \\
& DV_{m,t}(w) =  \sum_{l\in FL_m}\sum_{k\in\mathcal{K}_l}(1-z_{k,1,t}^l)\cdot x_{k,1}^l(w) && \forall m\in\mathcal{M}, t\in\mathcal{T}, w\in\mathcal{W}. \label{con:downMv} 
\end{align}

Constraints (\ref{cons_totalnumber_mv}) set the upper bound of MAUs used to cover all trips, where the parameter \text{B} represents the maximum number of available MAUs. Constraints (\ref{cons_number_mv}) are imposed to ensure that the cumulative number of MAUs exiting the depot $m$ cannot exceed its own storage and inflow. Constraints \eqref{cons_up_mv} and \eqref{con:downMv} are formulated to compute the inflows and outflows of MAUs at time interval $t$ and depot $m$ in scenario $w$, where $FL_m$ represents the sets of lines whose starting stop is depot $m$ and $LL_m$ refers to the set of lines for which the ending stops is depot $m$. $t^{pre}_m$ is the preparation time at depots of each MAU to complete a trip before it is reallocated by the next trip. Appendix \ref{sec:Linearization} presents the linearization of the nonlinear constraints \eqref{cons_up_mv} and \eqref{con:downMv}.

\subsection{Constraints on the interactions between vehicle scheduling and dynamic capacity allocation}\label{sec:dynamicAllocation}

We assume that rerouting operations are performed at the transfer stops only when the time window requirements are met. The reason for this is twofold: the above operations take a certain length of time to perform; on the other hand, if the decoupled units are stored for a considerable period of time in a transfer stop without being able to be re-routed, their utilization rate will be lower and they will take up the limited space in the transfer stop, bringing inconvenient for both operations and passengers. We define the variable $e^{g}_{k,k'}$ to indicate whether units assigned to trip $k$ can be decoupled and coupled to the MAV executing trip $k'$ at the transfer corridors $g$ or not. In other words, whether the decoupling and coupling behaviors are feasible in the time dimension. Recall that we have defined the variable $h^{g}_{k,k'}(w)$ to indicate whether MAUs on the MAV assigned to trip $k$ on line $l$ will be decoupled at transfer corridor $g=(l,i,l',i')\in \dot{\mathcal{A}}$ and rerouted to execute the trip $k'$ on line $l'$ in scenario $w$. We formulate the following constraints to describe the decoupling and coupling operations at transfer stops.
\begin{align}
    &e^{g}_{k,k'}=
    \begin{cases}
    1 & \text{if} \ \underline{\theta}_g\leq d^{l'}_{k',i'}-a^{l}_{k,i}\leq \overline{\theta}_g \\
    0 &\text{otherwise}\ \ 
\end{cases} \ \ \qquad\qquad \quad \forall g=(l,i,l',i')\in\dot{\mathcal{A}}, k\in\mathcal{K}_{l}, k'\in\mathcal{K}_{l'}.\label{cons_e}\\
&h^{g}_{k,k'}(w)\leq e^{g}_{k,k'} \qquad\qquad\qquad\qquad\qquad\qquad\qquad\forall w\in\mathcal{W}, g=(l,i,l',i')\in\dot{\mathcal{A}}, k\in\mathcal{K}_l, k'\in\mathcal{K}_{l'}.\label{cons_he}\\
& h^{g}_{k,k'}(w) + h^{g'}_{k',k}(w)\leq 1 \qquad\qquad\qquad\qquad\qquad\ \ \ \forall w\in\mathcal{W}, g=(l,i,l',i')\in \dot{\mathcal{A}}.\label{cons_h1h2}
\end{align}

Constraints (\ref{cons_e}) define the temporal bounds $[\underline{\theta}_g, \overline{\theta}_g]$ for rerouting the (de)coupling MAUs across lines, which ensure rerouting the decoupled MAUs from trip $k$ on line $l$ to serve trip $k'$ on line $l'$ is only possible if the time difference between the departure time of trip $k'$ (i.e., $d^{l'}_{k',i'}$) and the arrival time of trip $k$ (i.e., $a^{l}_{k,i}$) is within the predetermined time window. This set of nonlinear constraints is linearized in Appendix \ref{sec:Linearization}. Constraints (\ref{cons_he}) are formulated to ensure that only after the aforementioned temporal conditions of decoupling and rerouting are met, this operation is only practicable.  Constraints (\ref{cons_h1h2}) guarantee the unidirectional rerouting and coupling of the decoupled MAUs. The reason for developing constraints (\ref{cons_h1h2}) is that if MAUs composed in the MAV assigned to trip $k$ on line $l$ are decoupled, rerouted to line $l'$, and coupled on the MAV serving trip $k'$, this indicates that additional MAUs are required to successfully complete trip $k'$. Consequently, this also indicates that further decoupling of the MAV assigned to trip $k'$ is not necessary.

Finally, we formulate the following constraints to model the dynamics of the formations of MAVs during operations and to characterize the connections between the formations and the en-route (de)coupling operations at transfer stops. The binary variable $y^l_{k,i,q}(w)$ is defined to indicate whether the MAV allocated to perform the trip $k$ at stop $i$ on line $l$ in scenario $w$ is comprised by $q$ MAUs. If yes, we have $y^l_{k,i,q}(w) = 1$; otherwise, it equals to 0.
\begin{align}
 &x^l_{k,i}(w)=\sum_{q\in\mathcal{Q}}q \cdot y^l_{k,i,q}(w) \qquad\qquad\qquad\qquad\qquad\forall w\in\mathcal{W}, l\in{\mathcal{L}}, k \in {\mathcal{K}}_l, i \in {\mathcal{S}}_l.\label{cons_xq} \\
&x^l_{k,i-1}(w)=x^l_{k,i}(w) \ \qquad\qquad \qquad \qquad  \qquad \qquad \forall w\in\mathcal{W}, l\in{\mathcal{L}}, k \in {\mathcal{K}}_l, i \in \hat{\mathcal{R}}_{l}.\label{cons_x}\\
&x^l_{k,i}(w) - x^l_{k,i-1}(w) + 1\leq M_3(1-h^{g}_{k,k'}(w))   \nonumber\\ &\qquad\qquad\qquad\qquad\qquad\qquad \forall w\in\mathcal{W}, l\in{\mathcal{L}}, k\in\mathcal{K}_l,k'\in\mathcal{K}_l',i \in \mathcal{R}_{l}, g=(l,i,l',i')\in\dot{\mathcal{A}}. \label{cons_xh_1}\\
&x^l_{k,i}(w) = x^l_{k,i-1}(w) +\sum_{l'}\sum_{k'\in\mathcal{K}_{l'}} (h^{g'}_{k',k}(w) - h^{g}_{k,k'}(w))(x^{l'}_{k',i'-1}(w) - x^{l'}_{k',i'}(w))    \nonumber\\ &\qquad\qquad\qquad\qquad\forall  w\in\mathcal{W}, l\in{\mathcal{L}}, k\in\mathcal{K}_l,i \in \mathcal{R}_{l}, g=(l,i,l',i')\in\dot{\mathcal{A}}, g'=(l',i',l,i)\in\dot{\mathcal{A}}\label{cons_xh_2}.\\
&y^l_{k,i,q}(w)\in\{0, 1\} \qquad\qquad \qquad \qquad  \qquad\qquad \qquad\forall w\in \mathcal{W}, l\in{\mathcal{L}}, k \in {\mathcal{K}}_l, i \in {\mathcal{S}}_l, q \in \mathcal{Q}.  \label{eq:dominY}
\end{align}
Constraints (\ref{cons_xq}) are imposed to compute the real-valued number of MAUs comprising this MAV at stop $i$ on line $l$ (i.e., $x_{k,i}^{l}(w)$ ). Constraints (\ref{cons_x}) are developed to guarantee that the formation of an MAV executing trip $k$ remains unaltered upon arrival and departure at stops other than transfer stops. Moreover, given that MAVs allocated to trips have the flexibility to be decoupled or coupled at transfer stops, the formation of MAVs can potentially change as they depart from transfer stops. Following this, constraints (\ref{cons_xh_1}) link routing and decoupling operations at the transfer stop (i.e., when $h_{k,k'}^{g}(w)=1$) with modifications to the number of MAUs that make up the MAV to be decoupled. Nonlinear constraints (\ref{cons_xh_2}) are formulated to calculate the number of MAUs included in the MAV that is assigned to execute trip $k$ upon departing from transfer stop $i$ on line $l$. We derive the linearized form of this set of nonlinear constraints in Appendix \ref{sec:Linearization}. Constraints (\ref{eq:dominY}) define the domain.

In summary, we propose the following MILP model for the TT-VS-DCA problem in the modularized bus network with MAVs:
\begin{mini}|l|[0]
{\mathbf{z}, \mathbf{h},  \mathbf{x}, \boldsymbol{\kappa}}
{\varphi_1 \Psi^{pass}+ \varphi_2 \Psi^{oper}}{\label{stochastic_problem}}{} 
\addConstraint
{}
{\eqref{zbinary} - \eqref{eq:dominV}, \eqref{passenger_cost_sum}, \eqref{operation_cost_z}, \eqref{binary_z} - \eqref{cons_pi}, \eqref{eq:chi} - \eqref{cons_number_mv}, \eqref{cons_he} - \eqref{cons_xh_1},\eqref{eq:dominY}, \eqref{z_wait_linear} - \eqref{cons_hh_linear}.}
\end{mini}

\section{Theoretical properties of the TT-VS-DCA model}
\label{sec:properties}
In this section, we establish propositions based on the mathematical properties of the proposed model, which serve as the foundation for the solution methods introduced in the following section. The TT-VS-DCA problem studied in this work is NP-hard. As the problem size increases, the number of combinations across line $l \in \mathcal{L}$, trip $k \in \mathcal{K}_l$, stop $i \in \mathcal{S}_l$, and discretized time intervals $t \in \mathcal{T}$ grows exponentially, substantially increasing the computational complexity of the problem. To address this challenge, we propose the following propositions to strengthen the model.

Specifically, Proposition~\ref{propo_domin_reduction} introduces a time window-based reduction that restricts the feasible domain of timetabling variables. Proposition~\ref{propo_first_board} encodes boarding feasibility by eliminating logically invalid boarding times. Proposition~\ref{propo_transfer} ensures correct temporal sequencing for passenger transfers. Proposition~\ref{propo_coupling} provides a compact temporal window that enables the decoupling of an MAU at a transfer stop and its subsequent reassignment to another line, ensuring the temporal feasibility of its cross-line circulation.
\begin{proposition} (Time window-based search space reduction)\label{propo_domin_reduction}
    Constraints (\ref{cons_ve_timetable}) are valid for the TT-VS-DCA model (\ref{stochastic_problem}).
\begin{align}
        z^l_{k,i,t}=
        \begin{cases}
        0,\ \frac{\overline{\varsigma}^l_{k,i}}{\Delta}\leq t\leq \left|\mathcal{T}\right|\\
        1,\ 1\leq t\leq \frac{\underline{\varsigma}^l_{k,i}}{\Delta}-1 \label{cons_ve_timetable}
        \end{cases}, \qquad \forall l\in{\mathcal{L}}, k \in {\mathcal{K}}_l, i \in {\mathcal{S}}_l,
\end{align}
where $\underline{\varsigma}^l_{k,i}$ and $\overline{\varsigma}^l_{k,i}$ represent the earliest and latest departure times of the MAV assigned to trip $k$ leaving stop $i$ of line $l$, respectively. The specific values of $\underline{\varsigma}^l_{k,i}$ and $\overline{\varsigma}^l_{k,i}$ can be determined previously before optimization.
\end{proposition}

\textbf{Proof.}
See Appendix \ref{sec:proof}.

\begin{proposition}\label{propo_first_board}
    (Passenger boarding-related inequality) Constraints (\ref{cons_ve_boarding}) are a family of valid inequality for the TT-VS-DCA model related to the boarding process (\ref{stochastic_problem}).
    \begin{align}
& \overline{\varsigma}^l_{k,i}-u_p\geq \left|\mathcal{T}\right|\left(\xi^{l}_{p,k}(w)-\xi^{l}_{p,k-1}(w)-1 \right) & \forall w\in{\mathcal{W}}, p\in{\mathcal{P}_w}, l=l_p^1\in{\mathcal{L}}_p,k\in\mathcal{K}_l.\label{cons_ve_boarding}
\end{align}
\end{proposition}
\textbf{Proof.}
See Appendix \ref{sec:proof}.

\begin{proposition}\label{propo_transfer}
    (Transfer-related inequality) Constraints (\ref{cons_ve_transfer}) are a family of valid inequality for the TT-VS-DCA model (\ref{stochastic_problem}) associated with transfers.

\begin{align}
& \overline{\varsigma}^{l'}_{k',i'} - \underline{\varsigma}^{l}_{k,i} + \beta^{min}_{l,k,i} - \underline{\theta}_g \geq \left|\mathcal{T}\right|(\pi^{g}_{k,k'}-\pi^{g}_{k,k'-1} - 1) \quad \forall g = (l, l', i, i') \in \dot{\mathcal{A}}, k \in \mathcal{K}_l, k' \in \mathcal{K}_{l'}.
\label{cons_ve_transfer}
\end{align}
\end{proposition}
\textbf{Proof.}  
See Appendix \ref{sec:proof}.

\begin{proposition}\label{propo_coupling}
    (Coupling-related inequality) Constraints (\ref{cons_ve_couple}) are a family of valid inequality for the TT-VS-DCA model (\ref{stochastic_problem}).
\begin{align}
\begin{cases}
& \overline{\varsigma}^{l'}_{k',i'} - \underline{\varsigma}^{l}_{k,i} + \beta^{min}_{l,k,i} - \underline{\theta}_g \geq \left|\mathcal{T}\right|(e^{g}_{k,k'} - 1)\\
& \overline{\theta}_g - \underline{\varsigma}^{l'}_{k',i'} + \overline{\varsigma}^{l}_{k,i} - \beta^{max}_{l,k,i} \geq \left|\mathcal{T}\right|(e^{g}_{k,k'} - 1)
\end{cases}
\quad \forall g = (l, l', i, i') \in \dot{\mathcal{A}}, k \in \mathcal{K}_l, k' \in \mathcal{K}_{l'}.
\label{cons_ve_couple}
\end{align}
\end{proposition}
\textbf{Proof.}
See Appendix \ref{sec:proof}.

\section{Solution method}
\label{sec:algorithm}
In this section, we design the solution approaches to solve the TT-VS-DCA problem. First, we develop an integer L-shaped method with tailored valid inequalities in Section \ref{sec:decomposition}. However, through numerical experiments based on real-world operational data from the Beijing bus network with a 3-hour planning period, we observe that the proposed exact solution method still faces computational challenges. To address this, in Section \ref{sec:rolling}, we introduce a rolling-horizon optimization framework that dynamically incorporates the exact solution method, enabling a double decomposition of the problem. Finally, recognizing the need for real-time adaptability in practical operations with evolving passenger demand, we propose a learning-based real-time decision-making framework, as detailed in Section \ref{sec:real-time}.

\subsection{The integer L-shaped algorithm with tailored valid inequalities}
\label{sec:decomposition}
Since our proposed model (\ref{stochastic_problem}) is an MILP model, one way to find the optimal solution is to adopt the solution methods designed for linear programs, such as the built-in Branch-and-Bound or cutting-plane algorithms available in commercial softwares. Nevertheless, the quality of the obtained solution falls short of expectations when these methods are applied to a real-life instance based on a Beijing bus subnetwork that consists of two bidirectional lines, and the gap is greater than 42\% even after a 10-hour computation.

Another approach is to decompose the model using Benders decomposition method or the integer L-shaped method. Specifically, the integer L-shaped method is a powerful extension of the basic Benders decomposition approach, especially designed for addressing stochastic integer programming problems. In problems where there are binary first-stage variables and integer second-stage variables within a two-stage stochastic program, the integer L-shaped method has proven to be remarkably effective. Recall that our model has an endogenous sequential solution property, i.e., the optimization of vehicle scheduling and dynamic capacity allocation can only be achieved if timetables are determined. Therefore, decomposing our problem and solving it with the integer L-shaped algorithm promises to improve the solution efficiency and solution quality. 

We denote $\Psi_{wait}^1=\sum\limits_{w\in\mathcal{W}}\sum\limits_{p\in\mathcal{P}_w}n_p\cdot\rho_w\cdot \varepsilon_1\cdot\vartheta^T\cdot \Psi_p^1(w)$ and $\Psi_{wait}^2(w)=\sum\limits_{p\in\mathcal{P}_w} n_p\cdot\rho_w \cdot\vartheta^T\cdot\varepsilon_2\cdot \Psi_p^2(w)$. To implement the integer L-shaped method, we first reformulate the proposed model (\ref{stochastic_problem}) as follows:
\begin{eqnarray}\label{mp}
{\rm{[MP]}}\quad
\left\{
\begin{array}{ll}
&\min\limits_{\mathbf{z},\bm{\chi}} \  \Psi_{wait}^1 \\
&\mbox{s.t.}
\ \ \ 
\eqref{zbinary},
\eqref{binary_z} - \eqref{cons_pi}, \eqref{eq:chi}, \eqref{cons_ve_timetable}, \eqref{cons_ve_boarding} - \eqref{cons_ve_transfer},
\eqref{z_wait_linear},
\eqref{cons_indicator_transfer_linear},\\
& \quad \quad \ \ \eta\geq F(\mathbf{z},\bm{\chi},\bm{\zeta}),
\end{array}
\right.
\end{eqnarray}
where $F(\mathbf{z},\bm{\chi}, \bm{\zeta})$ is real-valued function of variables $\mathbf{z}$, $\bm{\chi}$, and $\bm{\zeta}$, which can be formulated as follows:
\begin{align}
    F(\mathbf{z},\bm{\chi},\bm{\zeta})&=\sum\limits_{w\in\mathcal{W}}\rho_w F^w(\mathbf{z},\bm{\chi},\bm{\zeta}) \nonumber\\ 
    &=\sum\limits_{w\in\mathcal{W}}\rho_w\bigg[\min\limits_{\mathbf{x},\mathbf{h}}\left\{\Psi_{wait}^2(w)+\Psi_{op}^w|\eqref{eq:dominH} - \eqref{eq:dominV},(\ref{operation_cost_z}), (\ref{cons_board})-(\ref{cons_number_mv}), \right. \nonumber\\
    &\phantom{=\sum\limits_{w\in\mathcal{W}}\rho_w\bigg[}\left.\qquad\qquad\qquad\qquad\quad(\ref{cons_he}) - (\ref{cons_xh_1}),\eqref{eq:dominY},(\ref{cons_ve_couple}), (\ref{z_transfer_linear}),(\ref{cons_vs_linear}) - (\ref{cons_hh_linear})\right\}\bigg].\label{sp}
\end{align}

In this study, we denote $F^w(\mathbf{z},\bm{\chi},\bm{\zeta})$ as the value function of the subproblem with respect to the scenario $w$. Given the timetables and the state of the passengers, the subproblem $F^w(\mathbf{z},\bm{\chi},\bm{\zeta})$ is to determine the optimal vehicle schedule and dynamic capacity allocations.

\begin{lemma}\label{lemma_lower_bound_1}
There exists a finite lower bound $L= \sum\limits_{l\in\mathcal{L}}\sum\limits_{i\in\mathcal{S}_l}\sum\limits_{k\in\mathcal{K}_l}\varphi_2\vartheta_{1,i}$ that satisfies $F(\mathbf{z},\bm{\chi},\bm{\zeta})\geq L$.
\end{lemma}
\textbf{Proof.}
See Appendix \ref{sec:proof}.

 As outlined in \cite{LAPORTE} and \cite{Angulo}, the core idea of the integer L-shaped method is to relax the constraint $\eta\geq F(\mathbf{z},\bm{\chi},\bm{\zeta})$. This approach entails adding cuts to more accurately approximate the shape of $F(\mathbf{z},\bm{\chi},\bm{\zeta})$. This iterative process continuously refines the solution space, providing an optimal solution that is incrementally tightly fitted. Specifically, the optimality cut proposed by \cite{LAPORTE} can be formulated for the explored problem as follows
\begin{align}\label{opt_cut}
\eta&\geq (F(\mathbf{z}^*,\boldsymbol{\chi}^*,\boldsymbol{\zeta}^*)-L)\bigg(\sum_{z_{k,i,t}^l\in\ H(z^*)}z_{k,i,t}^l+\sum_{\chi^l_{p,k}(w)\in\ G(\chi^*)}\chi^l_{p,k}(w)+\sum_{\zeta^g_{p,k,k'}(w)\in\ V(\zeta^*)}\zeta^g_{p,k,k'}(w)-\nonumber\\&\sum_{z_{k,i,t}^l\notin\ H(z^*)}z_{k,i,t}^l -\sum_{\chi^l_{p,k}(w)\notin\ G(\chi^*)}\chi^l_{p,k}(w)-\nonumber\\&\sum_{\zeta^g_{p,k,k'}(w)\notin\ V(\zeta^*)}\zeta^g_{p,k,k'}(w)-\left|H(z^*)\right|-\left|G(\chi^*)\right|-\left|V(\zeta^*)\right|\bigg)+F(\mathbf{z}^*,\boldsymbol{\chi}^*,\boldsymbol{\zeta}^*),
\end{align}
where $\eta$ represents the estimated value of $F(\mathbf{z},\boldsymbol{\chi},\bm{\zeta})$, and $H(z^*)=\{z^l_{k,i,t}|z^l_{k,i,t}=1\},G(\chi^*)=\{\chi^l_{p,k}(w)|\chi^l_{p,k}(w)=1\},V(\zeta^*)=\{\zeta^g_{p,k,k'}(w)|\zeta^g_{p,k,k'}(w)=1\}$. Moreover, $L$ represents the lower bound of $F(\mathbf{z},\boldsymbol{\chi},\bm{\zeta})$, whose values can be set based on Lemma \ref{lemma_lower_bound_1}. Since the aforementioned optimality cuts are not tight for solutions other than $(\mathbf{z}^*, \boldsymbol{\chi}^*,\bm{\zeta}^*)$, we employ the following continuous optimality cut the same as \cite{Angulo}:
\begin{align}\label{copt_cut}
\eta\geq \overline{\boldsymbol{\iota}}(\mathbf{z}-\mathbf{z}^*)+\hat{\boldsymbol{\iota}}(\boldsymbol{\chi}-\boldsymbol{\chi}^*)+\tilde{\boldsymbol{\iota}}(\boldsymbol{\zeta}-\boldsymbol{\zeta}^*)+F_{LP}(\mathbf{z}^*,\boldsymbol{\chi}^*,\boldsymbol{\zeta}^*),
\end{align}
where $\overline{\boldsymbol{\iota}}$, $\hat{\boldsymbol{\iota}}$ and $\tilde{\boldsymbol{\iota}}$ represents the subgradient vector of $F(\mathbf{z},\bm{\chi},\bm{\zeta})$ at $\mathbf{z}^*$,  $\boldsymbol{\chi}^*$, and $\boldsymbol{\zeta}^*$ respectively. $F_{LP}(\mathbf{z}^*,\boldsymbol{\chi}^*,\boldsymbol{\zeta}^*)$ represents the continuous relaxation of $F(\mathbf{z}^*,\boldsymbol{\chi}^*,\boldsymbol{\zeta}^*)$ that resets the binary variable $y^l_{k,i,q}(w)$ as a continuous variable. The subgradient optimality cut (\ref{copt_cut}) is not necessarily required for the convergence of the proposed method; however, it enhances the computing performance of this algorithm.

To ensure the feasibility of the solution, we introduced the following combinatorial cut
\begin{align}\label{feasibility_cut}
\sum_{z_{k,i,t}^l\in\ H(z^*)}(1-z_{k,i,t}^l)+\sum_{\chi^l_{p,k}(w)\in\ G(\chi^*)}(1-\chi^l_{p,k}(w))+\sum_{\zeta^g_{p,k,k'}(w)\in\ V(\zeta^*)}(1-\zeta^g_{p,k,k'}(w))\nonumber\\ +\sum_{z_{k,i,t}^l\notin\ H(z^*)}z_{k,i,t}^l+\sum_{\chi^l_{p,k}(w)\notin\ G(\chi^*)}\chi^l_{p,k}(w) + \sum_{\zeta^g_{p,k,k'}(w)\notin\ V(\zeta^*)}\zeta^g_{p,k,k'}(w)\geq 1.
\end{align}

Finally, we reformulate the MP based on the above analyses and construct the current problem (CP) at each node of the branch-and-cut search tree as 
\begin{equation}\label{rmp}
\left\{
\begin{aligned}
&\min_{z,\chi,\zeta,\eta}\Psi^{pass}+\eta \\
\mbox{s.t.}
&\text{ Optimality cuts } (\ref{opt_cut}),\\
&\text{ Subgradient optimality cuts } (\ref{copt_cut}),\\
&\text{ Feasibility cuts } (\ref{feasibility_cut}),\\
&\text{ Constraints } 
\eqref{zbinary}, \eqref{binary_z} - \eqref{cons_pi}, \eqref{eq:chi}, \eqref{cons_ve_timetable} - \eqref{cons_ve_transfer}, \eqref{z_wait_linear}, \eqref{cons_indicator_transfer_linear}.
\end{aligned}
\right.
\end{equation}

In summary, the main procedure of our designed decomposition-based integer L-shaped algorithm for the TT-VS-DCA problem can be presented as follows: (1) Relax CP (\ref{rmp}) by discarding constraints (\ref{zbinary}) and solving it to the optimum. It is worth noting that no optimality cut or subgradient cut has been generated yet. (2) Decide whether to proceed with the computation. If the stopping criteria are met, terminate the algorithm. Otherwise, we select a $z^l_{k,i,t}\notin\{0,1\}$ and then branch at this node. (3) Repeat the previous step until an integer solution $\mathbf{z}^*$ is obtained. In this process, it is necessary to determine whether to add optimality cuts (\ref{opt_cut}) and (\ref{copt_cut}) or feasibility cuts (\ref{feasibility_cut}).

In the specific implementation, we construct a branch-and-cut search tree using the state-of-the-art commercial solver GUROBI. The optimality cuts (\ref{opt_cut}) and (\ref{copt_cut}) as well as the feasibility cut (\ref{feasibility_cut}) are added to the CP step-by-step using the lazy constraint callback function. During the search process, the callback function is executed at the node where the relaxed CP is solved to the optimum, or a new MIP incumbent is found. The pseudocode of this algorithm is presented in Algorithm~\ref{table:algorithm} in Appendix \ref{sec:algorithmflow}. When the timetable is fixed, the vehicle schedule and dynamic-capacity allocation plan in each scenario can be optimized in parallel.

\subsection{A rolling-horizon optimization algorithm based on a double-decomposition framework}\label{sec:rolling}
The TT-VS-DCA problem is NP-hard, making the integer L-shaped method from Section \ref{sec:decomposition} inefficient for large-scale instances. A promising approach to address this challenge is to design a rolling-horizon optimization algorithm, which has been widely used to solve various traffic problems, as detailed in \cite{SANCHEZMARTINEZ20161}. This algorithm follows a progressive, periodically updated decision-making process, adding dynamism to the solution approach. 

The rolling-horizon optimization algorithm decomposes the TT-VS-DCA problem in both spatial and temporal dimensions, allowing the solution of a series of smaller MILPs instead of a single large MILP. Specifically, in the temporal dimension, the algorithm splits the entire study time horizon into multiple sequential stages, each treated as a separate subproblem. In the spatial dimension, it strategically considers only the stops involved within a specific stage during the optimization process. Furthermore, our rolling-horizon optimization algorithm is built upon a double-decomposition framework by incorporating the integer L-shaped method. It first decomposes the decision period into sequential stages with shorter time durations, and within each stage, the integer L-shaped method further decomposes the problem to improve computational efficiency.

Thereafter, we introduce this algorithm tailored to the TT-VS-DCA problem utilizing a rolling-horizon optimization framework in detail. The algorithm segments the studied time horizon into equal-length decision-making stages, each with two time domains: the operation control horizon (CH) and the information prediction horizon (PH). Both start simultaneously, but the PH, longer than the CH, allows for a broader scope of decisions and external data. Each PH intersects with CHs from adjacent stages, facilitating real-time updates. Within each phase, we optimize the TT-VS-DCA problem for the entire PH but employ solutions only within the CH. Information is sequentially transferred between stages, streamlining optimization and boosting algorithm efficiency. For clarity, we encapsulate the functions of the algorithm into the following modules:

\noindent\textbf{Initialization}: Input EPTN network topology, set of all passengers groups $\mathcal{P}_w$, control parameters of the system, model parameters. Initialize the time duration of each CH and PH. Also, the empty sets $\tilde{\mathcal{P}}$, and $\overline{\mathcal{P}}$ are initialized to store the set of passengers who do not complete transfers, and the set of passengers who do not get off MAVs, respectively. 

\noindent\textbf{Prediction modular}: We incorporate a prediction module into the algorithm to enable its extension to dynamic optimization. In this module, the travel time of vehicles between adjacent stops is predicted based on detailed road traffic conditions within a specific time horizon. Additionally, the module estimates passenger demand at each stop using data from monitoring equipment during the same period.

\noindent\textbf{Timetable simulation module}: Update the travel time of the vehicle between two adjacent stops. Simulate the vehicle arrival and departure times at each stop based on the control strategy imposed during the CH.

\noindent\textbf{Select module}: The timetable simulation module is first called to simulate the arrival and departure times of trips. Then both the MAVs assigned to trips and stops that are within the PH are selected. Denote the selection vehicles set as $\mathcal{K}^n\in\mathcal{K}$, and denote the selection stops set as $\mathcal{S}^n\in\mathcal{S}$. For each scenario $w$, we select all passengers whose arrival time is earlier than the ending time of PH from the set $\mathcal{P}_w$ and form the new set $\mathcal{P}_w^n$. Then, the model of TT-VS-DCA prblem for the current decision stage is created.

\noindent\textbf{Status updating module}: First, we check the status of each passenger group in each scenario $w$. Specifically, if the passenger group has completed all trips, then the group from $\mathcal{P}_w$ will be removed; if the passenger group has boarded but has not completed all transfers, then we place it in set $\tilde{\mathcal{P}}_w$; if the passenger has completed all transfers but has not yet alighted from the last leg of the trip, then we place it in $\overline{\mathcal{P}}_w$. In particular, for passengers that require placement in $\tilde{\mathcal{P}}_w$, we identify their completed legs, record the set of all vehicles they have boarded, and then remove the completed trips from the passenger group attribute. In the next optimization, by identifying the current state of passengers, the variables describing passengers’ state (i.e. $\xi_{p,k}^l(w)$, $\zeta_{p,k,k'}^{g}(w)$, $\chi^l_{p,k}(w)$) can be precisely constructed to complete the information transfer between different optimization time domains. Then, the number of MAVs arriving and alighting at each stop, the number of MAVs in operations, the formation of each operated MAVs, and the number of MAVs available at each depot and time are all recorded. Lastly, we update the time range covered by CH and PH.

\noindent\textbf{Optimization module}: Input $\mathcal{K}^n, \mathcal{S}^n, \mathcal{P}^n_w, \tilde{\mathcal{P}}_w$, and $\overline{\mathcal{P}}_w$. For passenger groups in $\overline{\mathcal{P}}_w$, since they have already completed their first boarding, we do not construct variable $\chi^l_{p,k}(w)$ to track their waiting time. However, we still need $\overline{\mathcal{P}}_w$ to accurately calculate the number of people alighting from the vehicle at each stop. For the passenger in $\tilde{\mathcal{P}}_w$, we selectively construct variable $\zeta_{p,k,k'}^{g}(w)$ to track the transfer possibly completed at that stage. For the passengers in $\mathcal{P}^n_w$, we leave out any special treatment and construct the corresponding variables according to the model requirements. Then, we run Algorithm 1 in Online Appendix G and output the optimization results, i.e., timetables, vehicle schedules, and the dynamic capacity allocation strategy in the current CH.

To sum up, the procedure of the overall algorithm is presented in Figure \ref{fig:RH_new} in Appendix \ref{sec:flowChat}. This procedure begins by introducing the initial system settings. The next step is to check that the start of the current iteration and the end of the entire optimization period coincide. If they align, the process is stopped and the optimization results are obtained. Otherwise, the prediction module is started as a precautionary measure to update the travel time between subsequent stops and the passenger demand at each stop. Subsequently, the timetable simulation module is utilized to obtain a partially regulated timetable including all the CH cycle control strategies from the previous sequence. Next, the selection module is run to select the stops and MAVs assigned to trips involved in this iteration. After that, the optimization module is run to get the current CH control strategy and the vehicle schedules. Further, we run the state update module to accurately identify the states of passengers and vehicles, and update the corresponding information which is the input for the next optimization iteration.

\subsection{Learning-based real-time decision-making framework}
\label{sec:real-time}

As passenger flows evolve over time, timetables and vehicle schedules derived from historical demand scenarios may become suboptimal when faced with real-time operational conditions. Moreover, as discussed in Section~\ref{sec:introduction}, MAVs are operated under a short decision-making time window. To address this and better respond to time-varying demand that may fall outside the historical passenger demand scenarios, we propose a real-time decision-making framework that fine-tunes the tactical timetable within a small adjustment range using real-time information and re-optimizes the vehicle schedule accordingly. The framework follows an RH approach, partitioning the studied time horizon into multiple decision-making stages with equal-length short durations. At each stage, a stochastic MILP is solved using the integer L-shaped method, where the timetable obtained from tactical-level planning is allowed to shift within an adjustment range (e.g., 1 or 2 minutes) and the vehicle schedules can be re-optimized.

In stochastic programming, incorporating more scenarios can improve solution quality but greatly increases computational difficulty (as will be shown in Section~\ref{stochastic_experiment_subsection}). However, real-time operations require finding solutions efficiently while maintaining high quality. To address this challenge, we design a learning-based scenario-retention method and embed it into our real-time decision-making framework. Our proposed learning-based real-time decision-making framework is illustrated in Figure~\ref{fig:Real_time} in Appendix~\ref{sec:Overflow}. At each decision-making stage: (i) the real-time information is updated; (ii) using predefined features, a trained machine learning (ML) model is used to rank all new demand scenarios that are out of the historical demand scenario set; (iii) we retain the top-$k$ most representative scenarios into a scenario subset; (iv) we solve the stochastic MILP using the integer L-shaped method over this scenario subset to optimize deviations of the tactical timetable and the vehicle schedule. This process is repeated throughout the entire studied horizon, ensuring that timetables and vehicle schedules in all decision-making stages remain consistent with real-time dynamics. By avoiding incorporating all demand scenarios into the integer L-shaped method at each stage and focusing computation over the most representative scenarios, the framework meets real-time computational requirements while preserving solution quality.

We now introduce our proposed learning-based scenario-retention method in detail. Table~\ref{tab:feature_definitions} summarizes the features used to characterize each demand scenario, which explicitly capture interactions between supply and demand. Intuitively, one might assume that scenarios with higher total demand are more critical. However, this criterion can be misleading. Consider the following illustrative example: suppose there are two stations. In scenario 1, the total demand is 50, evenly distributed with 25 passengers at each station. In scenario 2, the total demand is slightly lower at 45, but highly imbalanced with 40 passengers at the first station and only 5 at the second. A timetable which is optimized solely based on scenario 1 may fail to provide sufficient capacity for scenario 2, resulting in large passengers' costs or even infeasibility. This example highlights that the criticality of a scenario depends not only on the aggregate demand level, but also on the distribution of demand, the alignment of supply and demand, and their interactions. These relationships are highly nonlinear and difficult to capture analytically in complex and uncertain operating environments. Motivated by the strength of ML methods in capturing nonlinear patterns, we adopt an ML-based approach to learn the mapping between scenario features and their impact on system performance, thereby identifying the most critical scenarios for real-time optimization.

\begin{table}[htbp]
\centering
\begin{threeparttable}
\caption{Feature definitions and formulations used in the learning-based scenario-retention method.}
\label{tab:feature_definitions} 
\begin{tabular}{@{} c c c @{}}
\toprule
Feature name & Description & Formula \\
\midrule
\texttt{total\_waiting} & Total passengers waiting in the system &
$\sum\limits_{p \in \mathcal{P}_w} n_p$ \\
\texttt{max\_waiting} & Maximum passengers waiting at a single stop &
$\max\limits_{i \in \mathcal{S}} \left\{ \text{demand}_i\right\}$ \\
\texttt{avg\_gap} & Average of supply-demand gap across stops &
$\frac{1}{|\mathcal{S}|}\sum\limits_{i\in\mathcal{S}}\text{diff}_i$ \\
\texttt{worst\_gap} & Most negative (worst-case) supply-demand gap &
$\min\limits_{i \in \mathcal{S}} \left\{\text{diff}_i\right\}$ \\
\texttt{gap\_std} & Standard deviation of supply-demand gap &
$\sqrt{\frac{1}{|\mathcal{S}|}\sum_{i}(\text{diff}_i - \overline{\text{diff}_i})^2}$ \\
\texttt{demand\_std} & Standard deviation of demand across stops &
$\sqrt{\frac{1}{|\mathcal{S}|}\sum_{i}(\text{demand}_i - \overline{\text{demand}_i})^2}$ \\
\texttt{supply\_std} & Standard deviation of supply across stops &
$\sqrt{\frac{1}{|\mathcal{S}|}\sum_{i}(\text{supply}_i - \overline{\text{supply}_i})^2}$ \\
\bottomrule
\end{tabular}

\begin{tablenotes}
\footnotesize
\item \textit{Notes.} $\text{demand}_i =\sum\limits_{p \in \mathcal{PO}_i} n_p$ is the demand at stop $i$. $\text{supply}_i= \text{CAP}\cdot \hat{x}_{k,i}^l(w) - \hat{v}_{k,i}^l(w)$ is the supply at stop $i$. $\text{diff}_i= \text{supply}_i - \text{demand}_i$ is the supply-demand gap at stop $i$. $\overline{\text{diff}_i}=\frac{1}{|\mathcal{S}|} \sum_{i\in\mathcal{S}}\text{diff}_i$ is the average supply-demand gap across stops.  $\overline{\text{demand}_i}=\frac{1}{|\mathcal{S}|} \sum_{i\in\mathcal{S}}\text{demand}_i$ is the average demand across stops. $\overline{\text{supply}_i}=\frac{1}{|\mathcal{S}|} \sum_{i\in\mathcal{S}}\text{supply}_i$ is the average supply across stops.
 \end{tablenotes} 
\end{threeparttable}
\end{table}

In the learning-based scenario-retention method, we perform an \emph{offline training} procedure using historical data. In the offline training framework, passenger flow fluctuations are first simulated by perturbing station-level demand within historical demand scenarios, thereby generating a rich set of stochastic passenger flow scenarios that reflect possible real-time variations. Given that our proposed optimization method follows a RH framework (see Section~\ref{sec:rolling}), the training data generation process is embedded within this RH structure to ensure compatibility between the learning model and the underlying optimization procedure. Thereafter, at each decision-making stage, we randomly select a subset of scenarios and run the integer L-shaped method with this subset. At the beginning of each decision-making stage, we simulate the current operational plan and estimate station-level supply (i.e., available capacity). We then extract scenario-level features as defined in Table~\ref{tab:feature_definitions}. Next, we randomly choose one scenario from the selected subset as the input of the optimization procedure, solve the problem under this scenario, and evaluate the solution’s performance across the other selected scenarios. The resulting \emph{performance gap} between the pre- and post-optimization objectives is used as the ground-truth label for the supervised learning. To calibrate the learning model, we adopt the XGBoost algorithm. The main reasons are three-fold: (i) XGBoost captures complex nonlinear relationships between scenario features and solution quality; (ii) it supports efficient training with large-scale datasets and delivers high predictive performance; and (iii) it provides interpretable feature importance metrics, helping us understand the influence of each feature on scenario criticality.

\section{Numerical experiments}\label{sec:experiments}

In this section, we conduct a comprehensive set of numerical experiments to evaluate the performance and practical value of the proposed models and solution methods. The experiments are organized as follows. We begin by describing the test instances used throughout the experiments in Section~\ref{sec:instances}. We then analyze the benefits of incorporating dynamic capacity allocations and cross-line circulations in Section~\ref{sec:flexible_capacity}, and evaluate the impact of prioritizing passenger interests on both operational costs and the proportion of in-vehicle transfers in Section~\ref{sec:passenger_impact}. Next, we evaluate the benefits of our proposed integrated optimization method compared to a sequential approach (Section~\ref{sec:integrated_vs_sequential}), followed by a comparison of the performance of various solution methods, including our exact algorithm, rolling-horizon framework, and GUROBI (Section~\ref{sec:solution_comparison}).

To evaluate the robustness of our methods, we conduct two post hoc analyses: one testing the timetable and vehicle schedule and doing a sensitivity analysis of overload allowance in Section~\ref{sec:robustness_schedule}, and another testing the effectiveness of fine-tuning strategies embedded in our learning-based real-time decision-making framework in Section~\ref{sec:robustness_tuning}. We also investigate the value of the stochastic programming approach compared to deterministic models (Section~\ref{stochastic_experiment_subsection}). Finally, we evaluate the effectiveness of our learning-based real-time decision-making framework in Section~\ref{sec:real_time}, and analyze the trade-offs between operational and passengers' costs in Section~\ref{sec:tradeoff_analysis}. All experiments were implemented in Python and executed on a Windows 11 personal computer equipped with a 12th Gen Intel(R) Core(TM) i7-12700H processor and 64GB RAM.

\subsection{Instances}
\label{sec:instances}
The instances used for the experiments are derived from a virtual network and a part of the Beijing bus network, which are depicted in Figure \ref{fig:experimentNetwork} in  Appendix \ref{sec:instance}. The virtual network comprises three bidirectional lines. The Beijing bus subnetwork contains two bidirectional lines and 89 stops, where passenger demand, network structure, and running times between stops are all derived from real-world operating data. Based on these two networks, we generate 24 sets of instances with varying durations of time periods and numbers of scenarios to test the effectiveness of the proposed formulation and algorithms. Detailed information on each instance is provided in Table \ref{tab:instance_all} in Appendix \ref{sec:instance}. In addition, the five parameters used in the objective function \eqref{eq:objective_brief} are set as follows. In all experiments, the equivalent monetary value per unit of passengers' waiting costs, $\vartheta^T$, is set to 0.8(\$). The coefficients $\varepsilon_1$ and $\varepsilon_2$, which represent the weights on waiting time at origins and due to transfers, are set to 1 and 1.5, respectively. The coefficients $\varphi_1$ and $\varphi_2$, representing the weights on passengers' and operational costs, are both set to 1 in the experiments from Sections~\ref{sec:flexible_capacity} to \ref{sec:real_time}. A sensitivity analysis of these two coefficients is provided in Section~\ref{sec:tradeoff_analysis}.

\subsection{Benefits of incorporating dynamic capacity allocations and cross-line circulations}
\label{sec:flexible_capacity}

To gain insights into the effects of dynamic capacity allocations and cross-line circulations on the passengers' and operational costs, we perform a series of experiments. We define three operational strategies: (1) \textit{Fixed-capacity}, which refers to executing all trips using MAVs at their maximum capacity to ensure service quality; (2) \textit{Partially flexible-capacity}, which permits MAVs to be coupled or decoupled only at depots, maintains consistent en-route formations, and allows MAUs to be dispatched to different lines at depots (i.e., depot-based cross-line circulations); and (3) \textit{Completely flexible-capacity}, which allows MAVs to be flexibly coupled or decoupled at depots and transfer stops, as well as cross-line circulations of MAUs.

In Table \ref{tab:resultStrategy}, we present the results of adopting the aforementioned three operational strategies. We can observe that compared to the fixed-capacity operational strategy, both the partially and completely flexible-capacity strategies significantly reduce operational costs and the number of used MAUs. For example, for $instance\_rp\_2$, the fixed-capacity strategy requires 72 MAUs, while partial flexible-capacity and flexible-capacity only require 48 and 36 MAUs, resulting in reductions of 33.3\% and 50.0\%, respectively. We also observe that the completely flexible-capacity strategy slightly reduces total passengers' costs compared to the fixed-capacity strategy. This is mainly because the flexible-capacity strategy allows in-vehicle transfers, which do not incur transfer waiting costs that are included in the passengers' cost term of the objective function. In contrast, conventional transfers under the fixed-capacity strategy involve inconvenience, leading to positive transfer waiting costs. A third observation is that the completely flexible-capacity strategy outperforms the partially flexible-capacity with respect to the number of used units among all instances. The largest relative savings compared to the partially flexible-capacity operational approach occur in the $instance\_rp\_3$, where the completely flexible-capacity strategy saves over 28.1\% of units. In conclusion, the completely flexible-capacity operational strategy is beneficial to operators in terms of lowering operational costs and the number of required MAUs, thanks to its ability to dynamically allocate capacity at transfer stops and schedule cross-line MAU circulations during tirps.

\begin{table}[h]
\centering
\caption{The total passengers' costs, operational costs, and the number of used MAUs obtained by adopting various operational strategies.}
\label{tab:resultStrategy}
\begin{threeparttable}
\begin{tabular}{llrrr}
\toprule
Instance & Operational strategy & \begin{tabular}[c]{@{}l@{}}Total passengers'\\ costs (\$)\end{tabular} & \begin{tabular}[c]{@{}l@{}}Total operational\\ costs (\$) \end{tabular}& \begin{tabular}[c]{@{}l@{}}The number of \\ used MAUs\end{tabular} \\
\toprule
{instance\_rp\_1} & Fixed-capacity & 1541.8 & 198.9  & 48 \\
 & Partially flexible-capacity & 1545.7  & 115.0  & 25 \\
 & Completely flexible-capacity & 1539.9  & \textbf{100.2}  & \textbf{21} \\
  [1ex]
{instance\_rp\_2} & Fixed-capacity & 4106.3  & 300.4  & 72 \\
 & Partially flexible-capacity & 4122.4  & 197.4  & 48 \\
 & Completely flexible-capacity & 4047.1  & \textbf{163.3}  & \textbf{35} \\
   [1ex]
{instance\_rp\_3} & Fixed-capacity & 7042.6  & 403.3  & 96  \\
 & Partially flexible-capacity & 7064.5  & 283.0 & 71 \\
 & Completely flexible-capacity & 6981.1   & \textbf{227.8}   & \textbf{51} \\
  [1ex]
{instance\_rp\_4} & Fixed-capacity & 10413.0  & 506.2  & 96 \\
 & Partially flexible-capacity & 10459.2  & 372.3  & 79 \\
 & Completely flexible-capacity & 10472.1  & \textbf{310.5}  & \textbf{58} \\
[1ex]
{instance\_rp\_5} & Fixed-capacity & 15614.7  & 610.4  & 96 \\
 & Partially flexible-capacity & 15682.7  & 461.3  & 84 \\
 & Completely flexible-capacity & 15666.8  & \textbf{380.2}  & \textbf{61} \\[1ex]
{instance\_rp\_6} & Fixed-capacity & 19507.2   & 718.8  & 96 \\
 & Partially flexible-capacity & 19578.0   & 549.3   & 85 \\
 & Completely flexible-capacity & 19555.0   & \textbf{458.6}   & \textbf{69} \\
\bottomrule
\end{tabular}
\end{threeparttable}
\end{table}

\begin{figure}[h]
     \centering
     \begin{subfigure}{0.49\textwidth}
         \centering
         \includegraphics[width=\textwidth]{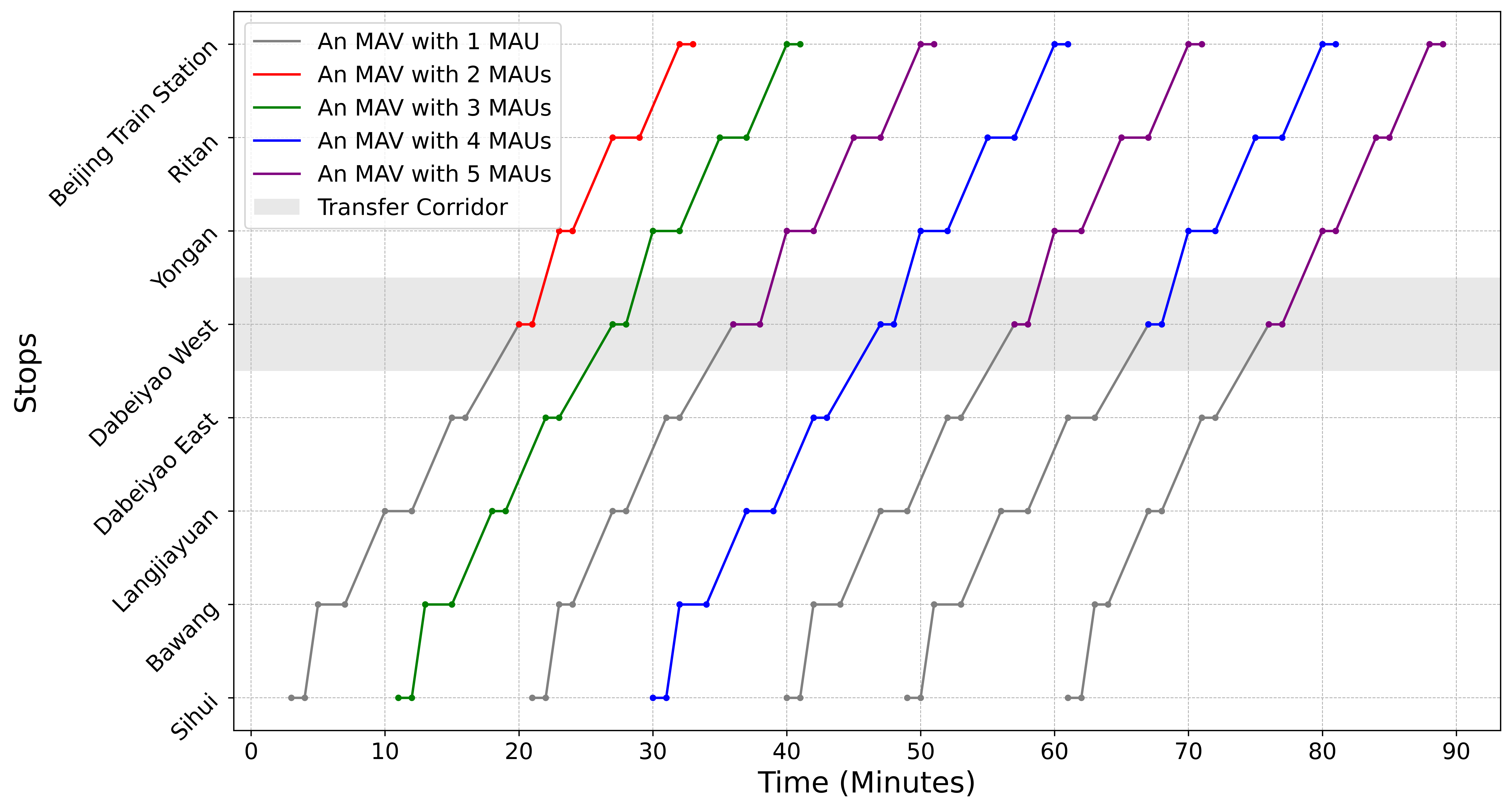}
         \caption{The timetable and formations of MAVs}
     \end{subfigure}
     \begin{subfigure}{0.5\textwidth}
         \centering
         \includegraphics[width=1\textwidth]{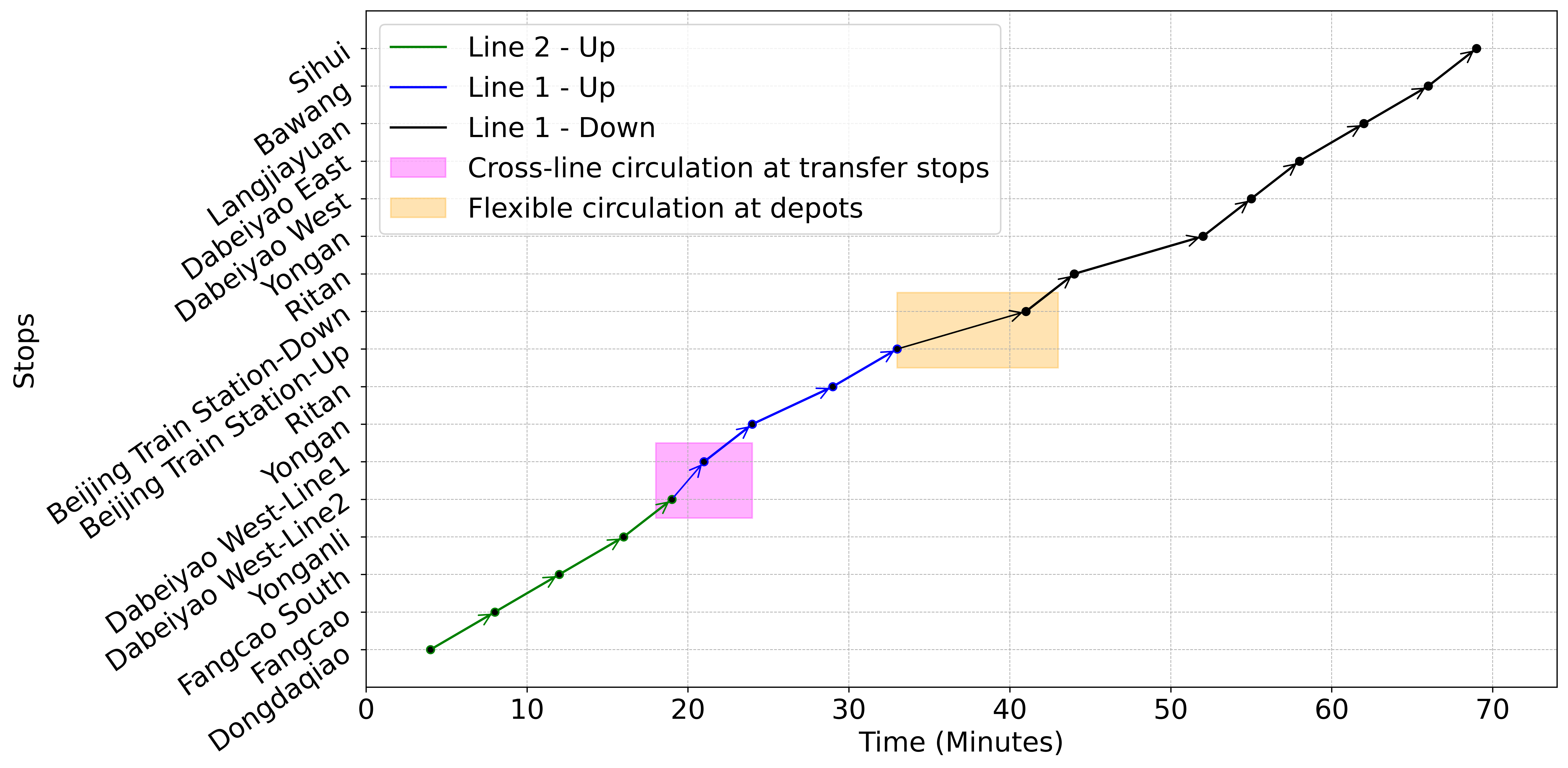}
         \caption{Movements of an MAU}
     \end{subfigure}
        \caption{The optimized timetable, formations of MAVs assigned to trips, and cross-line schedules of MAUs of the instance$\_$rp$\_$6. Note: Shaded area indicates the transfer stop.}
        \label{fig:resultTTVS}
\end{figure}

Figure \ref{fig:resultTTVS} delves into the solution for the $instance\_rp\_6$, where Figure \ref{fig:resultTTVS}(a) presents the optimal robust timetable and formations of MAVs designated to each trip for the east-west line, termed as Line 1, using the flexible-capacity operational strategy. Figure \ref{fig:resultTTVS}(b) displays the complete trajectory of an unit departing from the first stop in the upstream direction of the north-south line, denoted as Line 2. We can observe that the formation of the MAV assigned to the first trip is one unit before reaching the transfer stop Dabeiyao West due to low passenger demand. After the transfer stop, because of a marked surge in demand, an unit is insufficient. A vacant unit on Line 2 is coupled to the MAV assigned to this trip on Line 1, as shown in Figure \ref{fig:resultTTVS}(b). After completing the trip and return to the depot, this MAV is reallocated to a trip in the opposite direction. To sum up, the timetable, formations of MAVs at different stops and times, and cross-line circulations of MAUs are optimized according to the time-varying and unevenly distributed demand. 

\subsection{Impact of valuing the interests of passengers}
\label{sec:passenger_impact}
We next perform a series of
experiments base on $instance\_g\_1$ and $instance\_g\_2$ to assess the impact of valuing interests of passengers on passengers’ and operational costs, as well as the proportion of in-vehicle transfers among all transfers. Firstly, we set the weighting coefficient related to passengers' costs in the objective function, $\varphi_1$, to 0 and the weight associated with operational costs $\varphi_2$ as 1 to obtain the solution with the minimal operational costs. Thereafter, we incrementally increase the value of $\varphi_1$ to explore the effects of valuing passenger interests.

\begin{figure}[h]
\centering
\subfloat{
\begin{minipage}[]{0.5\linewidth}
 \centerline{\includegraphics[width=\textwidth]{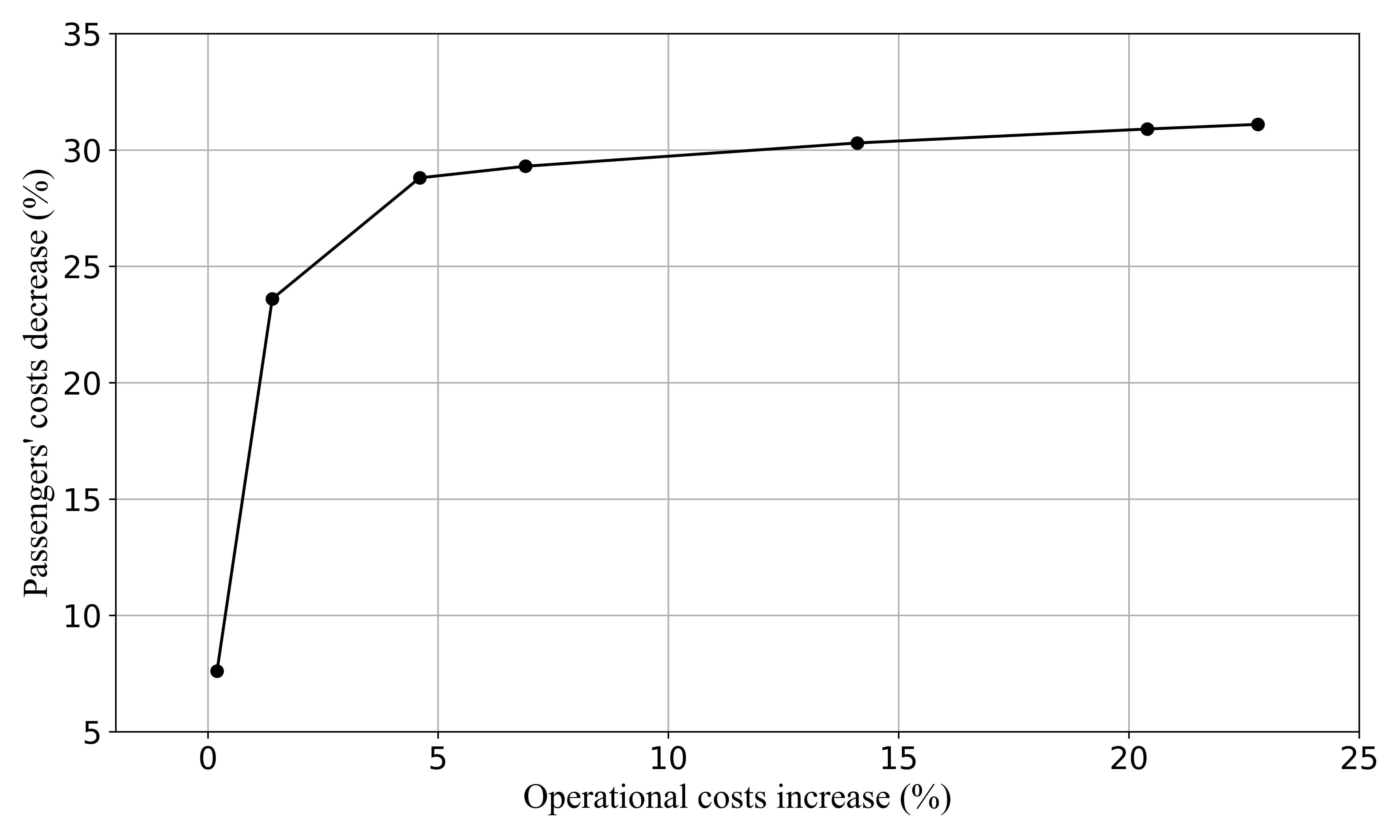}}
 \centerline{(a) $instance\_g\_1$}
\end{minipage}
}
\subfloat{
\begin{minipage}[]{0.5\linewidth}
 \centerline{\includegraphics[width=\textwidth]{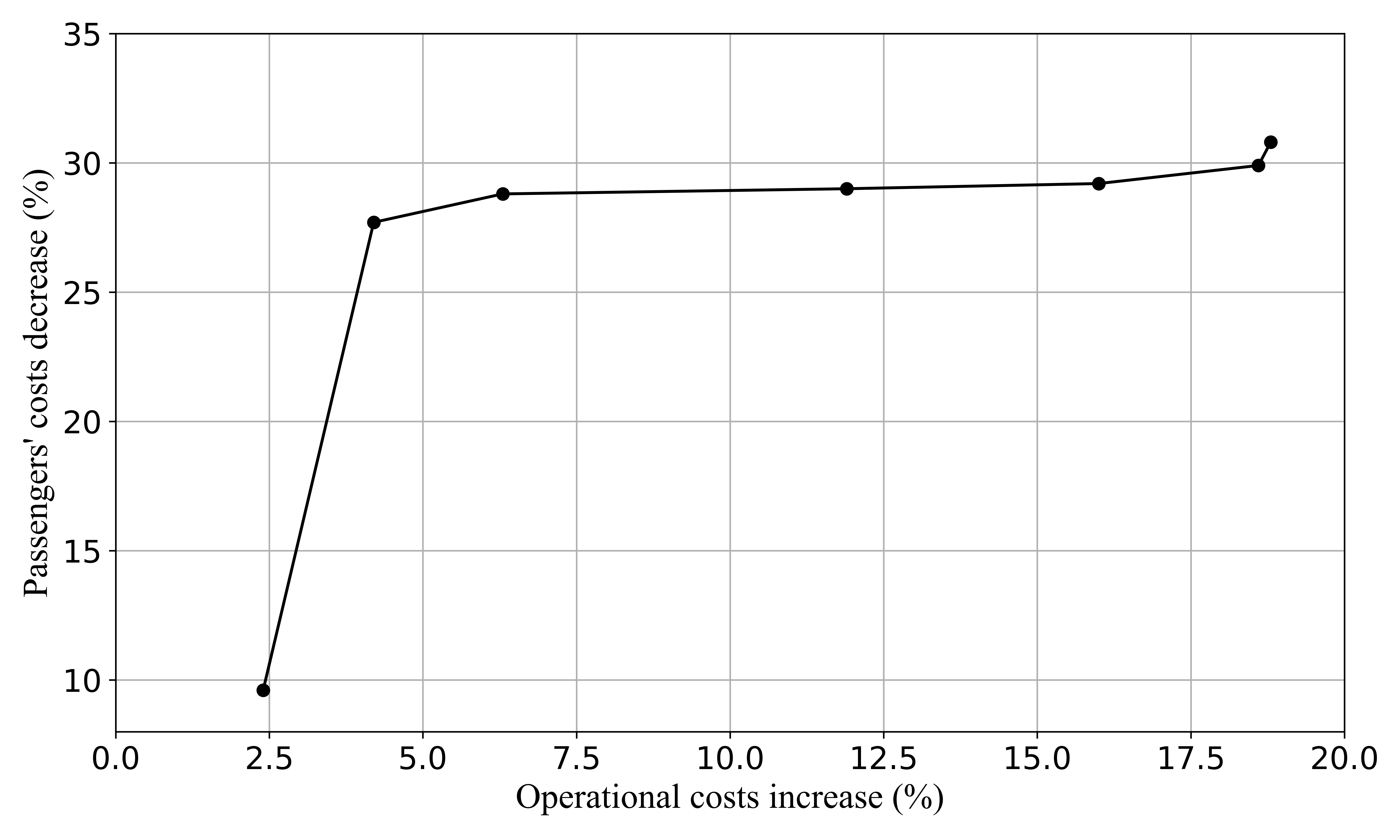}}
 \centerline{(b) $instance\_g\_2$}
\end{minipage}
}
\caption{Results of passengers’ and operational costs when the weighting coefficient related to passengers' costs in the objective function increases.}\label{fig:operational_costs_passengers_costs}
\end{figure}

As the value of $\varphi_1$ increases, the operational costs grow and passengers' costs decrease. Figure \ref{fig:operational_costs_passengers_costs} illustrates the increase in operational costs and decrease in passengers' costs for all Pareto-efficient solutions, relative to the solution with minimal operational costs. It can be observed that an increase in operational costs by approximately 5\% correlates with a reduction in total passengers' costs by about 30\% in both instances. However, as operational costs continue to increase, the subsequent reduction in total passengers' costs becomes marginal. These findings suggest while high values only have a little positive impact, a reasonable attention in the interests of passengers can significantly increase the quality of the service.

\begin{figure}[h]
\centering
\subfloat{
\begin{minipage}[]{0.5\linewidth}
 \centerline{\includegraphics[width=\textwidth]{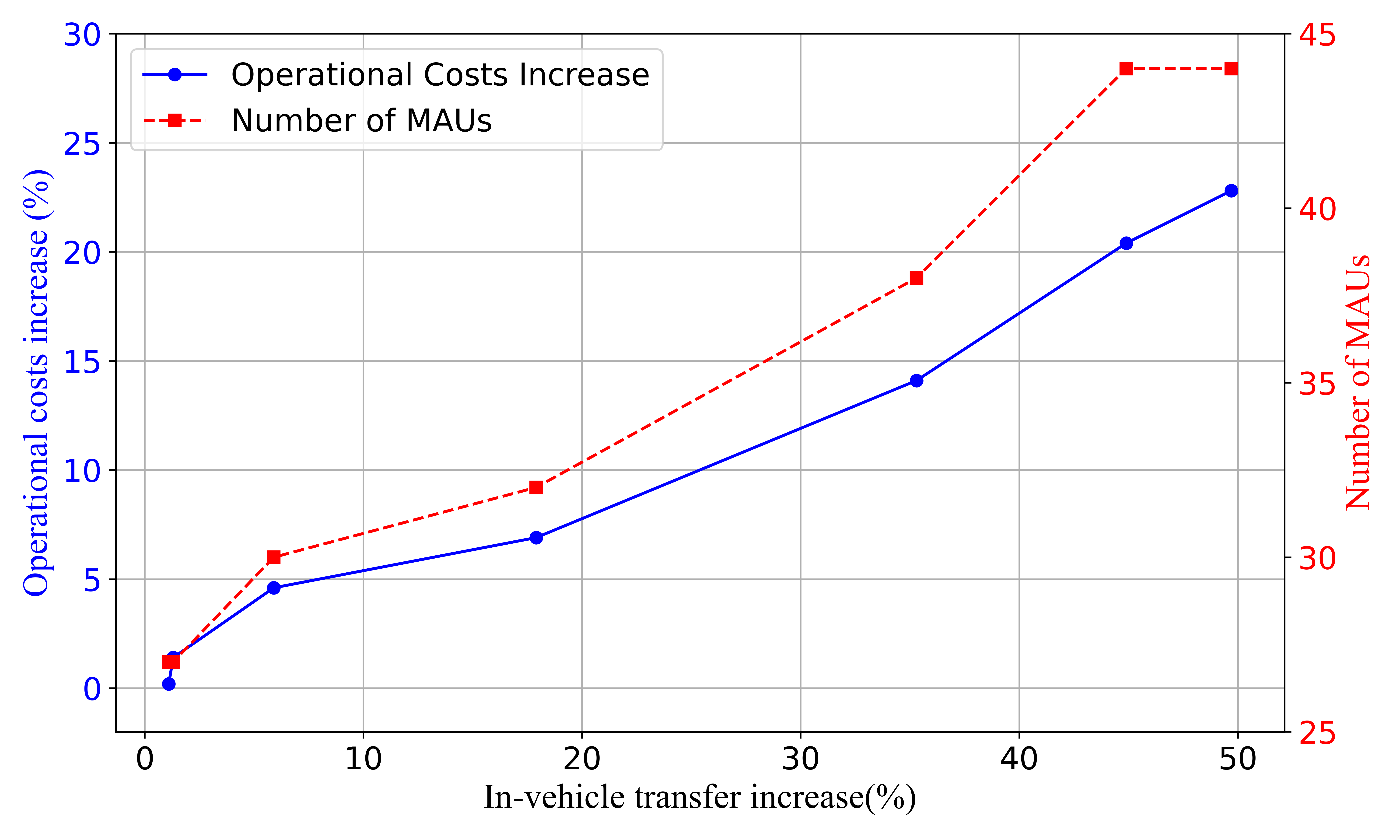}}
 \centerline{(a) $instance\_g\_1$}
\end{minipage}
}
\subfloat{
\begin{minipage}[]{0.5\linewidth}
 \centerline{\includegraphics[width=\textwidth]{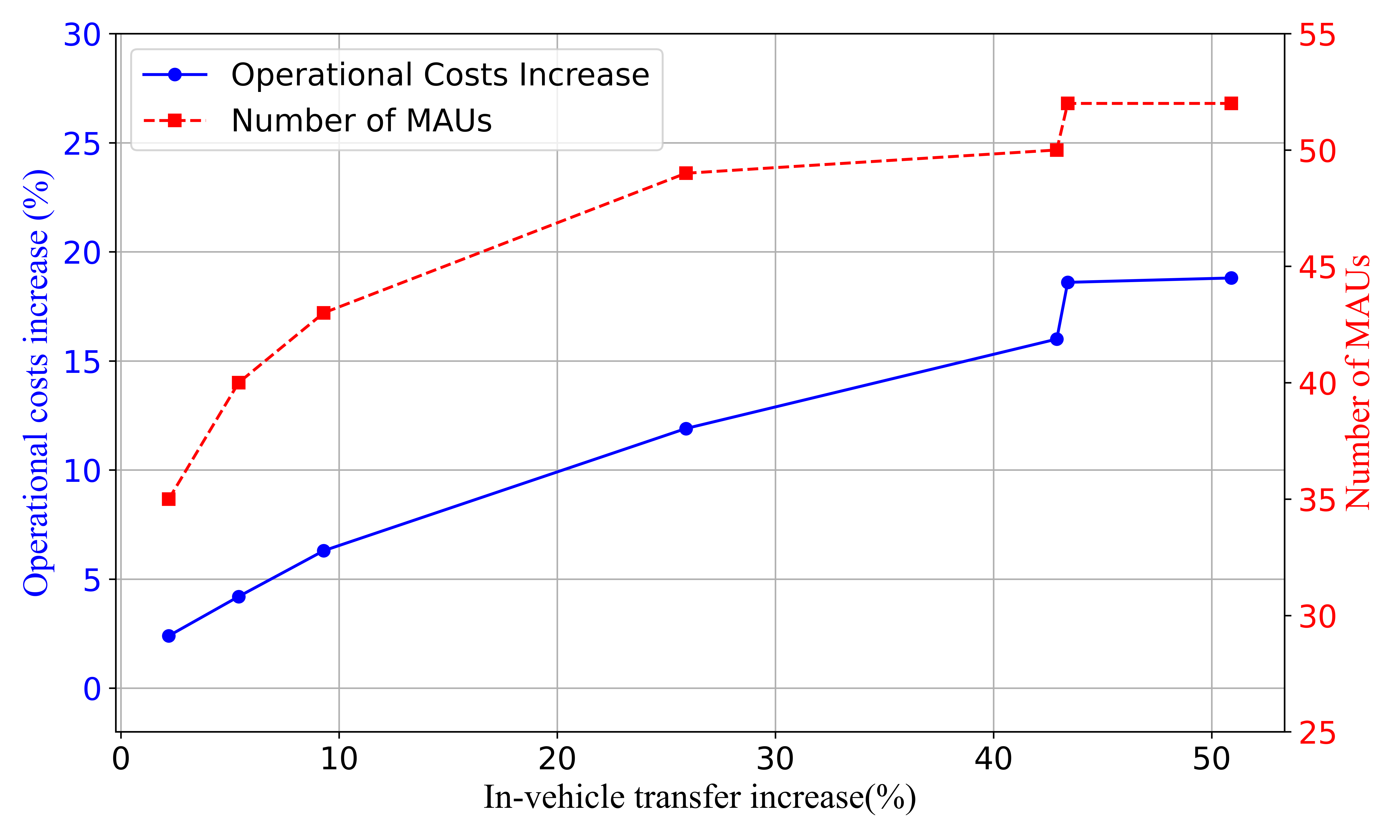}}
 \centerline{(b)  $instance\_g\_2$}
\end{minipage}
}
\caption{Results of the proportion of in-vehicle transfers among all transfers and the used MAUs when the weighting coefficient related to passengers' costs in the objective function increases.}\label{fig:in_vehicle_transfer}
\end{figure}

Another interesting aspect is the impact of prioritizing passengers' interests on the proportion of in-vehicle transfers among all passengers' transfers in the solution and the number of required MAUs. The obtained solutions are presented in Figure \ref{fig:in_vehicle_transfer}. In both instances, there is a clear trend indicating that greater proportions of in-vehicle transfers are closely associated with higher operational costs and more utilized MAUs. Specifically, as the proportion of in-vehicle transfers increases from 0\% to 50\%, there is a consistent and significant increase in operational costs — up to 25\% in $instance\_g\_1$ and 20\% in $instance\_g\_2$. Furthermore, in $instance\_g\_1$ and $instance\_g\_2$, the number of used MAUs expands from 26 to 45, and from 34 to 52, respectively, with relative percentage increases of almost 73\% and 53\%.

This trend can be explained by the role of transfer waiting time in passenger costs. In-vehicle transfers are not penalized, whereas conventional transfers incur additional penalties due to the inconvenience of alighting, walking between platforms, and re-boarding. These penalties are explicitly incorporated into passengers’ costs through weighted transfer waiting time, as defined in constraints (\ref{z_transfer}). As greater weight is placed on passenger interests, the model naturally favors in-vehicle transfers to reduce these costs. In summary, prioritizing passengers' interest increases the share of in-vehicle transfers and enhances transfer convenience, but this comes at the expense of higher operational costs and greater MAU usage.

\subsection{Benefits of the proposed integrated optimization method}
\label{sec:integrated_vs_sequential}
This section assesses the benefits of integrating timetabling and vehicle scheduling
decisions. We construct a set of experiments based on six instances from $instance\_rp\_1$ to $instance\_rp\_6$, where we use the sequential decision-making approach as a benchmark. Specifically, we propose a two-step sequential method that firstly determine the timetable. Then, the derived optimal timetable is treated as an input to optimize the vehicle schedule. We use GUROBI to solve our integrated optimization model~\eqref{stochastic_problem} and the sequential optimization formulation.

In Table~\ref{tab:solution_performance}, we present the objective function values, passengers' costs, relative deviations in passenger costs (denoted as $Dev^p$), and the number of utilized MAUs obtained from the two optimization methods. It can be observed that: (i) For smaller-scale instance, such as $instance\_rp\_1$, these two methods achieve same optimal solutions; however, as the scale of the instances slightly increases (e.g., $instance\_rp\_2$ and $instance\_rp\_3$), the sequential optimization method attains suboptimal solutions with a higher objective value compared to our proposed integrated optimization method. (ii) When the scale of instances is further enlarged, i.e., $instance\_rp\_4$, $instance\_rp\_5$ and $instance\_rp\_6$, the sequential method fails to find feasible solutions, which indicates that a given number of units cannot fully execute the timetable obtained from the optimization in the first step. In contrast, the proposed integrated optimization method is able to find the optimal timetable and corresponding vehicle schedule for each instance. In these larger instances, the number of used MAUs is 58, 61, and 69, respectively. (iii) When the sequential optimization method yields feasible solutions, the passenger costs obtained by the integrated optimization method are slightly higher, with deviations up to 0.12\%. These observations highlight the necessity for incorporating an integrated optimization approach to efficiently plan timetables and vehicle schedules in large-scale networks.

\begin{table}[h]
    \centering
    \caption{Performance comparison between the integrated and sequential timetabling and vehicle scheduling for MAVs.}
    \label{tab:solution_performance}
\resizebox{\textwidth}{!}{%
    \begin{tabular}{ccrrrccrrrrr}
    \toprule
        Instance & \begin{tabular}[c]{@{}c@{}}Optimization\\ method\end{tabular} & \begin{tabular}[c]{@{}c@{}}Objective\\ value (\$)\end{tabular}  & \begin{tabular}[c]{@{}c@{}}Passengers'\\ cost (\$)\end{tabular} & $Dev^p  (\%)$& \begin{tabular}[c]{@{}c@{}}The number of\\ used MAUs\end{tabular} & Instance & \begin{tabular}[c]{@{}c@{}}Solution\\ type\end{tabular} & \begin{tabular}[c]{@{}c@{}}Objective\\ value (\$)\end{tabular} & \begin{tabular}[c]{@{}c@{}}Passengers'\\ cost (\$)\end{tabular} & $Dev^p  (\%)$ &\begin{tabular}[c]{@{}c@{}}The number of\\ used MAUs\end{tabular} \\
    \midrule
    instance\_rp\_1 & Sequential & 3544.32 &2121.12 &-  & 21 & instance\_rp\_4 & Sequential & - & - & - \\
     & Integrated & 3544.32 &2121.12 & 0& 21 & & Integrated & 16682.1 &12513.17 & - & 58 \\
    \addlinespace[0.5em]
    instance\_rp\_2 & Sequential & 7349.31 &5049.79 & - & 36 & instance\_rp\_5 & Sequential &- & - &- & - \\
    & Integrated & 7313.58 &5053.34 & 0.07 & 35 &  & Integrated & 23270.8 &18783.85 & - & 61 \\
    \addlinespace[0.5em]
    instance\_rp\_3  & Sequential & 11633.17 &8791.77 & - & 59 & instance\_rp\_6 & Sequential &- & - &- & - \\
    & Integrated & 11538.13 & 8802.32 & 0.12 & 51 &  & Integrated & 28727.0 & 22727.0 & -& 69 \\
 \bottomrule
\end{tabular}%
}
\end{table}

On the other hand, sequential optimization has its own advantages. Specifically, the sequential method is advantageous in scalability and modeling flexibility, notably by allowing for the flexible modeling of over-saturation scenarios, where passengers unable to board the first available MAV may wait for subsequent ones. Although over-saturation did not occur in our computational results, it remains a possibility in practice. In contrast, the integrated optimization method excels in providing high-quality solutions, especially concerning fleet sizing, due to the inherent feedback between vehicle scheduling and timetabling. Nevertheless, incorporating these considerations into the integrated optimization framework significantly increases model complexity and computational burden, making it impractical at large scales. Therefore, the proposed RH optimization approach is motivated and useful, effectively balancing model complexity and scalability while retaining the benefits of the integrated optimization framework.

\subsection{Performance comparison among various solution methods}
\label{sec:solution_comparison}

To assess the effectiveness of the proposed solution methods, we perform the following two experiments. Firstly, we construct 12 instances based on the virtual network to compare the performance between GUROBI and the integer L-shaped method. Secondly, based on the Beijing bus subnetwork, we build six instances to evaluate the efficiency and solution quality of both the integer L-shaped and the rolling horizon solution methods by varying the control and predict horizons.

\subsubsection{Performance comparison between GUROBI and the integer L-shaped method. } 

In this experiment, we compare the performance between methods for solving the TT-VS-DCA model: a direct approach utilizing GUROBI, and the integer L-shaped method (denoted as IL method) with valid inequalities. Our termination criteria include: (i) achieving an optimality gap tolerance of 5\%, and (ii) reaching a 2-hour computational time limit.

In Table \ref{tab:algorithm_comparison}, we present the results of 12 instances, stating the objective value, computational time and the relative optimality gap. It can be observed that that the proposed IL method obtains solutions with smaller objective function values in all instances with respect to the solutions obtained by GUROBI. While GUROBI demonstrates smaller optimality gaps in smaller instances, such as $instance\_g\_1$ and $instance\_g\_7$, it requires longer computational times. As the size of instances increases, however, the superiority of the IL method over GUROBI becomes more evident. For example, when sloving $instance\_g\_5$ and $instance\_g\_10$, GUROBI struggles to find a good-quality solution within 2 hours. Although the IL method cannot achieve an optimality gap of less than 5\% for $instance\_g\_6$, $instance\_g\_11$, and $instance\_g\_12$ within two hours, it produces solutions with lower objective values than those generated by GUROBI in the same computational time. Overall, our proposed exact solution method has a significant advantage over GUROBI when addressing medium-scale problems.

\begin{table}[ht]
\centering
\caption{Performance comparison between GUROBI and the designed exact solution algorithm.}
\resizebox{\textwidth}{!}{%
\begin{tabular}{c c r r r c c r r r}
\toprule
Instance & \begin{tabular}[c]{@{}c@{}}Solution\\method\end{tabular} & \begin{tabular}[c]{@{}c@{}}Objective\\value (\$)\end{tabular} &\begin{tabular}[c]{@{}c@{}}Computational \\time (s) \end{tabular} & \begin{tabular}[c]{@{}c@{}}Gap\\ (\%) \end{tabular}  & Instance & \begin{tabular}[c]{@{}c@{}}Solution\\method\end{tabular} & \begin{tabular}[c]{@{}c@{}}Objective\\value (\$)\end{tabular} & \begin{tabular}[c]{@{}c@{}}Computational \\ time (s) \end{tabular}  &\begin{tabular}[c]{@{}c@{}}Gap\\ (\%) \end{tabular}  \\
\midrule
\multirow{2}{*}{instance\_g\_1} & GUROBI & 3075.2 & 452.7 & 4.3 & \multirow{2}{*}{instance\_g\_7} & GUROBI & 2727.8 & 668.3 & 1.7 \\
 & IL & 2990.2 & \textbf{344.1} & 4.8 &  & IL & 2708.8 & \textbf{564.4} & 4.4 \\\addlinespace[0.5em]
\multirow{2}{*}{instance\_g\_2} & GUROBI & 4985.6 & 1065.7 & 2.3 & \multirow{2}{*}{instance\_g\_8} & GUROBI & 4636.0 & 3046.2 & 4.7 \\
 & IL & 4976.8 & \textbf{588.5} & 4.7 &  & IL & 4523.3 & \textbf{756.8} & 5.0 \\\addlinespace[0.5em]
\multirow{2}{*}{instance\_g\_3} & GUROBI & 7254.3 & 2681.3 & 4.6 & \multirow{2}{*}{instance\_g\_9} & GUROBI & 6477.7 & 6799.3 & 3.7 \\
 & IL & 7105.1 & \textbf{1186.7} & 5.0 &  & IL & 6388.7 & \textbf{1670.9} & 5.0 \\\addlinespace[0.5em]
\multirow{2}{*}{instance\_g\_4} & GUROBI & 9663.5 & 5021.3 & 4.9 & \multirow{2}{*}{instance\_g\_10} & GUROBI & 11149.1 & 7200.0 & 25.7 \\
 & IL & 9450.6 & \textbf{1917.3} & 4.9 &  & IL & 8526.1 & \textbf{5453.4} & \textbf{4.7} \\\addlinespace[0.5em]
\multirow{2}{*}{instance\_g\_5} & GUROBI & 12639.9 & 7200.0 & 8.8 & \multirow{2}{*}{instance\_g\_11} & GUROBI & 15134.3 & 7200.0 & 30.9 \\
 & IL & 11915.6 & \textbf{4602.8} & \textbf{5.0} &  & IL & 10819.7 & 7200.0 & \textbf{5.2} \\\addlinespace[0.5em]
\multirow{2}{*}{instance\_g\_6} & GUROBI & 19023.6 & 7200.0 & 24.5 & \multirow{2}{*}{instance\_g\_12} & GUROBI & 18615.1 & 7200.0 & 30.7 \\
 & IL & 15035.9 & 7200.0 & \textbf{7.9} &  & IL & 13568.2 & 7200.0 & \textbf{7.3} \\
\bottomrule
\end{tabular}%
}
\label{tab:algorithm_comparison}
\end{table}

\subsubsection{Performance comparison between the integer L-shaped method and the rolling horizon algorithm under various settings.}\label{RH_IL}

In this experiment, we compare the integer L-shaped method and the rolling horizon algorithm combined with this exact solution method (denoted as RH+IL). We consider six different scales of instances based on the Beijing bus subnetwork. For each instance, we vary the lengths of both control horizon (CH) and predict horizon (PH) to explore the impact of increasing the length of each optimization span on the quality of the final solution. It is worth mentioning that the RH+IL algorithm is equivalent to the IL algorithm when the CH length equals the length of the entire investigated planning period. Besides, we use the formula $Opt. \ Gap = \frac{\text{(Objective value - Best-known objective value)}}{\text{(Objective value)}} \times 100 \ (\%) $ to calculate the optimality gap for the rolling horizon (RH) algorithm.

The results are presented in Table \ref{tab:resultAlgorithm}. We can obtain the following observations: (1) For small instances, the setting with $CH{=}PH$ spanning the entire studied time horizon attains the best or comparable gaps while requiring the highest computational time. However, as instance size grows, this setting becomes the worst performer within the time limit (e.g., 7{,}200\,s), yielding the largest optimality gaps due to the much larger problem size. (2) The quality of solutions generated by the RH+IL algorithm improves as the duration of CH and PH increase, especially when the scale of instances is relatively small. (3) For large-scale instances, such as $instance\_r\_5$ and $instance\_r\_6$, obtaining high-quality solutions within a reasonable computational time using only the integer-L shaped algorithm proves challenging. (4) Even with a shorter optimization period (i.e., CH = 20, PH = 40), the gap of the solution obtained by RH+IL with respect to the optimal solution (maximum of not more than 7\% in all the instances) is acceptable; the gap of the resulting solution is further reduced to not more than 5\% after increasing CH to 80 and PH to 100. This observation indicates the high quality of the solution obtained in our proposed rolling horizon  algorithm. If the real-time passenger demand can be inputted and optimized online using the algorithm proposed in this paper, the quality of the solution can be guaranteed even with a small optimization period, while its solution speed can also meet the needs of real-time control.

\begin{table}[]
\centering
\caption{The objective value, CPU time and the $Opt. \ Gap$ that are obtained by RH+IL algorithm under different lengths of CH and PH.}
\label{tab:resultAlgorithm}
\begin{threeparttable}
\resizebox{\textwidth}{!}{%
\begin{tabular}{lrlrrr}
\toprule
Instances & Clairvoyant$^{a}$ (\$) & Solution method & \begin{tabular}[c]{@{}c@{}}Objective \\ value (\$)\end{tabular} & \begin{tabular}[c]{@{}c@{}}Computational \\ time (s)\end{tabular} & $Opt. \ Gap$ (\%) \\ \toprule
{instance\_r\_1} & {958.54} &    RH with $CH^{b} = 20 \ min$ and $PH^{b} = 40 \ min$  &    967.9    &  59.7 & 1.0    \\ 
 &    &    RH with $CH = 40 \ min$ and $PH = 60 \ min$  &   966.7   & 71.7  &    0.9 \\ 
  &    &    RH with $CH = 80 \ min$ and $PH = 100 \ min$  &  963.7    & 356.4  &0.5    \\ 
  &    &    RH with $CH = 120 \ min$ and $PH = 120 \ min$  &  963.4   &  977.1 &  0.5  \\ 
  \cline{3-6}
  \addlinespace
  {instance\_r\_2} & {1965.7}  &    RH with $CH = 20 \ min$ and $PH = 40 \ min$  & 2011.5    &  227.4 &2.3\\
 &    &    RH with $CH = 40 \ min$ and $PH = 60 \ min$  &  2052.3   &  516.9 &4.4    \\ 
 &    &    RH with $CH = 80 \ min$ and $PH = 100 \ min$  & 1996.4    &  1800.1 &1.6   \\
  &    &    RH with $CH = 130 \ min$ and $PH = 130 \ min$  & 1983.1    &  7200 & 0.9    \\ \cline{3-6}
    \addlinespace
 {instance\_r\_3} & {3055.2} &    RH with $CH = 20 \ min$ and $PH = 40 \ min$  & 3172.7    &  317.8 & 3.9\\
 &    &    RH with $CH = 40 \ min$ and $PH = 60 \ min$  &  3187.7   &  1745.0 &  4.3   \\ 
  &    &    RH with $CH = 80 \ min$ and $PH = 100 \ min$  & 3176.9    &  3605.1 & 4.0    \\ 
    &    &    RH with $CH = 140 \ min$ and $PH = 140 \ min$  & 3161.1    &  7200.0 & 3.5    \\ \cline{3-6}
      \addlinespace
 {instance\_r\_4} & {5130.5} &    RH with $CH = 20 \ min$ and $PH = 40 \ min$  & 5261.5    &  1379.6 & 2.6 \\
 &    &    RH with $CH = 40 \ min$ and $PH = 60 \ min$  &  5264.8   &  3105.5 &2.6   \\ 
  &    &    RH with $CH = 80 \ min$ and $PH = 100 \ min$  & 5247.2    &  7210.0 & 2.3    \\
  &    &    RH with $CH = 150 \ min$ and $PH = 150 \ min$  & 5636.6    &  7200.0 & 9.9    \\  \cline{3-6}
        \addlinespace
   {instance\_r\_5} & {5934.7} &    RH with $CH = 20 \ min$ and $PH = 40 \ min$  & 6353.2    &  1551.0 & 7.0 \\
 &    &    RH with $CH = 40 \ min$ and $PH = 60 \ min$  &  6239.0   &  3332.0 &5.1   \\ 
  &    &    RH with $CH = 80 \ min$ and $PH = 100 \ min$  & 6131.6    &  7211.4 & 3.3    \\
  &    &    RH with $CH = 150 \ min$ and $PH = 150 \ min$  & 6688.2    &  7200.0 & 12.7    \\  \cline{3-6}
    \addlinespace
  {instance\_r\_6} & {8622.4}&    RH with $CH = 20 \ min$ and $PH = 40 \ min$  & 9158.2     &  1687.3 &6.2  \\
 &    &    RH with $CH = 40 \ min$ and $PH = 60 \ min$  & 9322.0      &1978.0    &8.1   \\ 
  &    &    RH with $CH = 80 \ min$ and $PH = 100 \ min$  & 9033.4     &  3614.1  &4.7     \\
  &    &    RH with $CH = 240 \ min$ and $PH = 240 \ min$  & 9677.7     & 7200.0 & 12.2     \\  
 \bottomrule
\end{tabular}%
}
         \begin{tablenotes}
        \footnotesize
    \item[a] The clairvoyant is obtained by GUROBI, with a termination condition of 168 hours of computation time or a gap value of \\ less than 0.1\%.
    \item[b] $CH$ refers to the control horizon of RH. $PH$ refers to the prediction horizon of RH.
      \end{tablenotes}
\end{threeparttable}
\end{table}

\subsection{Post hoc robustness evaluation and sensitivity analysis of overload allowance}
\label{sec:robustness_schedule}

In this section, we assess the robustness of the optimized timetable and vehicle schedule on out-of-sample demand scenarios and evaluate the potential of capacity relaxation. We first conduct a post hoc evaluation using new demand scenarios to examine how the solution performs under both demand decreases and increases. We then perform a sensitivity analysis by gradually relaxing the capacity constraints through an overload allowance parameter $\lambda$, to investigate how small relaxations can help absorb demand surges without requiring additional resources.

\subsubsection{Out-of-sample robustness under demand perturbations.}

We assess the robustness of the optimized timetable and vehicle schedule obtained from our proposed integrated optimization method through post hoc evaluation using new demand realizations. Based on the original demand scenario set of instance\_rp\_6 used during optimization, we generate five new stochastic scenario sets by reducing total passenger demand by 20\%, 15\%, 10\%, and 5\%, and increasing it by 5\%, respectively.

To evaluate the number of overloaded passengers, we relax the capacity constraints \eqref{cons_capacity} and instead incorporate an overload penalty term $\sum\limits_{w \in \mathcal{W}} \rho_w \sum\limits_{l \in \mathcal{L}} \sum\limits_{k \in \mathcal{K}_l} \sum\limits_{i \in \mathcal{S}_l} \max\left(0, v_{k,i}^l(w) - \text{CAP} \cdot x_{k,i}^l(w) \right)$ into the objective function, where $v_{k,i}^l(w)$ denotes the number of in-vehicle passengers, $\text{CAP} \cdot x_{k,i}^l(w)$ represents the available capacity of the MAV operating section $(i, i+1)$ on trip $k$ and line $l$, and $\rho_w$ is the probability of scenario $w$. Since overloading is not allowed during robustness evaluation, these overloaded passengers are interpreted as unserved. We then input the optimized timetable and vehicle schedule of instance\_rp\_6 into Formulation~\eqref{stochastic_problem} with the modified objective function while without the original capacity constraints.

Table~\ref{tab:robustness_eval} summarizes the robustness of the optimized solution for instance\_rp\_6 under five new demand scenario sets. The first column indicates the tested demand scenario set, the second column reports the expected number of overloaded passengers (i.e., unserved passengers), and the third column reports the percentage of overloaded situations. This percentage is calculated as the number of overloaded MAV-section traversals divided by the total number of such traversals, where each MAV-section traversal refers to an MAV operating over a specific section. From the results shown in Table~\ref{tab:robustness_eval}, we observe that under scenario sets with decreased passenger demand, the solution remains feasible with no unserved passengers. When the demand increases by 5\%, 0.72\% of MAV-section traversals become overloaded, leading to an expected number of 147 overloaded (i.e., unserved) passengers. 
\begin{table}[h]
\centering
\caption{Post hoc robustness evaluation of the optimized timetable and vehicle schedule under new demand scenario sets.}
\label{tab:robustness_eval}
\resizebox{\textwidth}{!}{%
\begin{tabular}{rrr}
\toprule
Demand scenario set & Expected number of overloaded (unserved) passengers & Overloaded MAV-section traversals (\%) \\
\midrule
20\% decrease & 0   & 0 \\
15\% decrease & 0   & 0 \\
10\% decrease & 0   & 0 \\
5\% decrease  & 0   & 0 \\
5\% increase  & 147.00 & 0.72 \\
\bottomrule
\end{tabular}%
}
\end{table}

\subsubsection{Sensitivity analysis to overload allowance.}

To address the potential unserved issues, one possible method is to relax the capacity limitation and allow overloading. Therefore, we conduct a sensitivity analysis on overload allowances under the new stochastic scenario set with a 5\% increase in passenger demand. To do so, we relax the original capacity constraints \eqref{cons_capacity} by introducing a soft overload threshold. Specifically, we define an overload factor $\lambda \geq 1$, which allows a certain level of overloading. Recall that $v_{k,i}^l(w)$ denote the number of in-vehicle passengers when the MAV assigned to trip $k \in \mathcal{K}_l$ departs from stop $i \in \mathcal{S}_l$ on line $l \in \mathcal{L}$ in scenario $w \in \mathcal{W}$. Let $\text{CAP} \cdot x_{k,i}^l(w)$ denote the nominal capacity of the MAV with the formation of $x_{k,i}^l(w)$ assigned to that trip. We allow $v_{k,i}^l(w)$ to exceed the nominal capacity by up to a factor of $\lambda$. The modified capacity constraints can be expressed as  
\begin{equation}\label{con_newcapacity}
    v_{k,i}^l(w) \leq \lambda \cdot \text{CAP} \cdot x_{k,i}^l(w), \quad \forall w \in \mathcal{W},\, l \in \mathcal{L},\, k \in \mathcal{K}_l,\, i \in \mathcal{S}_l.
\end{equation}
Here, if $\lambda =1$, the capacity is the nominal capacity without allowing overloading. 

\begin{figure}[h]
    \centering
    \includegraphics[width=0.8\linewidth]{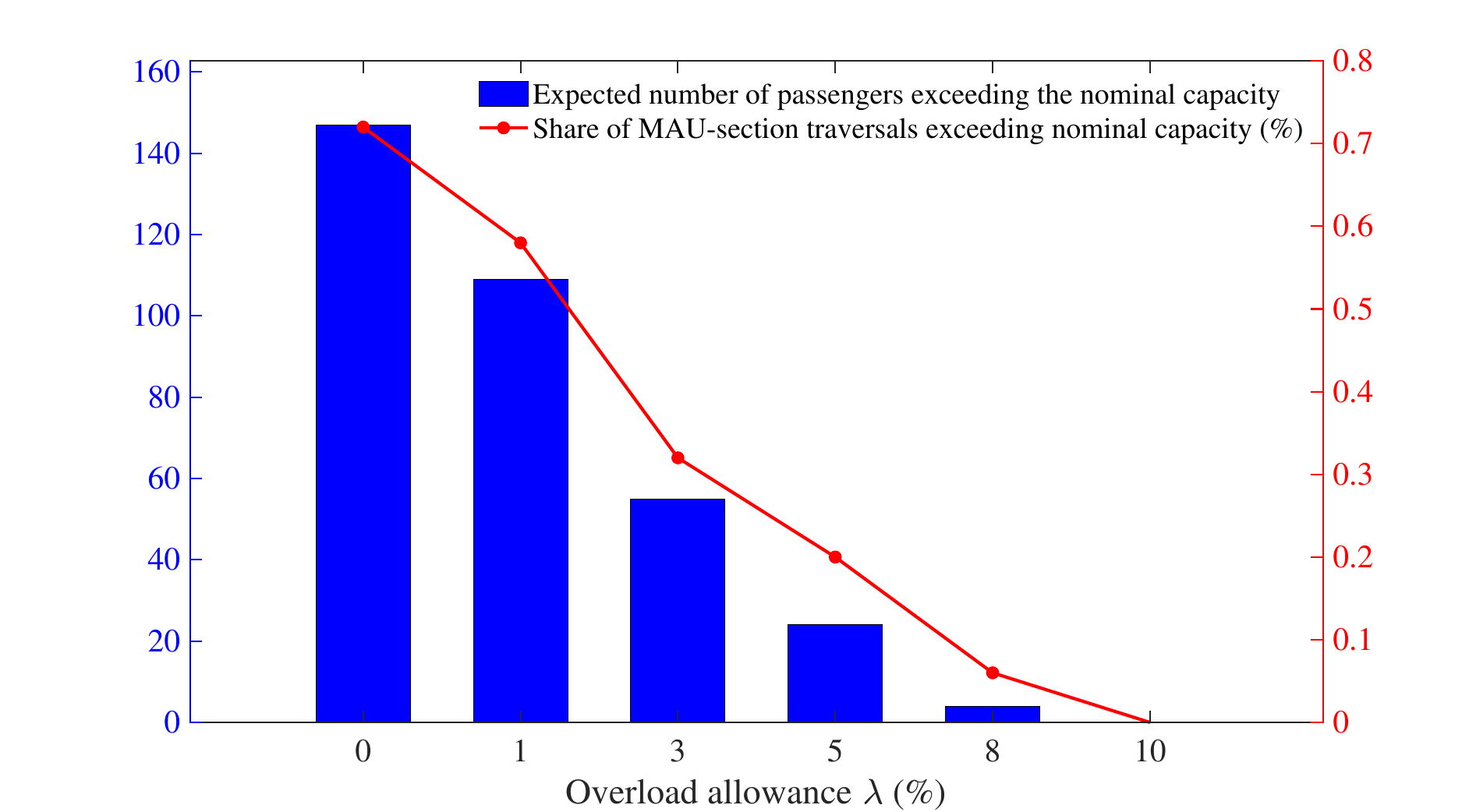}
\caption{Effect of overload allowance $\lambda$ on violations of the nominal capacity.}
    \label{fig:capacity}
\end{figure}

Figure~\ref{fig:capacity} shows the effect of the overload allowance $\lambda$ on two metrics, both evaluated against the nominal capacity (i.e., $\lambda \cdot \mathrm{CAP} \cdot x_{k,i}^l(w)$ with $\lambda = 1$): (i) the expected number of passengers exceeding the nominal capacity, i.e., i.e., $\sum\limits_{w \in \mathcal{W}} \rho_w \sum\limits_{l \in \mathcal{L}} \sum\limits_{k \in \mathcal{K}_l} \sum\limits_{i \in \mathcal{S}_l} \max\left(0, v_{k,i}^l(w) - \text{CAP} \cdot x_{k,i}^l(w) \right)$, and (ii) the proportion of overloaded MAU-section traversals where the number of in-vehicle passengers exceeds the nominal capacity (\%). The results demonstrate a clear monotonic improvement. Even modest overload allowances significantly reduce capacity violations. As $\lambda$ increases from 0 to 5\%, the expected number of overloaded passengers drops from 147 to 24, while the proportion of overloaded traversals declines from 0.72\% to 0.20\%. At $\lambda=10\%$, both indicators reach zero. These findings suggest that a 10\% overload allowance offers sufficient flexibility to absorb a 5\% demand surge without requiring additional capacity over the nominal capacity. Conversely, if strict adherence to nominal capacity is enforced, the results can guide the required capacity increase for each MAU to fully serve all passengers under such demand variations, based on the optimized timetable and vehicle schedule obtained from our proposed integrated optimization method.

\subsection{Post hoc robustness evaluation of fine-tuning strategies under demand surges and re-distributions}
\label{sec:robustness_tuning}
We now design two fine-tuning strategies to address capacity shortages under demand surges and conduct a post hoc robustness evaluation of them, aiming to provide practical and resilient solutions for real-world operations. The first strategy, Fix Timetable and Adjust Vehicle Schedule (AVS), retains the optimized timetable (denoted OT) generated from our integrated optimization method based on the in-sample demand scenario set, and re-optimizes only the vehicle scheduling decisions. The second strategy, Fine-tune Timetable and Adjust Vehicle Schedule (FTTAVS), uses OT as the baseline but permits each departure time to shift within a narrow range of ±1 minute, allowing joint re-optimization of the timetable and vehicle schedule. 

To evaluate these two strategies, we construct 48 new stochastic scenario sets based on instance\_rp\_6. Specifically, for each of the four levels of demand increase (5\%, 10\%, 15\%, and 20\%), we apply three distinct re-distribution strategies and generate 12 scenario sets. (i) The Uniform (U) strategy proportionally increases demand across all OD pairs, reflecting a general system-wide growth. (ii) The Head-biased (HB) strategy concentrates the increase on the top 10\% busiest OD pairs, simulating scenarios where already the busiest corridors become even more popular. (iii) In contrast, the Tail-biased (TB) strategy amplifies demand in the bottom 10\% least-used OD pairs, capturing emergent or previously underutilized travel patterns. These profiles allow us to test robustness under a broad range of practically relevant demand surges and spatial re-distributions. We test the two strategies across these 48 out-of-sample scenario sets. Each instance is denoted as D$k$\_{\text{profile}\_\text{strategy}}, where $k$ represents the level of demand increase, \enquote{profile} refers to the re-distribution strategy (U, HB, or TB), and \enquote{strategy} indicates the fine-tuning method (AVS or FTTAVS). For example, D10\_{\text{HB}}\_{\text{FTTAVS}} refers to the scenario with a 10\% head-biased demand surge evaluated using the FTTAVS strategy.

Table~\ref{tab:robust_metrics} summarizes the key robustness metrics across all demand surge and re-distribution scenarios. The second column reports the objective value, including passengers' and operational costs. The third column presents the expected number of unserved passengers across all scenarios--that is, passengers assigned to overloaded MAVs under the soft capacity limitation, which are penalized heavily and thus regarded as unserved. The fourth column, \textit{Top-5\% Avg. load ratio (\%)}, measures the average utilization of the most crowded 5\% of sections, where utilization $r^l_{k,i}(w)$ is computed as $r^l_{k,i}(w) = v^l_{k,i}(w) / (x^l_{k,i}(w)\, \text{CAP})$, with $v^l_{k,i}(w)$ denoting the number of in-vehicle passengers when the MAV assigned to trip $k$ leaves stop $i$ on line $l$ in scenario $w$, $x^l_{k,i}(w)$ representing the number of MAUs composing the MAV assigned to the section between stops $i$ and $i+1$, and \text{CAP} indicating the capacity per MAU. The fifth column, \textit{MAU assignment shift (\%)}, measures how widely the MAU formation assignments deviate from the baseline. It is calculated as the percentage of MAV-section pairs where the number of assigned MAUs differs from the baseline solution. The numerator is the total number of mismatched MAV-section assignments, and the denominator is the total number of all MAV-section assignments across the network. Finally, the last column reports the total number of used MAUs.

\begin{table}[h]
\centering
\caption{Post hoc robustness of fine-tuning strategies under demand surges and re-distributions.}
\label{tab:robust_metrics}
\begin{threeparttable}
\resizebox{\textwidth}{!}{%
\begin{tabular}{lrrrrr}
\toprule
Instance
& \begin{tabular}[c]{@{}c@{}}Objective \\ value (\$)\end{tabular}
& \begin{tabular}[c]{@{}c@{}} Expected number  of \\unserved  passengers\end{tabular}
& \begin{tabular}[c]{@{}c@{}}Top-5\% Avg. \\ load ratio (\%)\end{tabular}
& \begin{tabular}[c]{@{}c@{}} MAU assignment \\ shift (\%)\end{tabular}
& \begin{tabular}[c]{@{}c@{}}Total number \\ of used MAUs\end{tabular} \\
\midrule
D5\_U\_AVS          &  9532.88 & 0  & --   &  7.6 & 112 \\
D5\_U\_FTTAVS          &  9216.68 & 0  & 0.89 &  13.7 & 109 \\
D5\_HB\_AVS         &  9575.65 & 0  & 0.87  & 10.3 & 115 \\
D5\_HB\_FTTAVS         &  9291.48 & 0  & 0.89 &  16.2 & 115 \\
D5\_TB\_AVS         &  9511.42 & 0  & 0.87  &  7.4 & 117 \\
D5\_TB\_FTTAVS         &  9195.74 & 0  & 0.89 & 
14.3 & 111 \\
D10\_U\_AVS         &  9915.67 & 0  & 0.86  & 17.7 & 118 \\
D10\_U\_FTTAVS         &  9585.19 & 0  & 0.88 &  21.4 & 118 \\
D10\_HB\_AVS        &  9933.14 & 0  & 0.88 & 16.5 & 117 \\
D10\_HB\_FTTAVS        &  9642.35 & 0  & 0.89 &  25.3 & 116 \\
D10\_TB\_AVS        &  9864.62 & 0  & 0.88 & 14.2 & 115 \\
D10\_TB\_FTTAVS        &  9530.31 & 0  & 0.89 &  20.0 & 113 \\
D15\_U\_AVS        & 10239.50 & 0  & 0.88  & 21.3 & 120 \\
D15\_U\_FTTAVS         &  9929.97 & 0  & 0.88 &  26.7 & 124 \\
D15\_HB\_AVS        & 10291.87 & 0  & 0.90  & 22.7 & 127 \\
D15\_HB\_FTTAVS        & 10024.40 & 0  & 0.90 &  33.8 & 128 \\
D15\_TB\_AVS        & 10229.90 & 0  & 0.88 & 22.6 & 120 \\
D15\_TB\_FTTAVS        &  9905.54 & 0  & 0.88 &  25.8 & 119 \\
D20\_U\_AVS         & 10629.34 & 8  & 0.89  & 26.7 & 128 \\
D20\_U\_FTTAVS         & 10328.20 & 0  & 0.89 &  36.4 & 131 \\
D20\_HB\_AVS        & 10660.00 & 12 & 0.90  & 29.4 & 128 \\
D20\_HB\_FTTAVS        & 10356.40 & 0  & 0.90 &  34.8 & 133 \\
D20\_TB\_AVS        & 10614.61 & 0  & 0.89  & 29.3 & 126 \\
D20\_TB\_FTTAVS        & 10257.40 & 0  & 0.89 & 30.7 & 129 \\
\bottomrule
\end{tabular}%
}
\end{threeparttable}
\end{table}

From the results shown in Table~\ref{tab:robust_metrics}, we observe that across all instances, the FTTAVS strategy consistently yields lower objective values than the AVS strategy. Under severe surges, AVS fails to serve all passengers (e.g., 8 unserved in D20\_U\_AVS and 12 in D20\_HB\_AVS), whereas FTTAVS accommodates full demand. A second observation is that the Top-5\% average load ratio remains stable between 86\% and 90\% as demand increases from 5\% to 20\%, suggesting that both strategies effectively adapt vehicle schedules to the changing demand. Notably, FTTAVS achieves slightly higher Top-5\% load ratios than AVS, indicating a tighter capacity–demand alignment made possible by jointly fine-tuning the timetable and re-optimizing vehicle schedules. These gains are accompanied by broader reallocations: the MAU assignment shift is consistently larger under FTTAVS, reflecting more extensive re-optimization of vehicle schedules compared to AVS. Additionally, for demand surges of 5\% to 10\%, FTTAVS uses the same or fewer MAUs than AVS (e.g., D5\_U: 109 vs 112; D10\_HB: 116 vs 117; D10\_TB: 113 vs 115). However, when surges reach 15\% to 20\%, FTTAVS typically employs more MAUs (e.g., D20\_U: 131 vs 128; D20\_HB: 133 vs 128). This shift reflects the trade-off embedded in the objective function: under small surges, operational costs dominate, incentivizing fleet conservation; beyond a 15\% increase, passengers' costs become more prominent, and increasing fleet size becomes cost-effective to improve service quality. In conclusion, these findings suggest that both strategies exhibit a degree of robustness to demand surges. Jointly fine-tuning the timetable and re-optimizing the vehicle schedule consistently yields superior performance, particularly under more substantial surges.

\subsection{Value of the stochastic programming approach}\label{stochastic_experiment_subsection}
In this section, we vary the number of scenarios to explore the effects of the number of scenarios on the efficiency of the algorithm and the quality of solutions. We generated 30 passenger flow scenarios based on real-world operational data, and randomly selected 1, 2, 4, 6, 8, 12, 15 and 20 demand scenarios from them to form 8 instances, respectively. Next, each of these instances is solved using the proposed RH+IL algorithm. Then, the obtained timetables are fed into the 30 scenarios to evaluate the expected passengers' costs and operational costs. 

Table \ref{tab:StochasticNumber} shows the average performance of the above seven solutions over 30 scenarios. It can be seen that as the number of scenarios increases, the average performance of the resulting solution over the entire sample increases accordingly, but its solution time also tends to increase. Specifically, when the number of scenarios is very small (e.g., 1 and 2), although the solution time is quite short, too few scenarios do not reflect the general pattern of passenger flow, and therefore the average performance of the obtained timetable is poor. In constant, when the number of scenarios is overly abundant, the solution quality improves accordingly, but the solution time also becomes dramatically longer. Choosing an appropriate number of scenarios (i.e., 4, 6, or 8) improves the solution quality by at least 6.5\% relative to a deterministic optimization; meanwhile, it considerably improves the computational efficiency of solving the TT-VS-DAC problem relative to an excessive number of scenarios (i.e., 12, 15 or 20).
\begin{table}[h]
\centering
\caption{Results among the various numbers of stochastic scenarios.}
\label{tab:StochasticNumber}
\begin{tabular}{rrrrrr}
\toprule
\begin{tabular}[c]{@{}c@{}}Number of \\scenarios \end{tabular}& \begin{tabular}[c]{@{}c@{}}Total passengers'\\ costs (\$)\end{tabular} & \begin{tabular}[c]{@{}c@{}}Total operational \\costs (\$) \end{tabular}&\begin{tabular}[c]{@{}c@{}} Number of \\vehicles\end{tabular} & \begin{tabular}[c]{@{}c@{}}Objective function \\value (\$) \end{tabular}& \begin{tabular}[c]{@{}c@{}}Computational \\time (s) \end{tabular}\\ \midrule
1                  & 8511.2              & 2314.1                & 112                & 10825.3            & 106.2             \\
2                  & 8554.6             & 2334.3                & 113                & 10589.4            & 253.5             \\
4                  & 7791.3              & 2324.2                & 112                & 10115.5            & 524.4            \\
6                  & 7599.1              & 2321.6                & 115                & 9920.7            & 1047.7            \\
8                  & 7591.9              & 2316.4                & 112                & 9908.2            & 1896.3            \\
12                 & 7468.8              & 2313.9                & 113                & 9782.6            & 5160.3           \\
15                 & 7381.6              & 2316.0                & 110                & 9697.5            & 7202.1            \\ 
20                 & 7381.6              & 2314.1                & 112                & 9696.1            & 9711.2           \\ \bottomrule
\end{tabular}
\end{table}

\subsection{Performance of the learning-based real-time decision-making framework}
\label{sec:real_time}

To evaluate the effectiveness of the proposed learning-based real-time decision-making framework and the embedded fine-tuning strategy, we conduct two sets of experiments based on the Beijing bus subnetwork. First, we build 12 instances to compare the performance of our framework against benchmark methods. Second, we construct six instances to investigate the impact of four levels of re-optimization flexibility within the real-time decision-making process.

\subsubsection{Performance comparison of the learning-based scenario-retention method and benchmarks.}

To validate the effectiveness of our proposed learning-based scenario-retention method, we compare it against two benchmark approaches: the \textit{All scenario} scenario-retention method, which incorporates the entire stochastic demand scenario set at each decision-making stage; and the \textit{Greedy} scenario-retention method, which ranks all scenarios in descending order of the number of total passenger demand and selects the top-$k$ scenarios for optimization. For all methods, we impose a one-minute time limit per stage to ensure real-time applicability. If no feasible solution is found within this time limit, the process is terminated. This set of experiments is based on instance\_r\_6, with newly generated stochastic demand scenario sets. Specifically, we construct scenario sets containing 8, 15, 20, and 30 scenarios, respectively, and generate three distinct sets for each scenario size. These scenario sets are newly generated to reflect diverse demand realizations for evaluation purposes. In addition, the value of $k$ is set as three for both \textit{Greedy} and our learning-based scenario-retention methods. Each departure time is allowed to shift within a narrow adjustment range of ±2 minutes.

The results are summarized in Table~\ref{tab:performance_comparison}. In the \textit{Instance} column, the notation $|\mathcal{W}|$ refers to the number of stochastic demand scenarios included in the scenario set, and C1, C2, C3 represent three distinct sets generated independently for each scenario size. The $Dev$ column reports the relative difference in total passengers' and operational costs compared to the Greedy method. A negative $Dev$ value indicates that the corresponding scenario-retention
method yields lower total costs (i.e., superior solution quality), whereas a positive value suggests inferior performance. An asterisk \enquote{*} indicates that the method failed to produce a feasible solution within the time limit, while \enquote{\#} denotes model infeasibility.

\begin{table}[h]
\centering
\caption{Performance comparison between the different scenario-retention methods.}
\label{tab:performance_comparison}
\begin{threeparttable}
\resizebox{\textwidth}{!}{%
\begin{tabular}{ccrrrccrrr}
\toprule
Instance & \begin{tabular}[c]{@{}c@{}}Scenario-retention\\ method\end{tabular}  &{\begin{tabular}[c]{@{}c@{}}Passengers' \\  cost\end{tabular}} & {\begin{tabular}[c]{@{}c@{}}Operating \\  cost\end{tabular}} & {$Dev$ (\%)} &
Instance & \begin{tabular}[c]{@{}c@{}}Solution\\ method\end{tabular}  & {\begin{tabular}[c]{@{}c@{}}Passengers' \\  cost\end{tabular}} & {\begin{tabular}[c]{@{}c@{}}Operating \\  cost\end{tabular}} & {$Dev$ (\%)} \\
\midrule

\multirow{3}{*}{$|\mathcal{W}|$=8, C1}
  & Greedy        & 7104.89 & 2265.41 & \textendash     &
\multirow{3}{*}{$|\mathcal{W}|$=20, C1}
  & Greedy        & 6520.29 & 2258.80 & \textendash \\
  & All scenario  & 11490.49 & 2176.78 & 31.4 &
  & All scenario  & *      & *     & \textendash \\
  & Learning-based            & 6863.18 & 2260.95 & -2.6 &
  & Learning-based            & 6385.62 & 2230.76 & -1.9 \\

\addlinespace

\multirow{3}{*}{$|\mathcal{W}|$=8, C2}
  & Greedy        & 7411.18 & 2287.62 & \textendash     &
\multirow{3}{*}{$|\mathcal{W}|$=20, C2}
  & Greedy        & 5930.95 & 2215.54 & \textendash \\
  & All scenario  & 12547.57 & 2201.11 & 52.1 &
  & All scenario  & *      & *     & \textendash \\
  & Learning-based            & 6975.31 & 2292.77 & -4.4 &
  & Learning-based            & 5963.28 & 2216.69 & 0.4 \\

\addlinespace

\multirow{3}{*}{$|\mathcal{W}|$=8, C3}
  & Greedy        & 7119.63 & 2247.40 & \textendash    &
\multirow{3}{*}{$|\mathcal{W}|$=20, C3}
  & Greedy        & 6140.54 & 2227.32 & \textendash \\
  & All scenario  & 11591.99 & 2176.96 & 47.0 &
  & All scenario  & *      & *     & \textendash \\
  & Learning-based            & 6698.05 & 2248.48 & -4.5 &
  & Learning-based            & 5914.03 & 2217.00 & -2.8 \\

\addlinespace
\midrule

\multirow{3}{*}{$|\mathcal{W}|$=15, C1}
  & Greedy        & 6229.74 & 2225.63 & \textendash     &
\multirow{3}{*}{$|\mathcal{W}|$=30, C1}
  & Greedy        & 6960.45 & 2234.45 & \textendash \\
  & All scenario  & *       & *     & \textendash     &
  & All scenario  & *       & *      & \textendash \\
  & Learning-based            & 5566.56 & 2217.73 & -7.9 &
  & Learning-based            & 6528.14 & 2236.00 & -4.7 \\

\addlinespace

\multirow{3}{*}{$|\mathcal{W}|$=15, C2}
  & Greedy        & 6206.62 & 2220.60 & \textendash     &
\multirow{3}{*}{$|\mathcal{W}|$=30, C2}
  & Greedy        & \#       & \#       & \textendash  \\
  & All scenario  & *       & *      & \textendash     &
  & All scenario  & *       & *      & \textendash \\
  & Learning-based            & 5931.93 & 2224.79 & -3.2 &
  & Learning-based            & 6888.64 & 2267.09 & \textendash \\

\addlinespace

\multirow{3}{*}{$|\mathcal{W}|$=15, C3}
  & Greedy        & 6398.02 & 2262.52 & \textendash     &
\multirow{3}{*}{$|\mathcal{W}|$=30, C3}
  & Greedy        & \#       & \#      & \textendash \\
  & All scenario  & *       & *      & \textendash     &
  & All scenario  & *       & *      & \textendash \\
  & Learning-based            & 6375.61 & 2259.27 & -0.3 &
  & Learning-based            & 6965.53 & 2259.40 & \textendash \\

\bottomrule
\end{tabular}
}
\begin{tablenotes}
\footnotesize
\item [1] $|\mathcal{W}|$ refers to the number of stochastic demand scenarios included in the scenario set, and C1, C2, C3 represent three distinct \\ sets generated independently for each scenario size.
\item [2] * indicates no feasible solution was found within the allowed computation time; \# indicates the model is infeasible.
\item [3] $Dev (\%)$ represents the relative difference in total passengers' and operational costs compared to the Greedy method.
\end{tablenotes}
\end{threeparttable}
\end{table}

From the results shown in Table~\ref{tab:performance_comparison}, we can observe that the \textit{All scenario} scenario-retention method fails to find feasible solutions within the one-minute computational time limit in 9 out of the 12 instances. This outcome is mainly due to the excessive computational burden associated with incorporating the entire scenario set into the optimization within the one-minute solution time limitation. In contrast, our proposed learning-based scenario-retention method generates feasible solutions across all instances, demonstrating its suitability for real-time optimization. Moreover, it consistently outperforms the \textit{Greedy} benchmark in terms of total passengers' and operational costs. Across the 12 tested instances, our method yields lower total costs in 11 cases, with relative cost reductions ranging from 0.3\% (for $W=15$, C3) to 7.9\% (for $W=15$, C1). In several large-scale instances (e.g., $W=30$, C2 and C3), the \textit{Greedy} method fails to return feasible solutions. This is largely due to its strategy of selecting scenarios solely based on the number of passenger demand, which overlooks scenarios where local capacity constraints are more binding. These findings highlight the superior robustness and efficiency of our learning-based scenario-retention method in delivering high-quality solutions under real-time computational requirements.

\subsubsection{Impact of re-optimization flexibility in the learning-based real-time decision-making framework.}

Furthermore, to assess the impact of re-optimization flexibility within the learning-based real-time decision-making framework, we compare four strategies: (i) \textit{TT+VS}: both the timetable and vehicle schedules can be dynamically re-optimized; (ii) \textit{TT+VS (partial)}: same as (i), but boarding commitments from the previous stage are preserved, ensuring passengers already on board are not rerouted; (iii) \textit{VS}: the timetable is fixed, and only vehicle schedules are dynamically re-optimized; (iv) \textit{Fixed}: no re-optimization is allowed. In this set of experiments, we focus on instance\_r\_6. The timetable and vehicle schedule obtained from our integrated optimization approach are used as the baseline solution. Each departure time is allowed to shift within a narrow window of ±2 minutes. To evaluate the learning-based real-time framework under the four re-optimization strategies, we construct six new demand scenario sets, each containing 8 scenarios that follow the characteristics of instance\_r\_6.

Table~\ref{tab:rh_strategy_results} summarizes the results, reporting the relative deviations (\%) in total passengers' and operational costs compared to the \textit{Fixed} strategy, the presence of any unserved passengers, and the number of used MAUs. The negative value of $Dev.$ represents the reduction of the total costs. We can observe that fixing the timetable and vehicle schedule with no re-optimization leads to unserved passengers in four of six instances (i.e., C2, C3, C5, C6). The reason is that the space-time distribution of passenger demand changes in the new scenario sets while the operational plans are frozen and overloading is disallowed, so some passengers can no longer board. 

A second observation is that allowing only vehicle rescheduling (VS) resolves all infeasibilities and lowers total costs relative to the \textit{Fixed} strategy in every case from 0.81\% to 3.47\%. It also reallocates capacity more efficiently, typically using fewer MAUs than the \textit{Fixed} strategy (e.g., C2: 89 vs. 91; C3: 88 vs. 93; C5: 101 vs. 106; C6: 105 vs. 113). In addition, TT+VS (partial) yields additional gains over VS in most instances. To honor already-boarded passengers, this policy occasionally needs a slightly larger fleet than VS (e.g., C3 and C5), but it maintains passenger continuity and remains feasible throughout. Granting full flexibility to fine-tune the timetable and re-optimize vehicle schedules (i.e., TT+VS) delivers the best performance overall, with the largest cost reductions in every instance from 3.15\% to 8.14\%, and with the fleet size that is comparable to the other flexible strategies. In conclusion, relying on vehicle schedules and timetables optimized based on historical demand scenarios poses operational risks when facing real-time information in practical operations. Re-optimizing only the vehicle schedule can improve service quality by ensuring that all passengers are served. Preserving onboard commitments enhances the passenger experience, although it typically results in slightly higher total costs compared to the strategy that does not account for such commitments. Among these four strategies, the approach that dynamically re-optimizes both the timetable and vehicle schedules (i.e., TT+VS) consistently achieves the lowest total cost.

\begin{table}[h]
\centering
\caption{Performance comparison among various re-optimization strategies.}
\label{tab:rh_strategy_results}
\resizebox{\textwidth}{!}{%
\begin{tabular}{ccrcrccrcr}
\toprule
\cmidrule(lr){1-5}\cmidrule(lr){6-10}
Instance & Strategy &\begin{tabular}[c]{@{}c@{}}  $Dev. (\%)$ \\ (Total costs)\end{tabular} & \begin{tabular}[c]{@{}c@{}} Unserved \\ passengers\end{tabular} & \begin{tabular}[c]{@{}c@{}}Total number \\ of used MAUs\end{tabular} &
Instance & Strategy  &\begin{tabular}[c]{@{}c@{}}  $Dev. (\%)$ \\ (Total costs)\end{tabular} & \begin{tabular}[c]{@{}c@{}} Unserved \\ passengers\end{tabular} & \begin{tabular}[c]{@{}c@{}}Total number \\ of used MAUs\end{tabular} \\
\midrule
\multirow{4}{*}{C1}  & Fixed        & & no  & 85  &  
\multirow{4}{*}{C4}   & Fixed       &  & no  & 101 \\
      & VS          & -1.44 & no  & 85  &       & VS         & -0.81 & no  & 99  \\
    & TT+VS (partial)  & -3.40 & no  & 87  &       & TT+VS (partial)   &-2.19 & no  & 96  \\
 & TT+VS          & -3.92 & no  & 86  & & TT+VS          & -3.49 & no  & 97  \\
\addlinespace
\multirow{4}{*}{C2} & Fixed           & & yes & 91  &
\multirow{4}{*}{C5}& Fixed          & & yes & 106 \\
      & VS         & -1.68 & no  & 89  &       & VS        & -1.12 & no  & 101 \\
      & TT+VS (partial)  &-4.43 & no  & 89  &       & TT+VS (partial)   &-1.48 & no  & 104 \\
 & TT+VS          & -8.14& no  & 86  & & TT+VS           &-3.15 & no  & 101 \\
\addlinespace
\multirow{4}{*}{C3} & Fixed          & & yes & 93  & 
\multirow{4}{*}{C6}& Fixed           & & yes & 113 \\
      & VS        & -1.87 & no  & 88  &       & VS         & -3.47 & no  & 105 \\
          & TT+VS (partial) & -2.79 & no  & 91  &       
          & TT+VS (partial)  & -3.60 & no  & 99  \\
 & TT+VS           & -5.38 & no  & 90  &  & TT+VS         & -5.48 & no  & 102 \\
\bottomrule
\end{tabular}%
}
\end{table}

\subsection{Analysis of the benefits in the practical instance and the trade-offs}
\label{sec:tradeoff_analysis}
In this section, we explore the benefits that the proposed approach can provide to real-world operations by analysing the results of $instance\_r\_6$. We first explore the benefit of adopting our proposed flexible-capacity operational strategy. The results are presented in Table \ref{tab:DifferentStrategy}. It can be observed that adopting the flexible-capacity operational strategy can significantly reduce operational costs compared to the fixed-capacity operational strategy. Specifically, the implementation of the proposed flexible-capacity strategy approximately results in a 48\% reduction in operational costs and a 57\% decrease in the number of utilized MAUs. At the same time, passengers' costs also decrease from 6861.3 to 6757.5. The main reason is that our proposed method provides passengers with the convenience to make in-vehicle transfers, which reduces the waiting time for transfers in transfer corridors.

\begin{table}[h]
\caption{The total passengers’ costs, operational costs, and the number of used MAUs obtained by employing the fixed-capacity and flexible-capacity strategies in the real-world instance.}
\label{tab:DifferentStrategy}
\resizebox{\textwidth}{!}{%
\begin{tabular}{crrrr}
\toprule
                            & Operational strategy                   & Total passengers' costs (\$) & Total operational costs (\$) & The number of used MAUs \\
\midrule

& Fixed-capacity     & 6861.3      & 4381.2   & 264   \\
& Flexible-capacity  & 6757.5      & 2275.8   & 112                \\ \bottomrule
\end{tabular}%
}
\end{table}

Next, we use the weighted sum approach to analyze the tradeoff between operational and passengers' costs. The objective function (denote as $\Psi$ here) is reformulated as follows:
\begin{align}
    \Psi =\omega_1\frac{\Psi^{pass}-\underline{\Psi}^{pass}}{\overline{\Psi}^{pass}-\underline{\Psi}^{pass}} + \omega_2 \frac{\Psi^{oper}-\underline{\Psi}^{oper}}{\overline{\Psi}^{oper}-\underline{\Psi}^{oper}}.
\end{align}
To be specific, $\omega_1$ and $\omega_2$ are weighting coefficients designed to control the magnitude of the impact that normalised passenger travel costs and normalised bus operational costs have on the change in objectives, respectively. $\underline{\Psi}^{pass}$ is the lower bound on the total passenger travel cost, and $\overline{\Psi}^{pass}$ is the upper bound on the total passenger travel cost. $\underline{\Psi}^{pass}$ is theoretically minimised to 0, and $\overline{\Psi}^{pass}$ can be calculated by taking the maximum headway for both waiting time per passenger and transfer waiting time. Similarly, $\underline{\Psi}^{oper}$ is the lower bound on the total operating travel cost and $\overline{\Psi}^{oper}$ is the upper bound on the total operating travel cost. $\underline{\Psi}^{oper}$ can be calculated according to Lemma 1, and $\overline{\Psi}^{oper}$ can be computed by assuming that the maximum fleet is dispatched for each trip and does not allow for decoupling and coupling en-route.

\begin{figure}[]
\centering
\includegraphics[height=7cm]{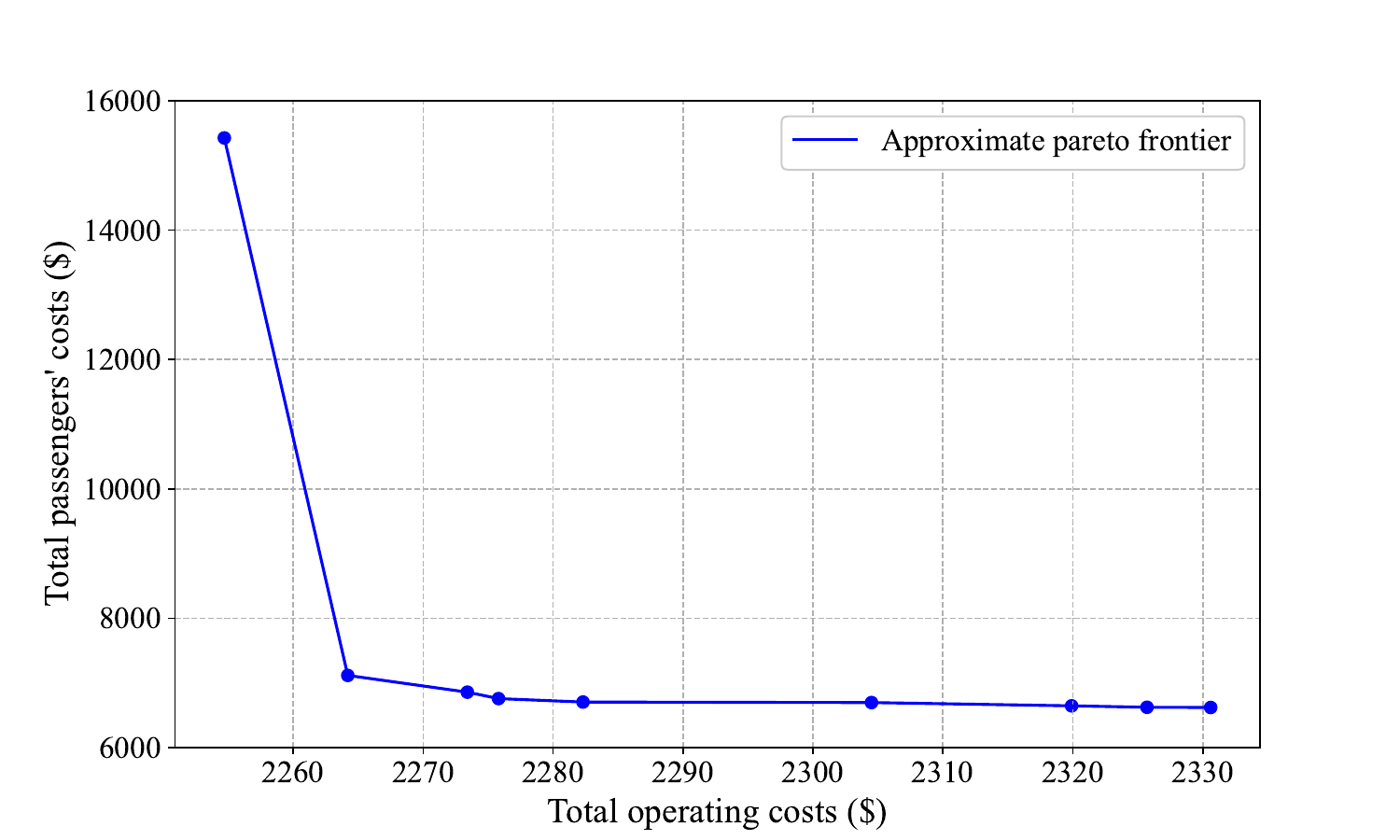}
\caption{Approximate Pareto frontier of the TT-VS-DCA model.}\label{fig:Approximate_pareto_frontier}
\end{figure}

To approximate the Pareto optimal point, we incrementally increase $\omega_1$ from 0 to 1 while simultaneously decreasing $\omega_2$ from 1 to 0 in increments of 0.1. The results obtained through the weighted sum method are depicted in Figure \ref{fig:Approximate_pareto_frontier}. Upon analyzing the results, two key observations emerge: (1) Excessive emphasis on reducing operational costs can precipitate a significant decline in service quality, thereby leading to a substantial rise in passengers' costs. (2) Reducing operational costs with little to no negative impact on service quality can be achieved with a slight increase in the weight assigned to operational costs.

\section{Conclusions}
\label{sec:conclusion}

We present a novel modeling approach for the TT-VS-DCA problem with time-varying and uncertain passenger demand. Our approach explicitly considers the trade-off between passengers’ costs and operational costs. MAUs can be flexibly (de)coupled at depots and transfer stops, enabling efficient circulations and rerouting across lines. To solve this integrated optimization problem, we propose an exact solution method with tailored valid inequalities and embed it in a rolling horizon framework to enhance computational efficiency. To facilitate real-time adjustments under evolving operational conditions, we develop a real-time decision-making framework to dynamically fine-tune timetables which are obtained based on historical passenger flow scenarios and optimize vehicle scheduling.

Computational experiments on both a virtual network and the Beijing bus subnetwork demonstrate the effectiveness of our approach. Several key takeaways are summarized as follows: (i) Our exact method outperforms GUROBI for medium-sized instances, while the rolling-horizon optimization framework combined with the integer L-shaped approach yields superior results for large-scale instances. (ii) Compared to the sequential optimization method, the proposed integrated approach consistently performs better across instances of varying scales. For smaller instances, it results in fewer required MAUs. As the network size increases, it remains capable of generating feasible timetables and vehicle schedules, whereas the sequential method often fails to find feasible solutions. (iii) Allowing MAVs to (de)couple and be rerouted at transfer stops, in addition to depots, further reduces both fleet size and operational costs. (iv) The optimized timetable and vehicle schedule, generated from historical demand scenarios using our integrated optimization approach, are robust under demand decreases in out-of-sample scenarios. Moreover, allowing overloading can handle a 5\% demand surge without requiring additional capacity resources. (v) Both the strategy that fixes the timetable while re-optimizing the vehicle schedule and the strategy that jointly fine-tunes these two schedules can handle the demand surges and redistributions and exhibit robustness. (vi) The proposed learning-based real-time decision-making framework outperforms benchmark methods, thanks to the novel scenario-retention approach using the machine learning model.

\bibliographystyle{informs2014trsc} 
\bibliography{cas-refs} 

\newpage

\begin{appendices}
	
\section{MAVs developed by Next Future Transportation Inc. }
\label{sec:exampleMAV}

\begin{figure}[h]
    \centering
    \includegraphics[height=5cm]{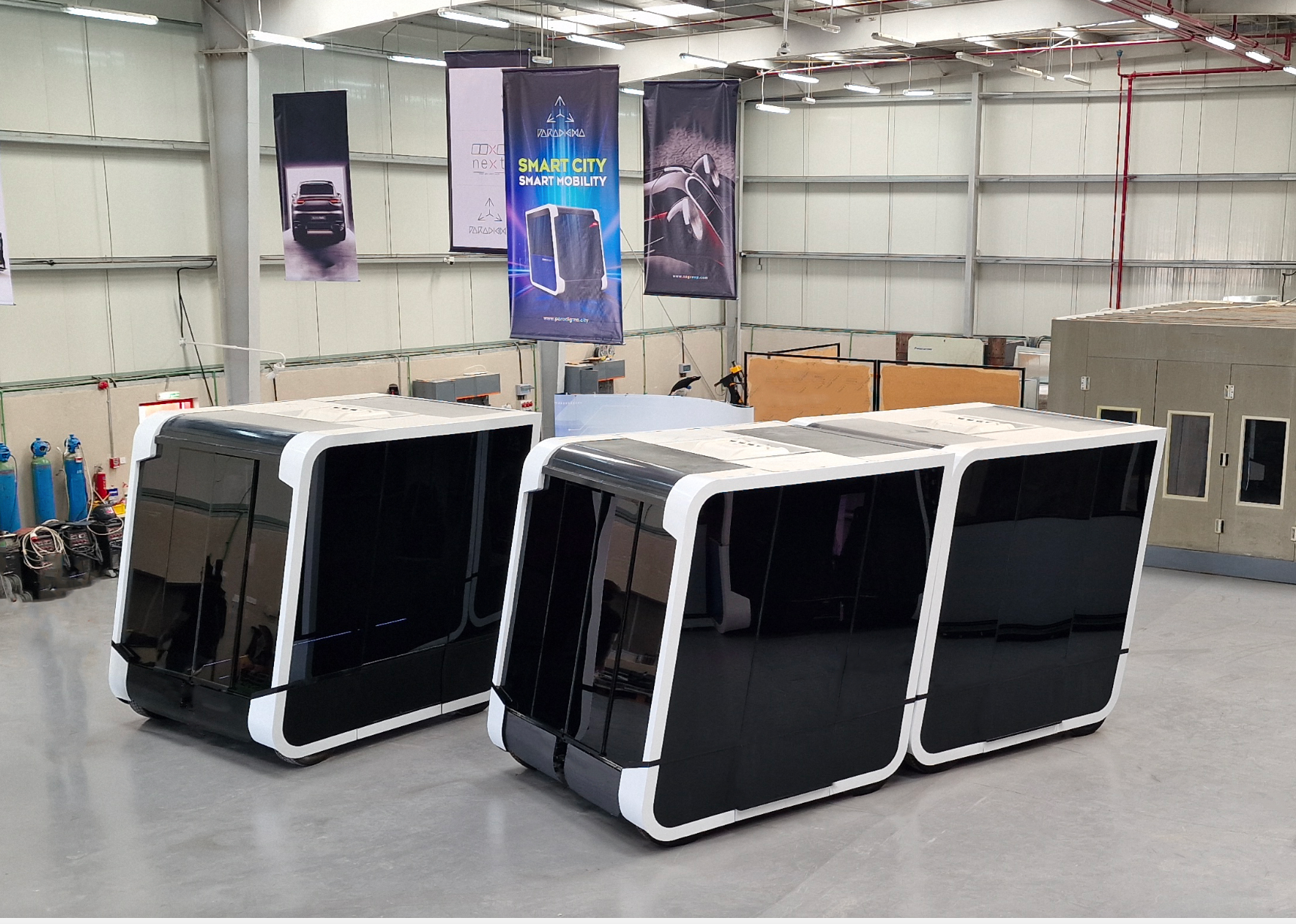}
    \caption{Modular autonomous vehicles developed by the Next Future Transportation Inc. \citep{nextMAV}.}
    \label{fig:MAV}
\end{figure}

\section{Overview of the closely related research}
\label{sec:research}
Table \ref{tab:literature} summarizes the research closely related to this study. 

\begin{table}[h]
\centering
\caption{Overview of the research closely related to this study.}
\label{tab:literature}
\begin{threeparttable}
\resizebox{\linewidth}{!}{%
\begin{tabular}{r|ccccccccccccc}
\hline
\multirow{2}{*}{Publication} & \multirow{2}{*}{\begin{tabular}[c]{@{}c@{}}Network\\ level\end{tabular}} & \multicolumn{2}{c}{Passenger demeand} &  & \multicolumn{2}{c}{Dynamic formation} &  & \multirow{2}{*}{\begin{tabular}[c]{@{}c@{}}One-line vehicle\\ circulation \\ scheduling\end{tabular}} &  & \multicolumn{2}{c}{\begin{tabular}[c]{@{}c@{}}Cross-line vehicle \\ circulation scheduling\end{tabular}}                                         & \multirow{2}{*}{\begin{tabular}[c]{@{}c@{}}Solution\\ method\end{tabular}} \\ \cline{3-4} \cline{6-7} \cline{11-12}
                             &            & Dynamic          & Uncertain          &  & At depots          & At stops         &  &                         &  & \begin{tabular}[c]{@{}c@{}}At depots\\ before/after trips\end{tabular} & \begin{tabular}[c]{@{}c@{}}At transfer stops\\ during trips\end{tabular} &                  \\ \hline
 \cite{Alfieri2006} &   \ding{55}& \ding{55}  &   \ding{55}      &  &  \ding{55}    &  \ding{55}  &  &  \ding{51}  &  &  \ding{55}  &  \ding{55}  & CS   \\
\cite{Schobel2009} &   \ding{51}& \ding{55}  &   \ding{55}      &  &  \ding{55}    &  \ding{55}  &  &  \ding{51}  &  &  \ding{55}  &  \ding{55}  & HE   \\
\cite{Ibarra-Rojas2014}&   \ding{51}& \ding{55}  &   \ding{55}      &  &  \ding{55}    &  \ding{55}  &  &  \ding{51}  &  &  \ding{55}  &  \ding{55}  & CS   \\
\cite{laporte2017multi}&   \ding{51}& \ding{51}  &   \ding{55}      &  &  \ding{55}    &  \ding{55}  &  &  \ding{51}  &  &  \ding{55}  &  \ding{55}  & CS   \\
\cite{Rolf2021} &   \ding{51}& \ding{55}  &   \ding{55}      &  &  \ding{55}    &  \ding{55}  &  &  \ding{51}  &  &  \ding{51}  &  \ding{55}  & CS   \\
\cite{RolfPartB} &   \ding{51}& \ding{55}  &   \ding{55}      &  &  \ding{55}    &  \ding{55}  &  &  \ding{51}  &  &  \ding{51}  &  \ding{55}  & CS   \\
\cite{Amberg2023}&   \ding{51}& \ding{55}  &   \ding{51}      &  &  \ding{55}    &  \ding{55}  &  &  \ding{51}  &  &  \ding{55}  &  \ding{55}  & CG   \\
\cite{CHEN20191}&   \ding{55}& \ding{51}  &   \ding{55}      &  &  \ding{51}    &  \ding{55}  &  &  \ding{55}  &  &  \ding{55}  &  \ding{55}  & DP  \\
\cite{Xiaowei2021}&   \ding{55}& \ding{51}  &   \ding{55}      &  &  \ding{51}    &  \ding{55}  &  &  \ding{55}  &  &  \ding{55}  &  \ding{55}  & DP  \\
\cite{Chen2022}&   \ding{55}& \ding{51}  &   \ding{55}      &  &  \ding{51}    &  \ding{51}  &  &  \ding{55}  &  &  \ding{55}  &  \ding{55}  & HE  \\
\cite{TIAN2023103986}&   \ding{55}& \ding{51}  &   \ding{55}      &  &  \ding{51}    &  \ding{51}  &  &  \ding{55}  &  &  \ding{55}  &  \ding{55}  & HE  \\
\cite{XIA2023}&   \ding{55}& \ding{51}  &   \ding{51}      &  &  \ding{51}    &  \ding{51}  &  &  \ding{55}  &  &  \ding{55}  &  \ding{55}  & IL  \\
\cite{Xia2024}&   \ding{55}& \ding{51}  &   \ding{51}      &  &  \ding{51}    &  \ding{51}  &  &  \ding{55}  &  &  \ding{55}  &  \ding{55}  & HY  \\
\textbf{This paper}&   \ding{51}& \ding{51}  &   \ding{51}      &  &  \ding{51}    &  \ding{51}  &  &  \ding{51}  &  &  \ding{51}  &  \ding{51}  & \begin{tabular}[c]{@{}c@{}}\textbf{IL/}\\ \textbf{IL + RH/} \\ \textbf{IL + RH + ML} \end{tabular}  \\
\hline
\end{tabular}}%
\begin{tablenotes}
    \footnotesize
    \item[a] ``Network level" denotes that the study focuses on a network consisting of at least two public transportation lines.
    \item[b] Abbrevations: CS = Commercial solver, HE = Heuristic, CG = Column generation, DP = Dynamic programming, \\
    IL = Integer L-shaped, HY = Hybird algorithm combining heuristics and a commercial solver,  RH = Rolling horizon framework, \\ ML = Machine Learning. 
\end{tablenotes}
\end{threeparttable}
\end{table}

\newpage

\section{An example of the vehicle scheduling and dynamic capacity allocation with cross-line circulations}
\label{sec:examples}
\begin{example}\label{ex:2}
Figure~\ref{ex:movement of mvs} illustrates vehicle scheduling and dynamic capacity allocation with cross-line circulations. An MAV consisting of two MAUs departs from stop 1-Line $a’$. Because many passengers board at stop 2-Line $a’$ with the destination of stop 3-Line $b’$, one MAU from the MAV operating on Line $a’$ is decoupled at transfer stop 2-Line $a’$. This MAU immediately couples with the MAV on Line $b’$, enabling onboard passengers to complete in-vehicle transfers. Lines $a’$ and $b’$ do not share a depot, showing that rerouting at a transfer stop is not restricted to lines served by the same depot.

When the three-unit MAV on Line $b’$ reaches the depot near stop 3-Line $b’$, the trip ends, and the next assignment of trips for the three MAUs are decided independently. This depot serves Lines $b$, $b’$, $c$, and $c’$. Therefore, each MAU may either remain at the depot, be scheduled for a subsequent trip on Line $b$ or $b’$, or be assigned to operate on Line $c$ or $c’$. In the latter case, the assignment occurs at the depot without any deadheading. To summarize, all cross-line circulations in this example occur at either transfer stops or depots, and MAUs never move empty between depots or from a depot to a transfer stop.
\begin{figure}[h]
\centering
\includegraphics[height=12cm]{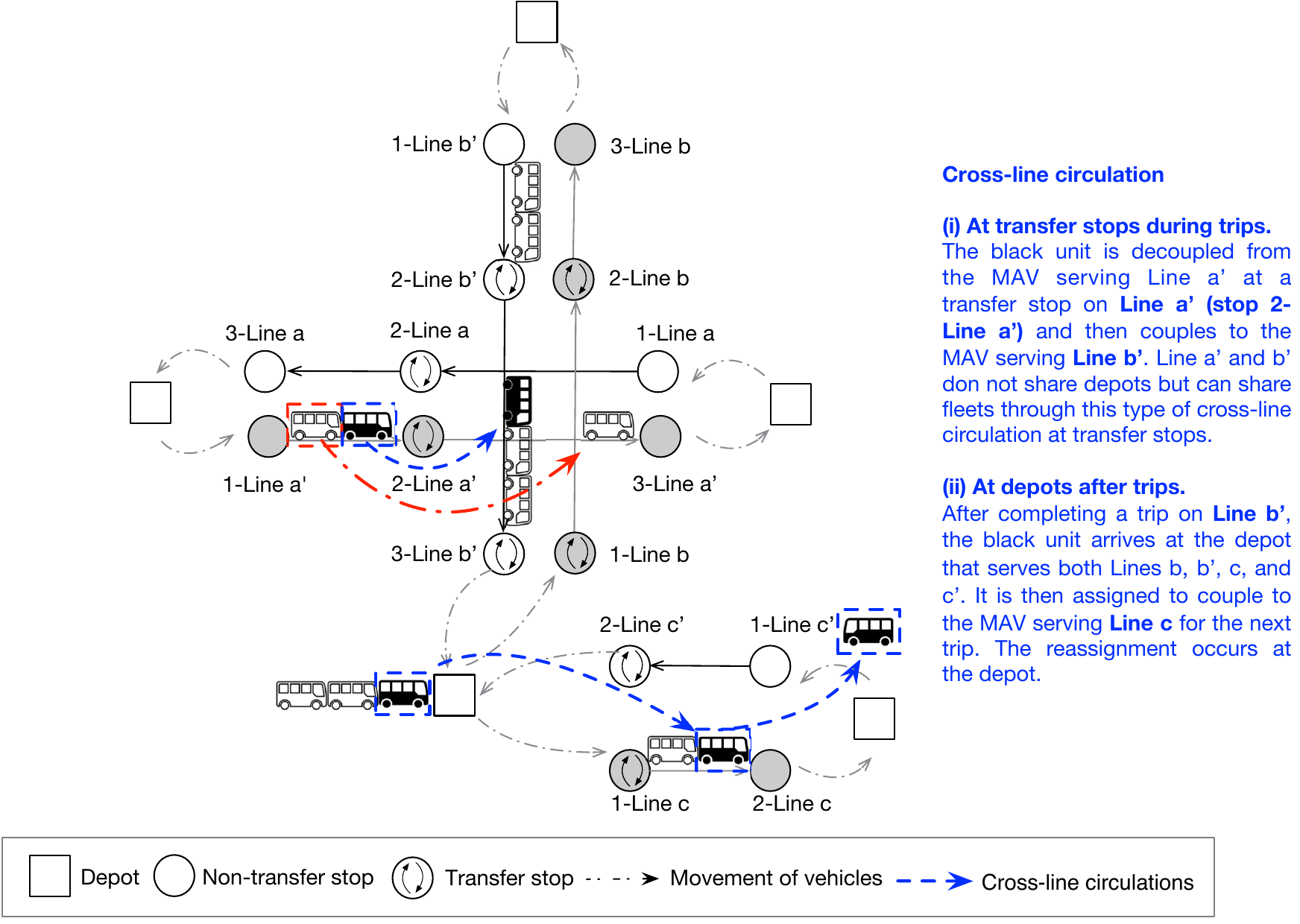}
\caption{An illustration of the vehicle scheduling and dynamic capacity allocation with cross-line circulations.}\label{ex:movement of mvs}
\end{figure}

\end{example}

\newpage

\section{Notations}
\label{sec:Notations}

Tables \ref{tab:notations} and \ref{tab:dependentVariables} introduce the notations of sets, parameters and variables used in our proposed mathematical formulation.
\begin{table}[h]
\centering
\caption{Notations of sets and parameters.}
\label{tab:notations}
\resizebox{\textwidth}{!}{%
\begin{tabular}{ll}
  \hline
  Notation & Description \\
\hline\noalign{\smallskip}
\multicolumn{2}{l}{\bf{Sets}}\\
$\mathcal{L}$& Set of lines, $\mathcal{L}=\{1,2,\ldots, \left|\mathcal{L}\right|\},$ indexed by $l$\\
$\mathcal{S}_l$& Set of stops on line $l$, $\mathcal{S}_l=\{1, 2, \ldots,  \left|\mathcal{S}_l\right| \},$ indexed by $i$ and $j$\\
$\mathcal{R}_l / \hat{\mathcal{R}}_{l}$& Set of transfer / non-transfer stops on line $l$, indexed by $i$ and $j$\\
$\dot{\mathcal{A}}$& Set of transfer corridors, $\dot{\mathcal{A}}=\{g: g=(l, l^{'}, i, i^{'}), i \in \mathcal{S}_{l}, i^{'} \in \mathcal{S}_{l'}, l \neq l^{'}, l, l^{'} \in \mathcal{L} \},$ indexed by $g$\\
$\mathcal{K}_l$& Set of trips on line $l$, $\mathcal{K}_l=\{1, 2, \ldots,  \left|\mathcal{K}_l\right| \},$ indexed by $k$\\
$\mathcal{T}$ & Set of discretized time intervals, $\mathcal{T}=\{1,2,\dots, \left|\mathcal{T}\right|\},$ indexed by $t$\\
$\mathcal{W}$ & Set of scenarios, $\mathcal{W}=\{1,2,\dots, \left|\mathcal{W}\right|\},$ indexed by $w$\\
$\mathcal{P}_w$ & Set that containing all passenger groups in scenario $w$, $\mathcal{P}=\{1,2,\dots, \left|\mathcal{P}_w\right|\},$ indexed by $p$\\
$R_p$& Set of transfer corridors of passenger group $p$\\
$\mathcal{L}_p$& Set of lines passed by passenger group $p$, $\mathcal{L}_p=\{l_0,l_1,\dots, \left|L_p\right|\},$\\
$\mathcal{M}$& Set of depots, $\mathcal{M}=\{1,2,\ldots, \left|\mathcal{M}\right|\},$ indexed by $m$\\
$FL_m$& Sets of lines whose starting stop is depot $m$\\
$LL_m$& Sets of lines whose ending stop is depot $m$\\
$\mathcal{Q}$ & Set of the number of MAUs that can be contained in an MAV, $\mathcal{Q}=\{1,2,\dots, \left|\mathcal{Q}\right|\},$ indexed by $q$\\
\multicolumn{2}{l}{\bf{Parameters}}\\
$\Delta$& Duration of each time interval\\
$o^l_{k, 1}$ & Originally scheduled time of trip $k$ on line $l$ arriving at first stop\\
$r^l_{k, i}$ & Running time from stops $i$ to $i+1$ of an MAV assigned to trip $k$ on line $l$\\
$\varepsilon_1$/$\varepsilon_2$ &Weighting coefficient related to the waiting times at origins / transfer corridors\\
$\varphi_1$/$\varphi_2$ &Weighting coefficient related to passengers' / operational costs\\
$\vartheta^T$ & The equivalent monetary value per unit of passengers' waiting costs (unit: \$)\\
$\vartheta_{q,i}$ & Cost-related parameter of using $q$ MAUs when MAV leave stop $i$ (unit: \$)\\
$\alpha_{l,k}^{min}/\alpha_{l,k}^{max}$ & Minimum / Maximum shifting time of trip $k$ on line $l$.\\
$\beta_{l,k,i}^{min}/\beta_{l,k,i}^{max}$ & Minimum / Maximum dwell time of trip $k$ at stop $i$ on line $l$\\
$h_{min}^l/h_{max}^l$ & Minimum / Maximum headway of line $l$\\
$\epsilon_1,\epsilon_2$ & Tiny positive constants\\
$M_1,M_2,M3$ & Big positive constants\\
$i^o_p/i^d_p$& Origin / destination stop of passenger group $p$\\
$u_p$& Arrival time of passenger group $p$ at their origin stop\\
$n_p$& Number of passengers in passenger group $p$\\
$\rho_w$ &The occurring probability of scenario $w$\\
$t^{pre}_m$ &The preparation time at depots of each MAU to complete a trip before it is allocated to the next trip\\
$\underline{\theta}_g/\overline{\theta}_g$ &Minimum / maximum transfer time at transfer corridor $g$\\
$\text{CAP}$ & Capacity of one MAU\\
$\text{B}$ & Maximum number of available MAUs \\
\hline
\end{tabular}%
}
\end{table}

\begin{table}[h]
\centering
\caption{Notations of variables.}
\label{tab:dependentVariables}
\resizebox{\textwidth}{!}{%
\begin{tabular}{ll}
  \hline
  Notation & Description \\
\hline\noalign{\smallskip}
\multicolumn{2}{l}{\bf{Decision variables}}\\
$\alpha_{k,1}^l$ & Shifting time of trip $k$ at the starting
stop on line $l$ \\
$\beta_{k,i}^l$ & Dwell time of trip $k$ at stop $i$ on line $l$ \\
$z_{k,i,t}^l$& Binary variable. $z_{k,i,t}^l=1$ indicates trip $k$ on line $l$ has not departed from stop $i$; otherwise, $z_{k,i,t}^l=0$ \\
$x^l_{k,i}(w)$ & Real-valued number of MAUs composed in the MAV allocated to trip $k$ to serve the section between stops \\& $i$ and $i+1$ on line $l$ \\
$y^l_{k,i,q}(w)$& Binary variable. $y^l_{k,i,q}(w) = 1 $ represents the number of MAUs composed in the MAV assigned to trip $k$ to \\& serve the segment between stops $i$ and $i+1$ on line $l$ in scenario $w$ is $q$; otherwise, $y^l_{k,i,q}(w) = 0$.\\
$h^{g}_{k,k'}(w)$&  Binary variable. $h^{g}_{k,k'}(w) = 1$ indicates MAUs on the MAV assigned to trip $k$ on line $l$ will be decoupled at \\& transfer corridor $g$ and rerouted to serve the $k'$th trip on line $l'$; Otherwise, $h^{g}_{k,k'}(w) = 0$.\\
$\kappa_m$ & Number of MAUs stored at each depot $m \in \mathcal{M}$ at the beginning of the study time horizon\\
\multicolumn{2}{l}{\bf{Intermediate variables}}\\
$\Psi^{pass}$ & Total passengers' costs\\
$\Psi^{oper}$ & Total operational costs\\
$\Psi_p^1(w)$ & Waiting time at the origin stop of passenger group $p$ in scenario $w$\\
$\Psi_p^2(w)$ & Waiting time at transfer corridors of passenger group $p$ in scenario $w$\\
$a_{k,i}^l$ & Arrival time of trip $i$ at stop $u$ on line $l$\\
$d_{k,i}^l$&  Departure time of trip $i$ at stop $u$ on line $l$\\
$\xi_{p,k}^l(w)$ & Binary variable.  $\xi_{p,k}^l(w) = 1$ represents passenger group $p$ is able to board the MAV assigned to trip $k$ on \\& line $l$ at their origin stop in scenario $w$; otherwise, $\xi_{p,k}^l(w)=0$\\
$\chi^l_{p,k}(w)$ & Binary variable. $\chi^l_{p,k}(w) = 1$ indicates passenger group $p$ boards the MAV assigned to trip $k$ on line $l$ \\& in scenario $w$; otherwise, $\chi^l_{p,k}(w)=0$\\
$\pi^{g}_{k,k'}$ & Binary variable. $\pi^{g}_{k,k'} = 1$ indicates passengers can transfer from the MAV assigned to trip $k$ to the MAV  \\& assigned to trip $k'$ at transfer corridor $g$; otherwise, $\pi^{g}_{k,k'}=0$ \\
$\zeta_{p,k,k'}^{g}(w)$ & Binary variable. $\zeta_{p,k,k'}^{g}(w) = 1$ represents passenger group $p$ transfers from the MAV assigned to trip $k$ to the \\& MAV assigned to trip $k'$ at transfer corridor $g$ in scenario $w$; otherwise, $\zeta_{p,k,k'}^{g}(w)=0$ \\
$b_{k,i}^l(w)$ & Number of passengers who board the MAV assigned to trip $k$ at stop $i$ on line $l$ in scenario $w$ \\
$c_{k,i}^l(w)$ & Number of passengers who alight from the MAV assigned to trip $k$ at stop $i$ on line $l$ in scenario $w$ \\
$v_{k,i}^l(w)$ & Number of in-vehicle passengers on the MAV assigned to trip $k$ at stop $i$ on line $l$ in scenario $w$ \\
$e_{k, k'}^{g}$& Binary variable. $e_{k, k'}^{g} = 1$ indicates rerouting operations at transfer corridor $g$ is available; otherwise, $e_{k, k'}^{g}=0$\\
$AV_{m,t}(w)$ & Cumulative inflow of MAUs at depot $m \in \mathcal{M}$ at time $t$ in scenario $w$\\
$DV_{m,t}(w)$ & Cumulative outflow of MAUs at depot $m \in \mathcal{M}$ at time $t$ in scenario $w$\\
\hline
\end{tabular}%
}
\end{table}

\newpage
\section{Linearization of the proposed formulation}
\label{sec:Linearization}

Constraints (\ref{z_wait}), (\ref{z_transfer}), (\ref{cons_indicator_transfer}), (\ref{cons_up_mv}), (\ref{con:downMv}), (\ref{cons_e}), and (\ref{cons_xh_2}) exhibit a nonlinear nature. We proceed to present their linearized forms in the following discussions. 

First, we introduce variables $\hat{\lambda}^l_{p,k,i}(w)$ to derive the equivalent linear form of constraints (\ref{z_wait}), that is, 
\begin{align}\label{z_wait_linear}
\begin{cases}
\Psi_p^1(w) = \sum\limits_{k\in {\mathcal{K}}_l}(\hat{\lambda}^l_{p,k,i}(w)-u_p\chi^l_{p,k}(w)) & \forall w\in{\mathcal{W}}, p\in\mathcal{P}_w, l=l_0\in\mathcal{L}_p,i=i^o_p,\\
     \hat{\lambda}^l_{p,k,i}(w)\leq M_4\chi^l_{p,k}(w) & \forall w\in{\mathcal{W}}, p\in\mathcal{P}_w, l=l_0\in\mathcal{L}_p, k\in {\mathcal{K}}_l,i=i^o_p,\\
    \hat{\lambda}^l_{p,k,i}(w)\leq d_{k,i}^{l}(w) & \forall w\in{\mathcal{W}}, p\in\mathcal{P}_w, l=l_0\in\mathcal{L}_p, k\in {\mathcal{K}}_l,i=i^o_p,\\
     \hat{\lambda}^l_{p,k,i}(w)\geq d_{k,i}^{l}(w) - M_4(1-\chi^l_{p,k}(w)) & \forall w\in{\mathcal{W}}, p\in\mathcal{P}_w, l=l_0\in\mathcal{L}_p, k\in {\mathcal{K}}_l,i=i^o_p,\\
     \hat{\lambda}^l_{p,k,i}(w)\in [0,\left|\mathcal{T}\right|] & \forall w\in{\mathcal{W}}, p\in\mathcal{P}_w, l=l_0\in\mathcal{L}_p, k\in {\mathcal{K}}_l,i=i^o_p.
    \end{cases}
\end{align}

By introducing variables $\tilde{\psi}_{p, k, k'}^{g}(w) =\zeta_{p,k,k'}^{g}(w)(1-h^g_{k,k'}(w))$, $\dot{\psi}_{p, k, k'}^{g}(w) =\tilde{\psi}_{p, k, k'}^{g}(w)\cdot d_{k',i'}^{l'}$ and  $\hat{\psi}_{p, k, k'}^{g}(w) =\tilde{\psi}_{p, k, k'}^{g}(w)\cdot a_{k,i}^{l}$, constraints (\ref{z_transfer}) can be linearized as below
\begin{align}
\begin{cases}
&\Psi_p^2(w) = \sum\limits_{g\in\mathcal{R}_p}\sum\limits_{k\in\mathcal{K}(i)}\sum\limits_{k'\in\mathcal{K}(i')}(\dot{\psi}_{p, k, k'}^{g}(w)-\hat{\psi}_{p, k, k'}^{g}(w))\\
&\tilde{\psi}_{p, k, k'}^{g}(w)\leq\zeta_{p,k,k'}^{g}(w)\\
&\tilde{\psi}_{p, k, k'}^{g}(w)\leq 1-h^g_{k,k'}(w)\\
&\tilde{\psi}_{p, k, k'}^{g}(w)\geq\zeta_{p,k,k'}^{g}(w)-h^g_{k,k'}(w)\\
    &\dot{\psi}_{p, k, k'}^{g}(w)\leq M_5\tilde{\psi}_{p, k, k'}^{g}(w)\\
    &\dot{\psi}_{p, k, k'}^{g}(w)\leq d_{k',i'}^{l'}\\
    &\dot{\psi}_{p, k, k'}^{g}(w)\geq d_{k',i'}^{l'} - M_5 (1-\tilde{\psi}_{p, k, k'}^{g}(w))\\
     &\hat{\psi}_{p, k, k'}^{g}(w)\leq M_5\tilde{\psi}_{p, k, k'}^{g}(w)\\
    &\hat{\psi}_{p, k, k'}^{g}(w)\leq a_{k,i}^{l}\\
    &\hat{\psi}_{p, k, k'}^{g}(w)\geq a_{k,i}^{l} - M_5 (1-\tilde{\psi}_{p, k, k'}^{g}(w))\\
    &\dot{\psi}_{p, k, k'}^{g}(w), \hat{\psi}_{p, k, k'}^{g}(w) \in [0, \left|\mathcal{T}\right|]
    \end{cases}\forall w\in\mathcal{W},p\in\mathcal{P}_w, g=(l,i,l',i')\in\dot{\mathcal{A}'},k\in\mathcal{K}_l,k'\in\mathcal{K}_l'.\label{z_transfer_linear}
\end{align}

Through defining $\gamma^{g}_{k,k'}=\pi^{g}_{k,k'}-\pi^{g}_{k,k'-1}$, constraints (\ref{cons_indicator_transfer}) can be linearized as follows:
\begin{align}
\begin{cases}
\zeta_{p,k,k'}^{g}(w)\leq\chi_{p,k}^l(w) &\forall w\in{\mathcal{W}} ,p\in\mathcal{P}_w,g=(l,i,l',i')\in\mathcal{R}_p, k\in {\mathcal{K}}_l, k'\in {\mathcal{K}}_{l'}, \\
\gamma^{g}_{k,k'}=\pi^{g}_{k,k'}-\pi^{g}_{k,k'-1} &\forall g=(l,i,l',i')\in\mathcal{R}_p, k\in {\mathcal{K}}_l, k'\in {\mathcal{K}}_{l'} \\
\zeta_{p,k,k'}^{g}(w)\leq\gamma^{g}_{k,k'}&\forall w\in{\mathcal{W}} ,p\in\mathcal{P}_w,g=(l,i,l',i')\in\mathcal{R}_p, k\in {\mathcal{K}}_l, k'\in {\mathcal{K}}_{l'},\\
\zeta_{p,k,k'}^{g}(w)\geq \chi_{p,k}^l(w)+\gamma^{g}_{k,k'}-1 
&\forall w\in{\mathcal{W}} ,p\in\mathcal{P}_w, g=(l,i,l',i')\in\mathcal{R}_p, k\in {\mathcal{K}}_l, k'\in {\mathcal{K}}_{l'}.\end{cases}\label{cons_indicator_transfer_linear}
\end{align}

Then, we introduce the integer variable $\mu_{k, i ,t}^l(w) = z_{k,i,t}^l x_{k,i}^l(w)$. With this modification, constraints (\ref{cons_up_mv}) and (\ref{con:downMv}) can be equivalently transformed into linearized forms as below:

\begin{align}\label{cons_vs_linear}
&\left\{ \begin{aligned}
& AV_{m,t}(w)=\sum_{l\in LL_s}\sum_{k\in\mathcal{K}_l}(x_{k,\left|\mathcal{S}_l\right|}^l(w)  -\mu_{k,\left|\mathcal{S}_l\right|,t-t^{pre}_i}^l(w)) && \forall m\in \mathcal{M}, t \in\mathcal{T}, w\in\mathcal{W},\\
& DV_{m,t}(w)= \sum_{l\in FL_m}\sum_{k\in\mathcal{K}_l}(x_{k,1}^l(w)-\mu_{k,1,t}^l(w)) && \forall m\in\mathcal{M}, t\in\mathcal{T}, w\in\mathcal{W}, \\
&\mu^l_{k, i ,t}(w)\leq x_{k,i}^l(w) && \forall l\in{\mathcal{L}}, k \in {\mathcal{K}}_l, i \in {\mathcal{S}}_l,w\in\mathcal{W}, \\
& \mu^l_{k, i ,t}(w)\geq x_{k,i}^l(w) - M_6 (1-z_{k,i,t}^l) && \forall l\in{\mathcal{L}}, k \in {\mathcal{K}}_l, i \in {\mathcal{S}}_l,w\in\mathcal{W}, \\
&\mu^l_{k, i ,t}(w)\leq M_6 z_{k,i,t}^l && \forall l\in{\mathcal{L}}, k \in {\mathcal{K}}_l, i \in {\mathcal{S}}_l,w\in\mathcal{W}, \\
& \mu^l_{k, i ,t}(w)\in [0, \left|\mathcal{Q}\right|] && \forall l\in{\mathcal{L}}, k \in {\mathcal{K}}_l, i \in {\mathcal{S}}_l,w\in\mathcal{W}.
\end{aligned} \right.
\end{align}

By introducing binary variables $\dot{e}^{g}_{k,k'}$ and $\hat{e}^{g}_{k,k'}$, constraints (\ref{cons_e}) can be linearized as follows:
\begin{align}
\begin{cases}
&M_7(\dot{e}^{g}_{k,k'}-1)\leq d^{l'}_{k',i'}-a^{l}_{k,i}-\underline{\theta}_{g}\leq M_7\dot{e}^{g}_{k,k'}- \epsilon \\
&M_7(\hat{e}^{g}_{k,k'}-1)\leq\overline{\theta}_{g} - d^{l'}_{k',i'}+a^{l}_{k,i}\leq M_7\hat{e}^{g}_{k,k'}- \epsilon\\
&e^{g}_{k,k'}\leq \dot{e}^{g}_{k,k'}\\
&e^{g}_{k,k'}\leq \hat{e}^{g}_{k,k'}\\
&e^{g}_{k,k'}\geq \dot{e}^{g}_{k,k'} + \hat{e}^{g}_{k,k'} - 1\\
&\dot{e}^{g}_{k,k'},\hat{e}^{g}_{k,k'}\in\{0,1\}
\end{cases}
&&\forall g=(l,i,l',i')\in\dot{\mathcal{A}}, k\in\mathcal{K}_{l}, k'\in\mathcal{K}_{l'}.\label{cons_e_linear}
\end{align}

Finally, we introduce four kinds of auxiliary variables, i.e. $h^1_{g,k,k'}(w)$, $h^2_{g,k,k'}(w)$, $h^3_{g,k,k'}(w)$ and $h^4_{g,k,k'}(w)$. The linearization of constraints (\ref{cons_xh_2}) can be achieved as follows:
\begin{align}
\begin{cases}
&x^l_{k,i}(w) = x^l_{k,i-1}(w) +\sum\limits_{l'}\sum\limits_{k'\in\mathcal{K}_{l'}(i')} (h^1_{g,k,k'}(w) \\ 
& \quad -h^2_{g,k,k'}(w) - h^3_{g,k,k'}(w)+h^4_{g,k,k'}(w))\\ 
&h^1_{g,k,k'}(w)\leq M_8 h^{g'}_{k',k}(w)\\
&h^1_{g,k,k'}(w)\leq x^{l'}_{k',i'-1}(w)\\
&h^1_{g,k,k'}(w)\geq x^{l'}_{k',i'-1}(w) + M_8(1-h^{g'}_{k',k}(w))\\
&h^2_{g,k,k'}(w)\leq M_8 h^{g'}_{k',k}(w)\\
&h^2_{g,k,k'}(w)\leq x^{l'}_{k',i'}(w)\\
&h^2_{g,k,k'}(w)\geq x^{l'}_{k',i'}(w) + M_8(1-h^{g'}_{k',k}(w))\\
&h^3_{g,k,k'}(w)\leq M_8 h^{g}_{k,k'}(w)\\
&h^3_{g,k,k'}(w)\leq x^{l'}_{k',i'}(w)\\
&h^3_{g,k,k'}(w)\geq x^{l'}_{k',i'}(w) + M_8(1-h^{g}_{k,k'}(w))\\
&h^4_{g,k,k'}(w)\leq M_8 h^{g}_{k,k'}(w)\\
&h^4_{g,k,k'}(w)\leq x^{l'}_{k',i'-1}(w)\\
&h^4_{g,k,k'}(w)\geq x^{l'}_{k',i'-1}(w) + M_8(1-h^{g}_{k,k'}(w))\\
&h^1_{g,k,k'}(w),h^2_{g,k,k'}(w),h^3_{g,k,k'}(w),h^4_{g,k,k'}(w)\in [0, \left|\mathcal{Q}\right|]
\end{cases}\quad \begin{array}{l}
\forall w\in\mathcal{W}, l\in{\mathcal{L}}, k\in\mathcal{K}_l,i \in \mathcal{R}_{l}, \\
g=(l,i,l',i')\in\dot{\mathcal{A}}, g'=(l',i',l,i)\in\dot{\mathcal{A}}.
\end{array}\label{cons_hh_linear}
\end{align}

\newpage
\section{Proofs of Propositions 1-4 and Lemma 1}
\label{sec:proof}

\textbf{Proof of Proposition \ref{propo_domin_reduction}.}
 Note that stop-skip and speed control strategies are not considered in this study. Hence, combining the constraints related to shifting time (\ref{cons_shift}) and the dwell time (\ref{cons_dwell}), we can thus obtain the following equation

\begin{align}\label{ve_departure_minmax}
\begin{cases}
    &\underline{\varsigma}^l_{k,i} = o^l_{k,i} + \alpha^{min}_{l,k} + \sum\limits_{j\leq i}\beta^{min}_{l,k,i}+ \sum\limits_{j\leq i}r^{l}_{k,j-1}\\
&\overline{\varsigma}^l_{k,i} = o^l_{k,i} + \alpha^{max}_{l,k} + \sum\limits_{j\leq i}\beta^{max}_{l,k,i}+ \sum\limits_{j\leq i}r^{l}_{k,j-1} 
\end{cases}, \qquad  \forall l\in{\mathcal{L}}, k \in {\mathcal{K}}_l, i,j \in {\mathcal{S}}_l.
\end{align}
Also note that prior to optimization, the value of the right-hand side of the equation (\ref{ve_departure_minmax}) is determined. For clarity, we assume that $\lceil{\frac{\underline{\varsigma}^l_{k,i}}{\Delta}}\rceil=\lfloor{\frac{\underline{\varsigma}^l_{k,i}}{\Delta}}\rfloor$ and $\lceil{\frac{\overline{\varsigma}^l_{k,i}}{\Delta}}\rceil=\lfloor{\frac{\overline{\varsigma}^l_{k,i}}{\Delta}}\rfloor$, which is easy to achieve during the data processing phase. Then, we have $\underline{\varsigma}^l_{k,i}\leq d^l_{k,i}\leq\overline{\varsigma}^l_{k,i}$. Next, we prove the two cases in Equation (\ref{cons_ve_timetable}) separately.

(i) Since the variable $z^l_{k,i,t}$ possesses the non-decreasing property as specified in constraints (\ref{binary_z}), we can deduce the subsequent inequality:
\begin{align}
\sum_{t\in\mathcal{T}}z^l_{k,i,t}=\sum_{t\in\mathcal{T}({\left|\mathcal{T}\right|})}z^l_{k,i,t}\geq\sum_{t\in\mathcal{T}({\left|\mathcal{T}\right|-1})}z^l_{k,i,t}\geq ... \geq\sum_{t\in\mathcal{T}({\frac{\underline{\varsigma}^l_{k,i}}{\Delta}-1})}z^l_{k,i,t},
\end{align}
where $\mathcal{T}(t)=\{1,2,...,t\}$. 

According to the relationship between $z^l_{k,i,t}$ and $d^l_{k,i}$ as shown in constraints (\ref{cons_couple}), we have
\begin{align}\label{eq_propostion_1_a}
    d^l_{k,i}=\Delta(1+\sum_{t\in\mathcal{T}}z^l_{k,i,t})\geq\Delta(1+\sum_{t\in\mathcal{T}({\frac{\underline{\varsigma}^l_{k,i}}{\Delta}-1})}z^l_{k,i,t}) \geq \underline{\varsigma}^l_{k,i}.
\end{align}
Because $z^l_{k,i,t}\leq 1$ holds, we can derive the following inequality:
\begin{align}\label{eq_propostion_1_b}
\sum_{t\in\mathcal{T}({\frac{\underline{\varsigma}^l_{k,i}}{\Delta}-1})}z^l_{k,i,t}\leq\frac{\underline{\varsigma}^l_{k,i}}{\Delta}-1.
\end{align}
Clearly, inequalities (\ref{eq_propostion_1_a}) and (\ref{eq_propostion_1_b}) hold simultaneously if and only if $z^l_{k,i,t}=1, \forall 1\leq t\leq \frac{\underline{\varsigma}^l_{k,i}}{\Delta}-1$.

(ii) From $d^l_{k,i}\leq\overline{\varsigma}^l_{k,i}$, we have
\begin{align}\label{eq_propostion_2_a}
    \sum_{t\in\mathcal{T}}z^l_{k,i,t}=\sum\limits_{t\in\mathcal{T}\setminus\mathcal{T}({\frac{\overline{\varsigma}^l_{k,i}}{\Delta}-1})}z^l_{k,i,t} + \sum\limits_{t\in\mathcal{T}({\frac{\overline{\varsigma}^l_{k,i}}{\Delta}-1})}z^l_{k,i,t} \leq \frac{\overline{\varsigma}^l_{k,i}}{\Delta}-1.
\end{align}
Similar to the derivation in (i), the inequality $\sum\limits_{t\in\mathcal{T}({\frac{\overline{\varsigma}^l_{k,i}}{\Delta}-1})}z^l_{k,i,t}\leq \frac{\overline{\varsigma}^l_{k,i}}{\Delta}-1$ must holds. From inequality (\ref{eq_propostion_2_a}), we can find that $\sum\limits_{t\in\mathcal{T}\setminus\mathcal{T}({\frac{\overline{\varsigma}^l_{k,i}}{\Delta}-1})}z^l_{k,i,t}\leq 0$. In other words, there exists a feasible solution for the TT-VS-DCA problem if and only if $z^l_{k,i,t}=0, \forall\frac{\overline{\varsigma}^l_{k,i}}{\Delta} \leq t\leq \left|\mathcal{T}\right|$.
\Halmos \endproof

While the path-based representation of passenger demand is convenient to denote passenger movements, this modeling approach inevitably introduces a substantial number of binary variables in the formulation. This in turn considerably exacerbates the difficulty of solving the model. To address this challenge, we propose the following three types of valid inequalities to narrow the search space of binary variables $\xi^{l}_{p,k}(w)$, $\pi^{g}_{k,k'}$, and $e^{g}_{k,k'}$, thus further accelerating the solution process.

\textbf{Proof of Proposition \ref{propo_first_board}.}
The intuitive interpretation of proposition \ref{propo_first_board} is that a passenger cannot board a vehicle whose latest arrival time exceeds his or her arrival time. That is, $\xi^{l}_{p,k}(w)-\xi^{l}_{p,k-1}(w)=1\Rightarrow \overline{\varsigma}^l_{k,i}\geq u_p$  always holds. By introducing the big-$M$(here taking the value $\left|\mathcal{T}\right|$), we can derive constraints (\ref{cons_ve_boarding}).
\Halmos \endproof

\textbf{Proof of Proposition \ref{propo_transfer}.}
According to the definition of the connected trip, the pair of trips $(k, k')$ is connectable only when there is enough transfer time for trip $k$ and $k^{'}$. Therefore, we have $\pi^{g}_{k,k'}-\pi^{g}_{k,k'-1}=1\Rightarrow \max\limits_{\mathbf{z}}\{d^{l'}_{k',i'}-a^{l}_{k,i}-\underline{\theta}_g\}\geq 0$. From constraints (\ref{cons_buffer}) and (\ref{cons_dwell}), we then find that $\max\limits_{\mathbf{z}}\{d^{l'}_{k',i'}-a^{l}_{k,i}-\underline{\theta}_g\}=\overline{\varsigma}^{l'}_{k',i'}-\min\limits_{\mathbf{z}}(d^l_{k,i}-\beta^l_{k,i})-\underline{\theta}_g=\overline{\varsigma}^{l'}_{k',i'}-\underline{\varsigma}^{l}_{k,i} + \beta^{min}_{l,k,i}-\beta^l_{k,i}$. Similar to proposition \ref{propo_first_board}, inequality (\ref{cons_ve_transfer}) holds.
\Halmos \endproof

\textbf{Proof of Proposition \ref{propo_coupling}.}
Based on constraints (\ref{cons_e}), we can easily derive that $e^{g}_{k,k'}=1\Rightarrow (\min\limits_{\mathbf{z}}\{d^{l'}_{k',i'}-a^{l}_{k,i}\} \leq \overline{\theta}_g\wedge\max\limits_{\mathbf{z}}\{d^{l'}_{k',i'}-a^{l}_{k,i}\} \geq \underline{\theta}_g)$. By conducting a similar derivation as proposition \ref{propo_transfer}, we have inequality (\ref{cons_ve_couple}).
\Halmos \endproof

\textbf{Proof of Lemma \ref{lemma_lower_bound_1}.} In this study, we do not consider the possibility of canceling trips. Consequently, each scheduled trip requires at least one MAU to execute. Suppose that each scenario contains a finite number of passengers and only a single MAU is needed for the MAV assigned to each trip, allowing all passengers boarding this trip, we have: $F(\mathbf{z},\bm{\chi},\bm{\zeta},\mathbf{x})\geq L= \sum\limits_{l\in\mathcal{L}}\sum\limits_{i\in\mathcal{S}_l}\sum\limits_{k\in\mathcal{K}_l}\varphi_2\vartheta_{1,i}.\nonumber$
\Halmos 
\endproof

\section{The integer L-shaped method for the TT-VS-DCA model}
\label{sec:algorithmflow}
Algorithm \ref{table:algorithm} presents the integer L-shaped method for the TT-VS-DCA problem.

\begin{algorithm}[H]
\caption{The integer L-shaped method for the TT-VS-DCA model}\label{table:algorithm}
\SetKwInOut{Input}{Input}
\SetKwInOut{Output}{Output}

\Input{The lower bound $L$ according to Lemma \ref{lemma_lower_bound_1} of the subproblem, root node of the relaxed current problem.}
\Output{The optimal solution of the TT-VS-DCA problem.}
\BlankLine
Create an empty list of nodes and put the root node into that list\;
Set the minimum gap $\epsilon$, upper bound $\mathrm{UB}=+\infty$ and lower bound $\mathrm{LB}=L$\;
Let $\mathrm{Gap} = (\mathrm{UB} - \mathrm{LB}) / \mathrm{UB}$\;
\While{$\mathrm{Gap} \leq \mathrm{Gap}^{min}$}{
    \If{the node list is empty}{
        \textbf{break}\;
    }
    \Else{
        Select a pendent node from the node list\;
        Solve the CP corresponding to the selected node\;
        \If{the CP is infeasible}{
            Fathom the current node and go to line 5\;
        }
        \Else{
            Let $(\mathbf{z}^*, \boldsymbol{\chi}^*,\boldsymbol{\zeta}^*, \eta^*)$ be the optimal solution of the CP\;
            \If{$\Psi_{wait}^1(\mathbf{z}^*, \boldsymbol{\chi}^*,\boldsymbol{\zeta}^*) + \eta \geq \mathrm{UB}$}{
                Fathom the current node and go to line 3\;
            }
            \Else{
                \If{some variables in $\mathbf{z}$ are non-binary}{
                    Select and branch a non-binary variable\;
                    Delete the current node and append the new two branch nodes to the pendent node list\;
                    Then go to line 5\;
                }
                \Else{
                    Compute the linear relaxation of $F(\mathbf{z}^*, \boldsymbol{\chi}^*,\boldsymbol{\zeta}^*)$ and get the corresponding optimal objective value $\mathrm{F_{LP}}$\;
                    \If{The linear relaxation of $F(\mathbf{z}^*, \boldsymbol{\chi}^*,\boldsymbol{\zeta}^*)$ is infeasible}{
                        Add the feasibility cut (\ref{feasibility_cut}) to CP and go to line 8\;
                    }
                    \Else{
                        \If{$\mathrm{F_{LP}} > \eta^*$}{
                            Add the subgradient cut (\ref{copt_cut}) to CP and go to line 8\;
                        }
                        \Else{
                            Call GUROBI to solve the subproblem $F(\mathbf{z}^*, \boldsymbol{\chi}^*,\boldsymbol{\zeta}^*)$ to optimal as $F_0$\;
                            \If{$F(\mathbf{z}^*, \boldsymbol{\chi}^*,\boldsymbol{\zeta}^*)$ is infeasible}{
                                Add the feasibility cut (\ref{feasibility_cut}) to CP and go to line 8\;
                            }
                            \Else{
                                \If{$\Psi_{wait}^1(\mathbf{z}^*, \boldsymbol{\chi}^*, \boldsymbol{\zeta}^*) + F_0 < \mathrm{UB}$}{
                                    Update $\mathrm{UB} = \Psi_{wait}^1(\mathbf{z}^*, \boldsymbol{\chi}^*, \boldsymbol{\zeta}^*) + F_0$\;
                                }
                                \If{$F_0 > \eta^*$}{
                                    Add the optimality cut (\ref{opt_cut}) to CP and go to line 8\;
                                }
                            }
                        }
                    }
                }
            }
        }
    }
}
\end{algorithm}

\section{The procedure of the overall algorithm}
\label{sec:flowChat}
\begin{figure}[h]
    \centering
    \includegraphics[height=12cm]{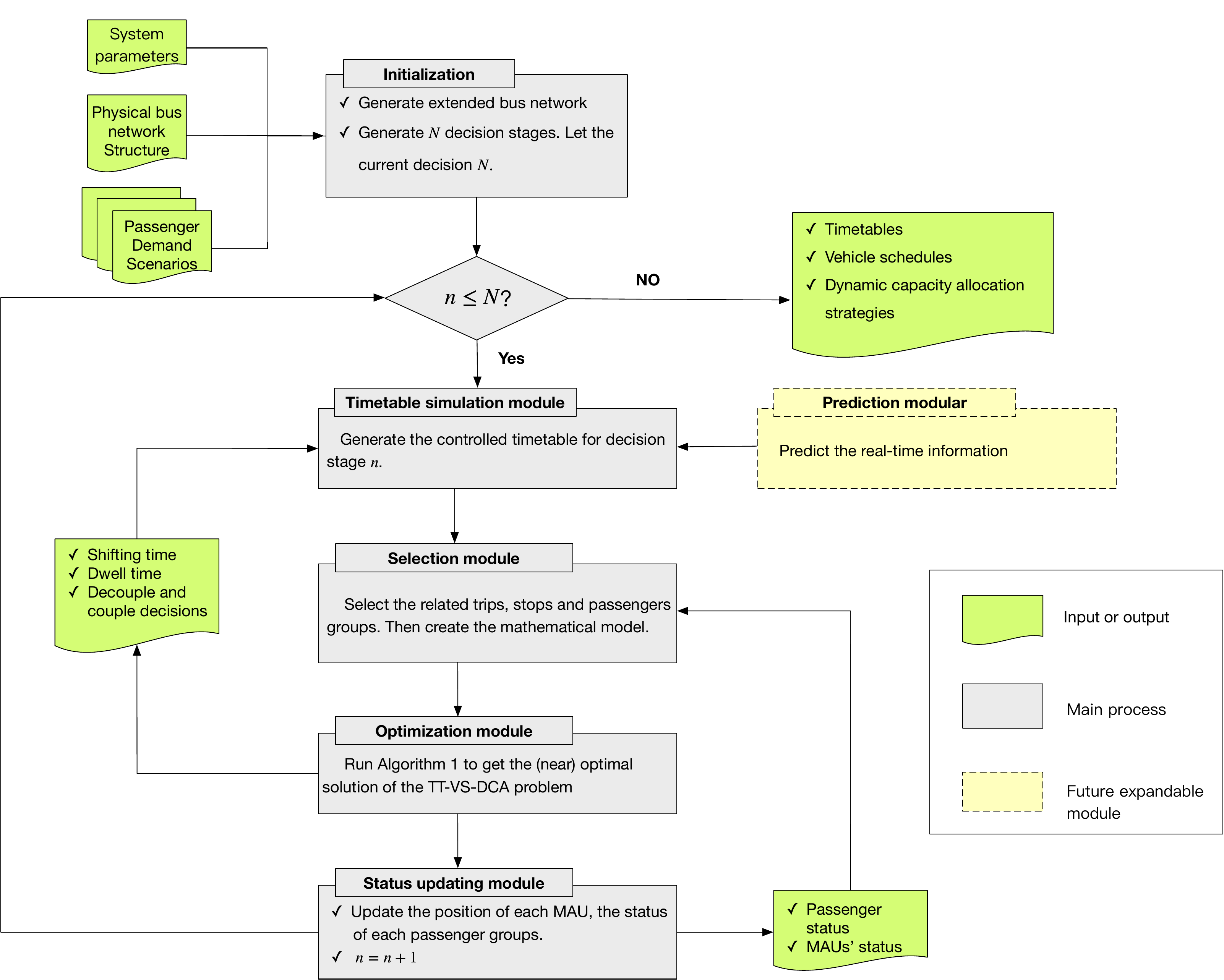}
    \caption{Rolling horizon optimization framework for the TT-VS-DCA problem.}
    \label{fig:RH_new}
\end{figure}

Figure \ref{fig:RH_new} shows the proposed rolling horizon optimization framework for the TT-VS-DCA problem.

\newpage
\section{Overflow of the learning-based real-time decision-making framework}
\label{sec:Overflow}

\begin{figure}[h]
\centering
\includegraphics[height=8cm]{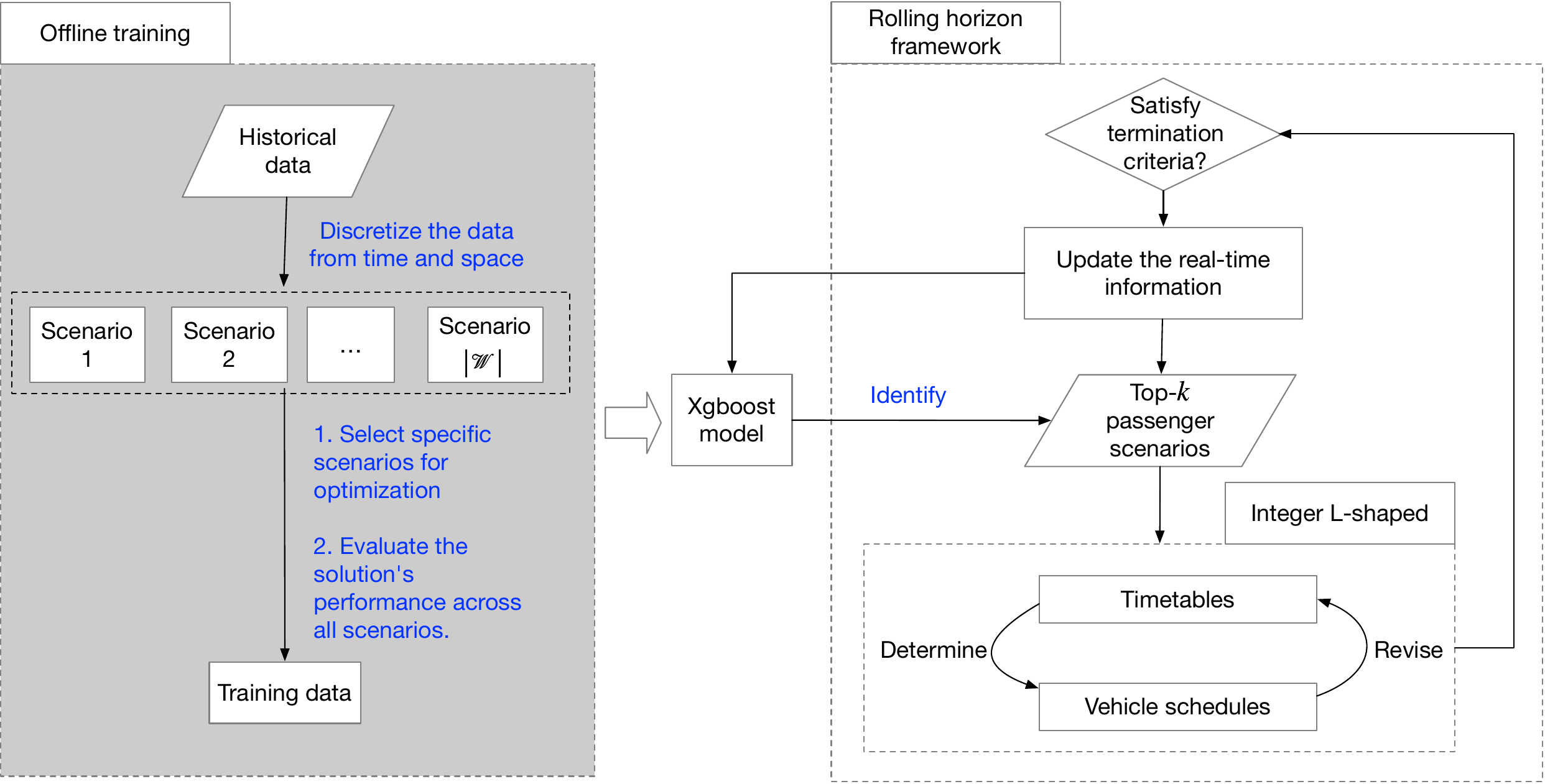}
\caption{Real-time decision framework of the TT-VS-DCA model.}\label{fig:Real_time}
\end{figure}

Figure \ref{fig:Real_time} illustrates our proposed learning-based real-time decision-making framework.

\newpage
\section{Experimental instances employed in the study}
\label{sec:instance}

\begin{figure}[h]
     \centering
     \begin{subfigure}{0.52\textwidth}
         \centering
         \includegraphics[width=0.9\textwidth]{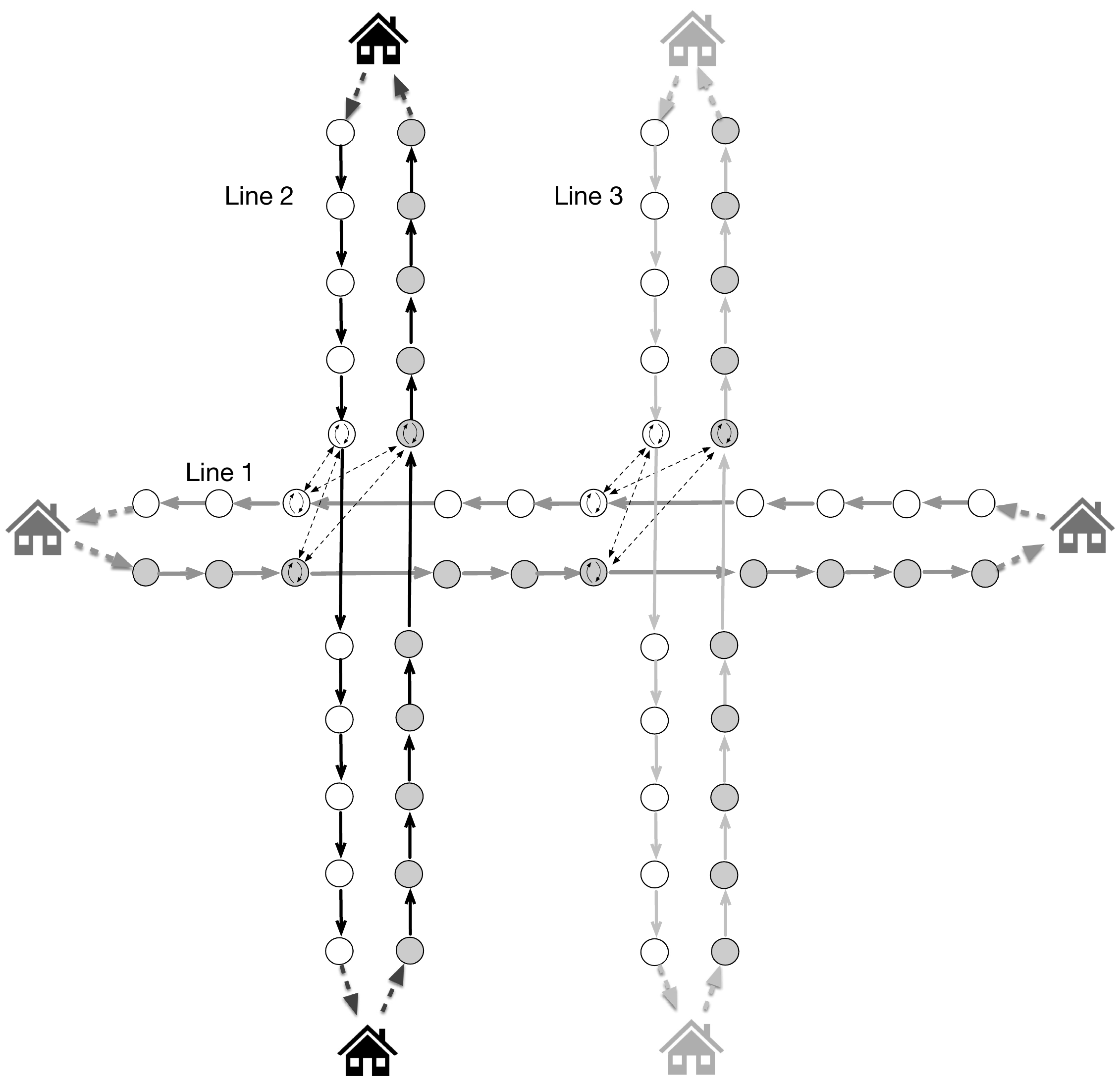}
         \caption{The virtual network}
     \end{subfigure}
     \begin{subfigure}{0.47\textwidth}\label{experimentNetwork_b}
         \centering
         \includegraphics[width=1\textwidth]{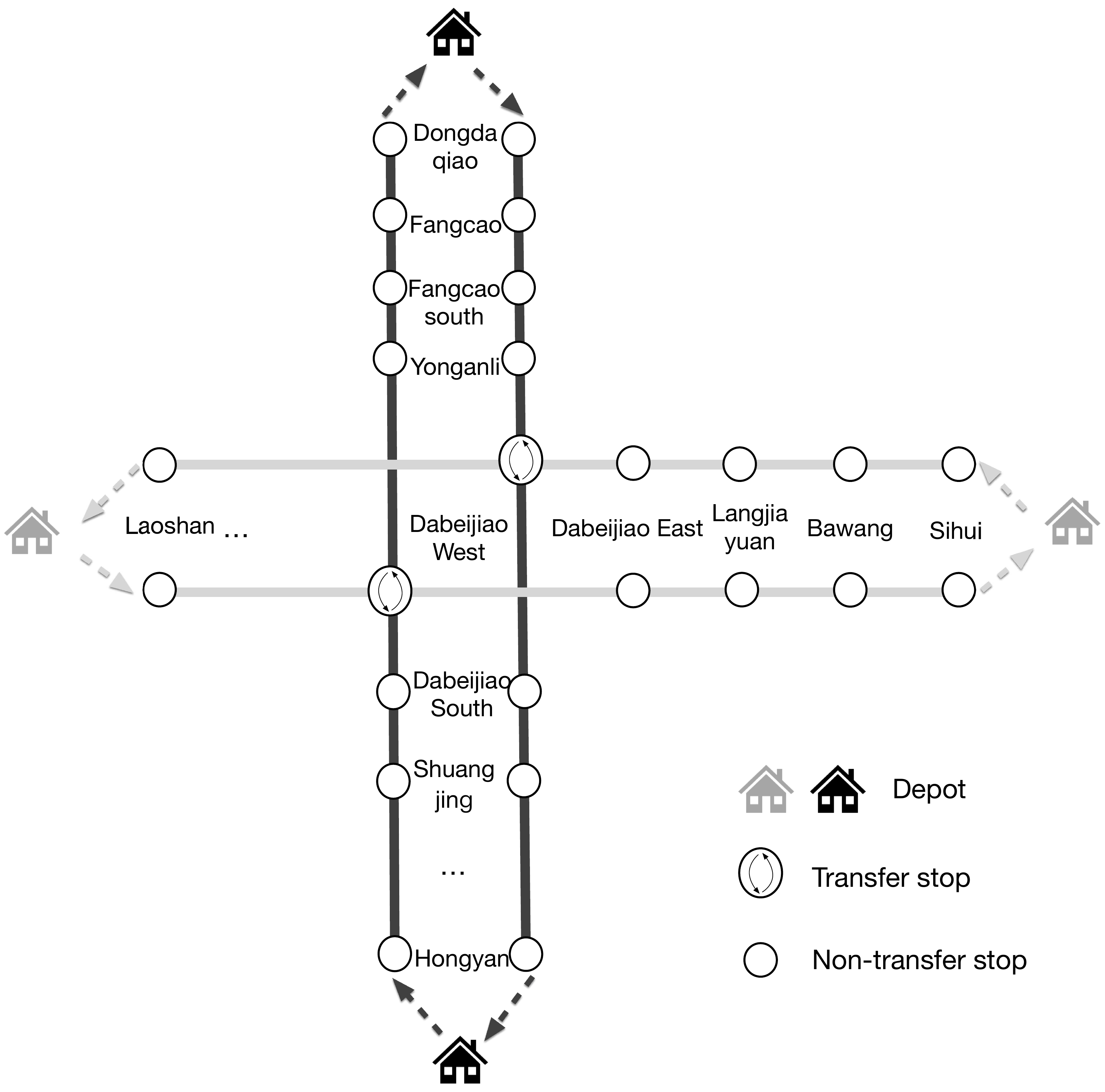}
         \caption{Beijing bus subnetwork}
     \end{subfigure}
        \caption{An illustration of the networks used in numerical experiments.}
        \label{fig:experimentNetwork}
\end{figure}

In this study, we generate various instances for conducting different experiments. Specifically, considering the network depicted in Figure \ref{fig:experimentNetwork}(a), we generate twelve instances, ranging from $instance\_g\_1$ to $instance\_g\_12$. Additionally, from the sub-network shown in Figure \ref{fig:experimentNetwork}(b), which includes only partial stops for each line, we generate six instances, named $instance\_rp\_1$ to $instance\_rp\_6$. These instances arise by progressively increasing both the number of trips and their study time horizon. Furthermore, we generate the other six instances, i.e., from $instance\_r\_1$ to $instance\_r\_6$, based on the network shown in Figure \ref{fig:experimentNetwork}(b), enabling a thorough exploration of the effectiveness of our approaches in real-world problems.

\begin{table}[h]
\centering
\caption{Characteristics of the instances used in the numerical experiments.}
\renewcommand{\arraystretch}{1.2}
\begin{tabular}{crrrrrcrrrrr}
\hline
Instance  & $\left|\mathcal{L}\right|$ & $\left|\mathcal{S}\right|$  & $\left|\mathcal{K}\right|$  & $\left|\mathcal{T}\right|$   & $\left|\mathcal{W}\right|$ & Instance   & $\left|\mathcal{L}\right|$ & $\left|\mathcal{S}\right|$  & $\left|\mathcal{K}\right|$  & $\left|\mathcal{T}\right|$   & $\left|\mathcal{W}\right|$  \\ \hline
instance\_g\_1 & 3 & 48 & 18 & 60  & 8 & instance\_g\_7  & 3 & 48 & 18 & 60  & 12 \\
instance\_g\_2 & 3 & 48 & 24 & 80  & 8 & instance\_g\_8  & 3 & 48 & 24 & 80  & 12 \\
instance\_g\_3 & 3 & 48 & 30 & 100 & 8 & instance\_g\_9  & 3 & 48 & 30 & 100 & 12 \\
instance\_g\_4 & 3 & 48 & 36 & 120 & 8 & instance\_g\_10 & 3 & 48 & 36 & 120 & 12 \\
instance\_g\_5 & 3 & 48 & 42 & 140 & 8 & instance\_g\_11 & 3 & 48 & 42 & 140 & 12 \\
instance\_g\_6 & 3 & 48 & 48 & 160 & 8 & instance\_g\_12 & 3 & 48 & 48 & 160 & 12 \\ \addlinespace[0.1em]
instance\_rp\_1 & 2 & 32 & 8 & 70  & 8 & instance\_rp\_2  & 2 & 32 & 12 & 80  & 8 \\
instance\_rp\_3 & 2 & 32 & 16 & 90  & 8 & instance\_rp\_4  & 2 & 32 & 20 & 100  & 8 \\
instance\_rp\_5 & 2 & 32 & 24 & 110  & 8 & instance\_rp\_6  & 2 & 32 & 28 & 120  & 8 \\ \addlinespace[0.1em]
instance\_r\_1 & 2 & 89 & 8 & 120 & 8 & instance\_r\_4 & 2 & 89 & 24 & 150 & 8 \\
instance\_r\_2 & 2 & 89 & 12 & 130 & 8 & instance\_r\_5 & 2 & 89 & 28 &160  & 8 \\
instance\_r\_3 & 2 & 89 & 16 & 140 & 8 & instance\_r\_6 & 2 & 89 & 50 & 240 & 8 \\ 
\hline
\end{tabular}
\label{tab:instance_all}
\end{table}
\end{appendices}
	\end{document}